\documentclass[11pt]{amsart}

\usepackage{amssymb,amsfonts,amsmath}
\usepackage{amsthm, stmaryrd}
\usepackage{graphicx}
\usepackage[all,arc]{xy}
\usepackage{hyperref}
\hypersetup{colorlinks,allcolors=violet}
\usepackage{enumerate}
\usepackage{bbm}
\usepackage{dsfont}
\usepackage{mathrsfs}
\usepackage{tikz-cd}
\usetikzlibrary{shapes}
\usepackage{import}
\usepackage{float}
\restylefloat{table}
\usepackage{multirow}
\usepackage{pdflscape}
\usepackage{listofitems}
\usepackage{standalone}

\usepackage[margin=1in]{geometry}
\usepackage{mathtools}
\usetikzlibrary{decorations.pathmorphing}
\usepackage{mathrsfs}
\newtheorem{thm}{Theorem}[section]
\newtheorem*{thm*}{Theorem}

\newtheorem{innercustomthm}{Theorem}
\newenvironment{customthm}[1]
  {\renewcommand\theinnercustomthm{#1}\innercustomthm}
  {\endinnercustomthm}

\newtheorem{cor}[thm]{Corollary}
\newtheorem{prop}[thm]{Proposition}
\newtheorem{lem}[thm]{Lemma}

\theoremstyle{definition}
\newtheorem{defn}[thm]{Definition}

\newtheorem{exmp}[thm]{Example}

\newtheorem{rem}[thm]{Remark}
\newtheorem*{thm1.2}{\textrm{Theorem 1.2}}

\theoremstyle{remark}
\newcommand{\Mbar}{\overline{\mathcal{M}}}
\newcommand{\M}{\mathcal{M}}

\newtheorem{sch}[thm]{Scholium}

\newcommand{\sfch}{\mathsf{ch}}
\newcommand{\sfe}{\mathsf{e}}

\newcommand{\Z}{\mathbb{Z}}

\newcommand{\QQ}{\mathbb{Q}}

\renewcommand{\P}{\mathbb{P}}

\newcommand{\Aut}{\operatorname{Aut}}

\newcommand{\val}{\operatorname{val}}

\newcommand{\ch}{\operatorname{ch}}

\def\C{\mathbb{C}}

\def\E{\mathbb{E}}

\def\P{\mathbb{P}}

\def\Z{\mathbb{Z}}

\def\M{\mathcal{M}}

\newcommand{\Cat}{\mathsf{Cat}}
\newcommand{\Dih}{\mathsf{Dih}}

\newcommand{\fS}{\mathbb{S}}

\newcommand{\pp}{\mathfrak{p}}
\newcommand{\qq}{\mathfrak{q}}

\newcommand{\Ind}{\operatorname{Ind}}
\newcommand{\Res}{\operatorname{Res}}
\newcommand{\Spec}{\operatorname{Spec}}

\newcommand{\Log}{\operatorname{Log}}

\newcommand\cycle[2][\,]{%
  \readlist\thecycle{#2}%
  (\foreachitem\i\in\thecycle{\ifnum\icnt=1\else#1\fi\i})%
}

\makeatletter
\let\c@equation\c@thm
\makeatother
\numberwithin{equation}{section}

\graphicspath{ {images/} }
\title{The $S_n$-equivariant Euler characteristic of $\Mbar_{1, n}(\P^r, d)$}

\author[S. Kannan]{Siddarth Kannan}\address{Department of Mathematics, Massachusetts Institute of Technology}
\email{\url{spkannan@mit.edu}}
\author[T. Song]{Terry Dekun Song}\address{Department of Pure Mathematics and Mathematical Statistics, University of Cambridge, Cambridge, CB3 0WA}\email{\url{ds2016@cam.ac.uk}}

\begin{document}

\begin{abstract}
We compute the $S_n$-equivariant topological Euler characteristic of the Kontsevich moduli space $\Mbar_{1, n}(\P^r, d)$. Letting $\Mbar_{1, n}^{\mathrm{nrt}}(\P^r, d) \subset \Mbar_{1, n}(\P^r, d)$ denote the subspace of maps from curves without rational tails, we solve for the motive of $\Mbar_{1, n}(\P^r, d)$ in terms of $\Mbar_{1, n}^{\mathrm{nrt}}(\P^r, d)$ and plethysm with a genus-zero contribution determined by Getzler and Pandharipande \cite{GetzlerPandharipande}. Fixing a generic $\C^\star$-action on $\P^r$, we derive a closed formula for the Euler characteristic of $\Mbar_{1, n}^{\mathrm{nrt}}(\P^r, d)^{\C^\star}$ as an $S_n$-equivariant virtual mixed Hodge structure, which leads to our main formula for the Euler characteristic of $\Mbar_{1,n}(\P^r, d)$. Our approach connects the geometry of torus actions on Kontsevich moduli spaces \cite{graberpand} with symmetric functions in Coxeter types $A$ and $B$ \cite{Macdonald}, as well as the enumeration of graph colourings with prescribed symmetry \cite{Hanlon}.
\end{abstract}

\maketitle\thispagestyle{empty}

\section{Introduction}
Let $\Mbar_{g, n}(\P^r, d)$ denote the Kontsevich moduli space of stable maps of degree $d$ from $n$-pointed curves of genus $g$ to $\mathbb{P}^r = \mathbb{P}^r_{\C}$. It is a fundamental object in enumerative geometry and in the study of algebraic curves. The fibres of the map $\Mbar_{g}(\P^r, d) \to \Mbar_{g}$ are studied in Brill--Noether theory. While the moduli space is smooth and irreducible when $g = 0$, it is usually reducible, non-equidimensional, and highly singular when $g > 0$. These pathologies make it difficult to study the geometry of $\Mbar_{g, n}(\P^r, d)$ in positive genus. From the enumerative perspective, these issues have been circumvented via torus localisation by Graber and Pandharipande \cite{graberpand}, following earlier work of Kontsevich \cite{KontsevichTorusActions}. They express the Gromov--Witten invariants of $\P^r$ in terms of computable integrals on the space $\Mbar_{g, n}(\P^r, d)^{\C^\star}$ of torus-fixed stable maps. The action of $\C^\star$ is induced by a generic one on $\P^r$, with $r + 1$ isolated fixed points.

In this paper, we enhance this perspective, shifting from integrals to $S_n$-equivariant motivic invariants of the moduli space $\Mbar_{g, n}(\P^r, d)^{\C^\star}$, where the symmetric group $S_n$ permutes the marked points. This space admits a combinatorial stratification, where each stratum is isomorphic to a product of moduli spaces of stable curves, modulo the action of the automorphism group of a finite graph. In the case $g = 1$, we solve the combinatorics of the strata, at the motivic level, in terms of plethysm of symmetric functions. Our main application is a formula for the $S_n$-equivariant topological Euler characteristic
\begin{equation}\label{eq:frobchardefn}
\chi^{S_n}\left(\Mbar_{1, n}(\P^r, d)\right) := \sum_{i \geq 0} (-1)^i \ch_n(H^i(\Mbar_{1, n}(\P^r, d); \QQ)) \in \Lambda,
 \end{equation}
 where
 \[ \Lambda := \lim_{\longleftarrow} \QQ[[x_1, \ldots, x_n]]^{S_n} \]
 is the completed ring of symmetric functions, and $\ch_n(V) \in \Lambda$ denotes the Frobenius characteristic of an $S_n$-representation $V$. The invariant $\chi^{S_n}\left(\Mbar_{1, n}(\P^r, d)\right)$ is the enrichment of the ordinary topological Euler characteristic to a virtual character of the symmetric group.
 
For a fixed genus $g$, let
\[ \mathsf{b}_{g, r} := \sum_{n \geq 0} \sum_{d\geq 0}\chi^{S_n}(\Mbar_{g, n}(\P^r, d)) q^d \in \Lambda[[q]] \]
be the generating function for the $S_n$-equivariant Euler characteristics. When $g = 0$, a recursive formula for $\mathsf{b}_{0, r}$ has been given by Getzler and Pandharipande \cite{GetzlerPandharipande}. As noted in \S\ref{introsec:relatedwork} below, their result is actually strong enough to determine the $S_n$-equivariant Hodge decomposition of the moduli space. In positive genus, the pathologies of the space make it a subtle problem to calculate the coarser invariant $\mathsf{b}_{g,r}$.

We determine $\mathsf{b}_{1, r}$ in terms of $\mathsf{b}_{0, r}$, as well as the two terms $\mathsf{a}_0, \mathsf{a}_1 \in \Lambda$, where
\begin{equation}\label{eq:sfa} 
\mathsf{a}_g := \sum_{n > 2 - 2g} \chi^{S_n}(\M_{g, n}).
\end{equation}  Closed formulas for $\mathsf{a}_g$ when $g \leq 1$ have been given by Getzler: see \cite[Theorem 5.7]{GetzlerGenusZero} for $\mathsf{a}_0$ and \cite[Corollary 5.4]{GetzlerMHM} for $\mathsf{a}_1$. The terms $\mathsf{a}_0$ and $\mathsf{a}_1$ are combined in our formula to determine $\mathsf{b}_{1,r}^{\mathrm{nrt}}$, which is defined by $$\mathsf{b}_{1,r}^{\mathrm{nrt}}:=\sum_{n\geq 0}\sum_{d\geq 0}\chi^{S_n}(\Mbar_{1,n}^{\mathrm{nrt}}(\mathbb{P}^r,d))q^d,$$ where $\Mbar_{1,n}^{\mathrm{nrt}}(\mathbb{P}^r,d)\subset \Mbar_{1,n}(\mathbb{P}^r,d)$ is the subspace of maps from smooth elliptic curves or nodal cycles of smooth rational curves. Such curves are characterised as the ones which do not contain any \textit{rational tails}, as the superscript indicates. 
To state our main theorem, we recall that $\Lambda = \QQ[[p_1, p_2, \ldots]],$ where $p_n$ is the $n$th power sum symmetric function. For a graded symmetric function $f \in \Lambda[[q]]$ we put $f' := \frac{\partial f}{\partial p_1}$ and $\dot{f} := \frac{\partial f}{\partial p_2}$. We use the symbol $\circ$ to denote plethysm of graded symmetric functions. For $f \in \Lambda[[q]]$, we set $\psi_n(f) := p_n \circ f$. We use the notation $\varphi$ for Euler's totient function and $\mu$ for the M\"obius function on the integers.
\begin{customthm}{A}\label{thm:main_formula} 
For $r \geq 1$, we have
\begin{equation*}
\mathsf{b}_{1, r}  = \mathsf{b}_{1, r}^{\mathrm{nrt}} \circ \left(p_1 + \frac{\mathsf{b}_{0, r}'}{r + 1} \right),
\end{equation*}
where
\begin{align*}
 \mathsf{b}_{1, r}^{\mathrm{nrt}}= & (r + 1)\left(\mathsf{a}_1 + \frac{\dot{\mathsf{a}}_0(\dot{\mathsf{a}}_0 + 1) + \frac{1}{4}\psi_2(\mathsf{a}_0'')}{1 - \psi_2(\mathsf{a}_0'')} - \frac{1}{2} \sum_{n \geq 1} \frac{\varphi(n)}{n}\log(1 - \psi_n(\mathsf{a}_0''))\right)  \\&+
\sum_{d \geq 2} q^d \sum_{k = 2}^{d} \left(\eta_{k,d}(r) \frac{(1 + 2 \dot{\mathsf{a}}_0)^2}{(1 - \psi_2(\mathsf{a}_0''))^{k/2 + 1}} +\sum_{j \mid k} \theta_{j,k,d}(r) \frac{1}{(1 - \psi_j(\mathsf{a}_0''))^{k/j}}\right).
 \end{align*}
The factors $\eta_{k,d}(r)$ and $\theta_{j, k,d}(r)$ are polynomials in $\QQ[r]$ given by the formulas
 \[\eta_{k,d}(r) = \frac{1}{4}\sum_{j \mid k}    \sum_{\substack{{\ell}\\{2\ell \mid j \text{ and } k \mid d\ell}}} \mu\left( \frac{j}{2\ell}\right)\binom{\frac{d\ell}{k} - 1}{\ell - 1}(r^{\ell + 1} + r^{\ell}), \]
 and
 \[\theta_{j,k,d}(r) = \frac{\varphi(j)}{2k}\sum_{i \mid \frac{k}{j}}\sum_{\substack{{\ell}\\{\ell \mid i\text{ and }k\mid d\ell}}} \mu\left( \frac{i}{\ell} \right) \binom{\frac{d\ell}{k} - 1}{\ell - 1}(r^\ell + (-1)^\ell r). \]
\end{customthm}

Theorem \ref{thm:main_formula} determines the generating function $\mathsf{b}_{1,r}$ after substituting the known formulas for $\mathsf{a}_0$, $\mathsf{a}_1$, and $\mathsf{b}_{0, r}'$ discussed above. Sample calculations obtained using Theorem \ref{thm:main_formula} are included in Tables \ref{table:degree1-frob-chars}, \ref{table:degree2-frob-chars}, and \ref{table:degree3-frob-chars} in \S\ref{sec:data}. For fixed $d$ and $n$, the multiplicity of each irreducible character of $S_n$ in $\chi^{S_n}(\Mbar_{1, n}(\P^r, d))$ is a polynomial function of $r$, of degree at most $d + 1$, with rational coefficients. In \cite[Theorem C]{KSGraphEnumeration}, we show further that this holds for $\chi^{S_n}(\Mbar_{g, n}(\P^r, d))$ for any $g$.

The formula for $\mathsf{b}_{1, r}^{\mathrm{nrt}}$ in Theorem \ref{thm:main_formula} follows from our calculation of the $S_n$-equivariant Serre characteristic\footnote{Our notation follows that of \cite{GetzlerPandharipande}. The Serre characteristic is a refinement of what has been called the Hodge--Deligne polynomial or $E$-polynomial.} of the torus fixed locus
\begin{equation}\label{eq:serrenrtdef}
\mathsf{e}^{S_n}\left(\Mbar_{1, n}^{\mathrm{nrt}}(\P^r, d)^{\C^\star}\right) := \sum_{i \geq 0} (-1)^i\left[H^i_c(\Mbar_{1, n}^{\mathrm{nrt}}(\P^r, d)^{\C^{*}})\right] \in K_0(\mathsf{M}) \otimes_{\Z} \Lambda,
\end{equation}
which is valued in the Grothendieck group of $S_n$-representations in the abelian category ${\mathsf{M} := \mathsf{MHS}_{\QQ}}$ of mixed Hodge structures over $\QQ$. See Theorem \ref{thm:main-tech} below.

We use tools from both algebraic geometry and algebraic combinatorics to study $\Mbar_{g, n}(\P^r, d)^{\C^\star}$. We define the key tools and outline how they are applied in more detail below.

\subsection{Graded $\mathbb{S}$-spaces and plethysm}\label{introsec:gradedSspaces}
An \textit{$\fS$-space} $\mathcal{X}$ is the data of an $S_n$-variety $\mathcal{X}(n)$ for each integer $n \geq 0$, while a \textit{graded $\fS$-space} is an $S_n$-variety $\mathcal{X}(n, d)$ for each pair of integers $n,d \geq 0$.  The Grothendieck ring $K_0(\mathsf{V},\fS)$ of $\fS$-spaces and the Grothendieck ring $K_0(\mathsf{V}, \mathbb{S})[[q]]$ of graded $\mathbb{S}$-spaces are introduced in \cite{GetzlerPandharipande}. The Serre characteristic
\begin{equation}\label{eq:SerreChar}
\mathsf{e}(\mathcal{X}):= \sum_{d \geq 0} q^d \sum_{n \geq 0} \mathsf{e}^{S_n}(\mathcal{X}(n, d)), 
\end{equation}
is a ring homomorphism ${K_0(\mathsf{V}, \fS) \to K_0(\mathsf{M}) \otimes \Lambda[[q]]}$. There is a composition operation on $K_0(\mathsf{V}, \mathbb{S})$ such that the Serre characteristic $\mathsf{e}$ becomes a map of composition algebras, namely 
\[ \mathsf{e}(\mathcal{X} \circ \mathcal{Y}) = \mathsf{e}(\mathcal{X}) \circ \mathsf{e}(\mathcal{Y}), \]
where the composition on the right is plethysm of symmetric functions. Write $\Mbar_{g}(r)$ for the graded $\mathbb{S}$-space whose $(n, d)$-component is $\Mbar_{g, n}(\P^r, d)$. A genus-$g$ stable map $f: C\to \mathbb{P}^r$ is said to \textit{have no rational tail} if the domain curve $C$ does not contain any non-trivial subcurve of arithmetic genus $g$. This is a combinatorial condition on dual graphs and defines a locally closed $\mathbb{S}$-subspace \[ \Mbar_{g}^{\mathrm{nrt}}(r) \subset \Mbar_{g}(r). \]

In \S\ref{section-rt}, we prove that attaching maps from trees of rational curves along marked points translates to the following formula involving the composition structure on $K_0(\mathsf{V}, \fS)[[q]]$.
    \begin{thm*}[Proposition \ref{prop-cut}]
    For a graded $\mathbb{S}$-space $\mathcal{X}$, define
    \[D \mathcal{X}(n, d) := \Res^{S_{n + 1}}_{S_n} \mathcal{X}(n + 1, d). \]
    Then in $K_0(\mathsf{V},\mathbb{S})[[q]],$ we have the formula
        \[ [\Mbar_{g}(r)] = [\Mbar_{g}^{\mathrm{nrt}}(r)] \circ \left(  [\Spec\C] + \frac{[D\Mbar_0(r)]}{[\P^r]}\right),\]
        where $\Spec \C$ has the trivial action of $S_1$, in $d$-degree $0$.
    \end{thm*}
 
    Since a recursive formula for $\mathsf{e}(\Mbar_0(r))$ has been determined in \cite{GetzlerPandharipande}, this formula reduces the calculation of $\mathsf{e}(\Mbar_{g}(r))$ to that of $\mathsf{e}(\Mbar_g^{\mathrm{nrt}}(r))$.

    \subsection{$\mathbb{C}^\star$-localisation and type-$B$ symmetric functions}\label{introsec:localisation}
    In genus one, we express the geometry of the space of torus-fixed maps from curves without rational tails using type-$B$ combinatorics. The full graded $\mathbb{S}$-space $\Mbar_{1}^{\mathrm{nrt}}(r)$ has subtle geometry: it includes the space of maps $\M_{1,n}(\P^r, d)$ from smooth elliptic curves for all $n$ and $d$, as well as maps to $\P^r$ from arbitrary cycles of smooth pointed rational curves. The picture simplifies upon passing to the $\C^\star$-fixed locus. We set
    \begin{equation}\label{eq:scriptBdefn}
        \mathscr{B}_r := \mathsf{e}\left({\Mbar_{1}^{\mathrm{nrt}}(r)^{\C^\star}}\right) \in \Lambda[[q]]
    \end{equation}
    for the Serre characteristic of the locus of $\C^\star$-fixed genus one maps without rational tails. The starting point for our calculation of $\mathscr{B}_r$ is the graph stratification of $\Mbar_{g, n}(\P^r, d)^{\C^\star}$ described by Graber and Pandharipande \cite{graberpand}. The strata correspond to $(\P^r)^{\C^\star}$-coloured, edge-weighted graphs of total weight $d$, which we call \textit{localisation graphs}. If $\mathcal{G}$ is a localisation graph, then the action of $\Aut(\mathcal{G})$ on the set of half-edges $H(\mathcal{G})$ of $\mathcal{G}$ induces an embedding
    \begin{equation}\label{eq:flagsofgraph}
    \Aut(\mathcal{G}) \hookrightarrow  \mathrm{Perm}(H(\mathcal{G})).
    \end{equation}
    If we restrict to $\Mbar_{1}^{\mathrm{nrt}}(r)$, the graph $\mathcal{G}$ must be a cycle of length $k \geq 2$. Thus (\ref{eq:flagsofgraph}) factors through the hyperoctahedral group $B_k := {S_2 \wr S_k}$. Therefore, when specialised to $\Mbar_{1}^{\mathrm{nrt}}(r)^{\C^\star}$, equation (\ref{eq:flagsofgraph}) assigns to each localisation graph $\mathcal{G}$ a character of $B_k$, which encodes the action of $\Aut(\mathcal{G})$ on $H(\mathcal{G})$. We find that the $S_n$-action on the stratum determined by $\mathcal{G}$ depends only on this character. We package localisation graphs and their characters into a graded $\mathbb{B}$-space $\mathsf{Dih}_r$, where $\mathbb{B}$-spaces are the natural type-$B$ analogues of the $\mathbb{S}$-spaces discussed in \S\ref{introsec:gradedSspaces}. The moduli in the stratum defined by $\mathcal{G}$ comes from the choice of a pointed chain of rational curves for each vertex. The moduli space of chains is an $S_2 \times \mathbb{S}$-space $\mathsf{Cat}$, where $S_2$ acts by switching the orientation of the chain. We prove the following motivic formula for the graded $\mathbb{S}$-space $\Mbar_{1}^{\mathrm{nrt}}(r)^{\C^\star}$.
    \begin{thm*}[Theorem \ref{thm-locnrt}] We have the following equality in $K_0(\mathsf{V}, \fS)[[q]]$:
    \[ [\Mbar_{1}^{\mathrm{nrt}}(r)^{\C^\star}] = [(\P^r)^{\C^\star}] \cdot  
 [\Mbar_{1}^{\mathrm{nrt}}(0)] + [\mathsf{Dih}_r ]\circ_{S_2} [\mathsf{Cat}]. \]
    Above, $\circ_{S_2}$ is the composition of a graded $\mathbb{B}$-space with an $S_2 \times \fS$-space, discussed in \S\ref{sec-GradedPlethysm}. 
\end{thm*}

The term $\Mbar_1^{\mathrm{nrt}}(0)$ parameterises maps to a point, and its Serre characteristic has been computed \cite{GetzlerSemiClassical, semiclassicalremark}, as has the Serre characteristic of $\Cat$ \cite{BergstromMinabe1}. Theorem \ref{thm-locnrt} thus reduces the calculation of $\mathscr{B}_r$ to the calculation of $\mathsf{e}(\mathsf{Dih}_r)$ and its wreath product plethysm with $\mathsf{e}(\Cat)$. The Serre characteristic $\mathsf{e}(\mathsf{Dih}_r)$ is valued in ${K_0(\mathsf{M}) \otimes \Lambda(S_2)[[q]]}$, where $\Lambda(S_2)$ is the ring of \textit{wreath product symmetric functions} studied by Macdonald \cite{Macdonald}.

\subsection{Enumerative combinatorics of cycles}\label{introsec:enumerative}
In \S\ref{sec-calculation}, we use graph enumeration techniques to compute $\mathsf{e}(\mathsf{Dih}_r)$ and in \S\ref{sec-cat} we compute ${\mathsf{e}(\mathsf{Dih}_r) \circ_{S_2} \mathsf{e}(\Cat)}$. These calculations lead to our main technical result, stated as Theorem \ref{thm:main-tech} below. Together with Proposition \ref{prop-cut}, it implies Theorem \ref{thm:main_formula}. We set the notation
 \[ \mathscr{A}_g := \sum_{n > 2 - 2g} \mathsf{e}^{S_n}(\M_{g, n}), \]
 so $\mathscr{A}_g$ is the natural lift of the series $\mathsf{a}_g$ in (\ref{eq:sfa}) to a virtual mixed Hodge structure in $K_0(\mathsf{M}) \otimes \Lambda$.
\begin{customthm}{B}
We have
\begin{align*}
\mathscr{B}_r &= (r + 1)\left(\mathscr{A}_1 + \frac{\dot{\mathscr{A}}_0(\dot{\mathscr{A}}_0 + 1) + \frac{1}{4}\psi_2(\mathscr{A}_0'')}{1 - \psi_2(\mathscr{A}_0'')} - \frac{1}{2} \sum_{n \geq 1} \frac{\varphi(n)}{n}\log(1 - \psi_n(\mathscr{A}_0''))\right) \\&\phantom{space}+ \sum_{d \geq 2} q^d \sum_{k = 2}^{d} \left(\eta_{k, d}(r) \frac{(1 + 2 \dot{\mathscr{A}}_0)^2}{(1 - \psi_2(\mathscr{A}_0''))^{k/2 + 1}} + \sum_{j \mid k} \theta_{j,k,d}(r)\frac{1}{(1 - \psi_j(\mathscr{A}_0''))^{k/j}}\right).
\end{align*}
\end{customthm}
 The Serre characteristic of $\mathsf{Dih}_r$ is the solution to a graph enumeration problem: it is a count of isomorphism classes of localisation graphs, where each graph $\mathcal{G}$ is weighted by the type-$B$ character of the action of $\Aut(\mathcal{G})$ on $H(\mathcal{G})$. We can write this as
\[ \mathsf{e}(\mathsf{Dih}_r) = \sum_{d \geq 2} q^d \sum_{k = 2}^{d} \sum_{[H] \leq D_k} \gamma_{r, d}^{[H]}(k) \ch_{2\wr k}(\Ind_{H}^{B_k} \mathrm{Triv}_H ). \]
The innermost sum is over conjugacy classes of subgroups of the dihedral group $D_k$, viewed as a subgroup of $B_k$. The factor $\gamma_{r, d}^{[H]}(k)$ denotes the number of isomorphism classes of localisation graphs such that the action of $\Aut(\mathcal{G})$ on $H(\mathcal{G})$ has the same character as the action of $H$ on the half-edges of a fixed $k$-cycle.

In \S\ref{sec-calculation}, we derive an explicit formula for $\mathsf{e}(\mathsf{Dih}_r)$ by computing the numbers $\gamma_{r, d}^{[H]}(k)$ for all possible values of the parameters, as well as the associated characters. To count decorations with prescribed automorphisms, we use standard techniques in combinatorial enumeration. The most important step is to apply M\"obius inversion to the subgroup lattice $\mathcal{L}(D_k)$ of $D_k$, which gives formulas for $\gamma_{r, d}^{[H]}(k)$ in terms of the following ingredients:
\begin{itemize}
    \item binomial coefficients which encode partitions of the stable map degree along the edges,
    \item chromatic polynomials of paths and cycles which count $(\P^r)^\mathbb{C^\star}$-colourings of quotients of localisation graphs by automorphisms, and
    \item the M\"obius function $\mu_{D_k}$ of $\mathcal{L}(D_k)$.
\end{itemize}

The combinatorial factors $\theta_{j,k,d}(r)$ and $\eta_{k,d}(r)$ in Theorem \ref{thm:main_formula} arise from combining these formulas with the associated characters, as well as the wreath product plethysm with $\mathsf{e}(\mathsf{Cat})$. All of the relevant calculations are assembled to prove Theorems \ref{thm:main_formula} and \ref{thm:main-tech} in \S\ref{sec-cat}.

\subsection{History and context}\label{introsec:relatedwork}
We are inspired by work of Getzler \cite{GetzlerGenusZero, GetzlerSemiClassical} and Getzler--Kapranov \cite{GetzlerKapranov} from the 1990s, using symmetric functions to work systematically with $S_n$-equivariant graph stratifications of moduli spaces of stable curves. Getzler and Kapranov work in the setting of \textit{modular operads}, and give a formula \cite[Theorem 8.13]{GetzlerKapranov} which encodes how the two generating functions
\begin{equation}\label{eq:GKgenfuns}
\mathsf{m}:= \sum_{2g - 2 + n >0}\mathsf{e}^{S_n}(\M_{g, n}) q^{g - 1} \quad \mbox{and} \quad \overline{\mathsf{m}}:= \sum_{2g - 2 + n >0}\mathsf{e}^{S_n}(\Mbar_{g, n}) q^{g - 1}
\end{equation}
in $K_0(\mathsf{M}) \otimes \Lambda[[q]]$ determine one another, in terms of plethysm and other algebraic operations on $\Lambda[[q]]$. Getzler computed $\mathsf{m}$ for $g \leq 1$ geometrically, and he then used the symmetric function formalism to compute $\overline{\mathsf{m}}$ for $g \leq 1$ \cite{GetzlerGenusZero, GetzlerSemiClassical, GetzlerMHM}. In higher genus, the framework developed by Getzler and Kapranov has become a standard tool for passing cohomological information between moduli spaces of smooth and stable curves: see \cite[\S 2]{BergstromTommasi} or \cite[\S 2.5]{PW1}. Gorsky has computed the $S_n$-equivariant Euler characteristic $\chi^{S_n}(\M_{g, n})$ \cite{Gorsky} for $g > 1$ in terms of the orbifold calculation of Harer and Zagier \cite{HarerZagier}. Gorsky's formula and Getzler's results can in principle be passed through Getzler and Kapranov's formula to determine $\chi^{S_n}(\Mbar_{g, n})$ for arbitrary $g$ and $n$, though this requires daunting calculations \cite{Diaconu}. 

When it comes to Kontsevich moduli spaces, the most complete results on their topology and intersection theory are limited to genus zero. The $S_n$-equivariant Serre characteristics of both $\M_{0, n}(\P^r, d)$ and $\Mbar_{0, n}(\P^r, d)$ were computed by Getzler and Pandharipande \cite{GetzlerPandharipande}, by introducing symmetric functions in the stable map setting. Since $\Mbar_{0, n}(\P^r, d)$ is smooth and proper, its Serre characteristic  determines its Hodge decomposition. Due to its role in enumerative geometry, there has been much interest in the intersection theory of $\Mbar_{0,n}(\P^r,d)$ \cite{PandharipandeQ, BehrendOH, BF, OpreaPrTaur, OpreaFlags, MustataMustata1,MustataMustata2, KM}. The rational cohomology of the interior is also well-understood \cite{NonlinearGrassmannian, FarbWolfson}, as is the integral Chow ring for odd values of $d$ \cite{CavalieriFulghesu}. There has also been recent progress on the cohomology of moduli spaces of logarithmic stable maps to $\P^r$ in genus zero \cite{KannanP1, DivisorsAndCurves, Devkota}. 


In positive genus, most of the literature on the cohomology of $\M_{g, n}(\P^r, d)$ and $\Mbar_{g, n}(\P^r, d)$ is focused on specific classes of interest in Gromov--Witten theory \cite{graberpand, HodgeIntegrals, Lambdag}. There are few known calculations that apply to the entire cohomology ring, likely due to the difficult singularities and excess components of the moduli space. Vakil and Zinger desingularised the so-called \textit{main component} of $\Mbar_{1, n}(\P^r, d)$ \cite{VakilZinger} to produce a smooth compactification $\widetilde{\M}_{1, n}(\P^r, d)$ of $\M_{1, n}(\P^r, d)$ with tractable geometry. Their work was later given a modular interpretation in logarithmic geometry by Ranganathan, Santos--Parker, and Wise \cite{rspw} and extended to genus two by Battistella and Carocci \cite{BattistellaCarocci}. We used the compactification $\widetilde{\M}_{1, n}(\P^r, d)$ to compute the boundary complex of $\M_{1, n}(\P^r, d)$ in \cite{vzdualcomplex}, and we deduced the vanishing of the weight-zero compactly-supported cohomology $W_0H^\star_c(\M_{1, n}(\P^r, d);\QQ)$. The second-named author computes the pure-weight rational cohomology $W_{\star}H^\star(\Mbar_{1,n}(\mathbb{P}^r,d); \mathbb{Q})$ in \cite{ds}. Bae has recently initiated the study of tautological rings of Kontsevich spaces in positive genus \cite{Bae}.

Topologists have proven many interesting theorems about moduli spaces of maps, going back to Segal \cite{Segal}. Results in this direction usually focus on the space of maps from a fixed curve \cite{Guest, Banerjee}. Tosteson \cite{tosteson} has proven one of the few positive genus structural results about the rational cohomology of $\Mbar_{g, n}(\P^r, d)$ of which we are aware: for fixed $g, r, d$, and $i$, he proves that the sequence of $S_n$-representations determined by $H_i(\Mbar_{g, n}(\P^r, d); \QQ)$ satisfies a form of representation stability, in the sense of Sam and Snowden \cite{SamSnowden}. His theorem implies several constraints on the asymptotic behaviour of these representations.

On the combinatorial side, our technique of enumerating graphs with specified automorphisms by counting quotient graphs is inspired by the calculation of the $S_n$-equivariant top-weight Euler characteristic of $\M_{g, n}$ in \cite{CFGP}. Our method is reminiscent of P\'olya enumeration \cite{Polya,PolyaRead}. Adjacent problems on the enumeration of graph colourings and necklaces have been addressed \cite{Moreau, Robinson, Hanlon}. Our use of wreath product symmetric functions is inspired by work of Macdonald \cite{MACDONALD1980}, Petersen's combinatorial proof \cite{semiclassicalremark} of Getzler's formula for $\mathsf{e}^{S_n}(\Mbar_{1, n})$, and his work on operad structures for moduli spaces of admissible $G$-covers \cite{PetersenOperadGCovers}.

\subsection{Further directions} We outline four natural directions in which to generalise Theorems \ref{thm:main_formula} and \ref{thm:main-tech}. 


\subsubsection*{Virtual Euler characteristics} It is possible that our approach to the topological Euler characteristics applies to the virtual Euler characteristics of $\Mbar_{1,n}(\mathbb{P}^r,d)$ \cite{GiventalWDVV}, \cite{Lee}. 
Among the key ingredients of our method, torus localisation \cite{GiventalII, AndersonChenTseng} and $S_n$-actions \cite{GiventalI} have been explored in quantum K-theory, and the virtual Euler characteristic satisfies cut-and-paste properties relative to the forgetful morphism $\Mbar_{1,n}(\mathbb{P}^r,d)\to \mathfrak{M}_{1,n}$. Hence, there is reason to expect that analogous calculations could contribute to positive genus K-theoretic Gromov--Witten invariants.
\subsubsection*{Higher genus calculations }
The Graber--Pandharipande description of $\Mbar_{g, n}(\P^r, d)^{\C^\star}$ implies that its Serre characteristic is determined by products of Serre characteristics of moduli spaces of stable curves, modulo finite groups. In \cite{KSGraphEnumeration}, we prove a graph sum formula expressing the generating function $\sum_{g,n}\chi^{S_n}(\Mbar_{g, n}(\P^r, d)) t^{g-1}$ in terms of the generating function $\overline{\mathsf{m}}$ from (\ref{eq:GKgenfuns}) and certain graph sums. The present work calculates the $t$-degree $0$ term by computing the genus one graph sums explicitly after separating rational tails. Algorithms to compute the graph sums in higher genus would make many more calculations possible.
\subsubsection*{Comparing compactifications} There are alternative modular compactifications of $\M_{1, n}(\P^r, d)$, most notably the moduli space of quasimaps $\overline{\mathcal{Q}}_{1,n}(\mathbb{P}^r,d)$ \cite{MarianOpreaPandharipande} and the Vakil--Zinger desingularisation of $\Mbar_{1, n}(\P^r, d)$ \cite{VakilZinger, rspw, vzdualcomplex}. It would be interesting to compare the topology and combinatorics of these spaces. The locus $\Mbar_{1}^{\mathrm{nrt}}(r)$ is common to all three compactifications, so one could hope for plethystic formulas comparing the Serre characteristics. We expect that our techniques could eventually give motivic analogues of wall-crossing formulas relating Gromov--Witten and quasimap invariants \cite{Qmapwallcrossing}, similar to the motivic wall-crossing formulas for Hassett spaces obtained in \cite{hodgehl}. 

\subsubsection*{Combinatorial complexity of strata} The enumeration of coloured cycles in our work gives the count of $\C^*$-fixed strata in $\Mbar_{1,n}^{\mathrm{nrt}}(\mathbb{P}^r,d)$. It would be interesting if our techniques could be adapted to count all $\C^\star$-fixed strata in $\Mbar_{1,n}(\mathbb{P}^r,d)$.
Such enumeration problems encode the combinatorial complexity of the compactification. In the case of stable curves $\mathcal{M}_{g,n}\subset\Mbar_{g,n}$,  Aluffi--Marcolli--Nascimento \cite{AluMarNas} enumerate the boundary strata of fixed dimension when $g = 0$, and Getzler enumerates all boundary strata when $g = 1$ \cite{GetzlerSemiClassical}. To our knowledge, this direction has not been pursued for mapping spaces.

\subsection*{Acknowledgements} We thank Dhruv Ranganathan and Burt Totaro for edifying conversations. We also thank Sarah Brauner, Melody Chan, Megan Chang-Lee, Dan Petersen, Dhruv Ranganathan, and Burt Totaro for providing feedback on a draft of this article. We would like to acknowledge helpful comments from anonymous referees. SK is supported by NSF DMS-2401850, and he is grateful for the hospitality extended by members of the DPMMS at the University of Cambridge, where much of this work was carried out. TS is supported by a Cambridge Trust international scholarship. 

\section{Symmetric functions and wreath products}\label{section-sym}
In this section we recall the background on symmetric functions and their wreath product analogues that we will use throughout the paper. We refer to \cite{MACDONALD1980, Macdonald} and \cite[§7]{GetzlerKapranov} for more detailed discussions of this material. Given a partition $\lambda \vdash n$, we write $\lambda_i$ for the number of times $i$ appears in $\lambda$, so $\sum_{i= 1}^{n} i\lambda_i = n$, and we write $|\lambda|:=n$. Given a permutation $\sigma \in S_n$, we write $\lambda(\sigma) \vdash n$ for the cycle type of $\sigma$.

\subsection{Symmetric functions and $S_n$-equivariant characteristics}\label{subsec-sym}We use $\Lambda$ to denote the ring of completed\footnote{Many authors write $\hat{\Lambda}$ for this ring to indicate that it is the degree completion of the usual ring of symmetric functions.} symmetric functions $$\Lambda := \lim_{\longleftarrow} \QQ[[x_1, \ldots, x_n]]^{S_n}.$$
For $i>0$, let $p_i:=\sum_{k>0}x^i_k$, then $$\Lambda = \QQ[[p_1,p_2\dots ]].$$ It is convenient to identify $p_0=1$. We set $\Lambda_n \subset \Lambda$ for the homogeneous degree-$n$ part of $\Lambda$, so $p_n \in \Lambda_{n}$ for all $n \geq 0$. If $V$ is an $S_n$-representation in the category of $\QQ$-vector spaces, we define its \textit{Frobenius characteristic} by
\begin{equation}
\ch_n(V) := \frac{1}{n!} \sum_{\sigma \in S_n} \mathrm{Tr}_V(\sigma) p_{\lambda(\sigma)} \in \Lambda, 
\end{equation}
where $p_\lambda := \prod_{i} p_{i}^{\lambda_i}$ for $\lambda \vdash n$. This gives a homogeneous symmetric function of degree $n$. If $\mathrm{Triv}_n$ is the trivial representation of $S_n$, we set
\begin{equation}\label{eq:homogeneousdefn}
h_n := \ch_n(\mathrm{Triv}_n) \in \Lambda_n
\end{equation}
for the $n$th \textit{homogeneous} symmetric function. The elements $h_i$ freely generate $\Lambda$ as a power series ring. If $V = V_{\lambda}$ is the Specht module corresponding to $\lambda$, we set
\begin{equation}\label{eq:schurdefn}
s_\lambda := \ch_n(V_\lambda) \in \Lambda_n
\end{equation}
for the \textit{Schur function} corresponding to the partition $\lambda$. There is an inner product $\langle-, -\rangle$, called the \textit{Hall inner product}, encoding the inner product of characters of $S_n$. In particular, the Schur functions $s_\lambda$ for $\lambda \vdash n$ form an orthonormal basis for the homogeneous $\Lambda_n$ with respect to this inner product. For an arbitrary $S_n$-representation $V$, we have that
\begin{equation}\label{eq:Hall}
\langle s_\lambda , \ch_n(V) \rangle = \dim_{\QQ}(\mathrm{Hom}_{S_n}(V_{\lambda}, V)),
\end{equation}
i.e. the inner product of $\ch_n(V)$ with $s_\lambda$ is given by the multiplicity of the irreducible $V_{\lambda}$ in $V$.
Now suppose that $V$ is an $S_n$-representation in the abelian category $\mathsf{M}$ of mixed Hodge structures over $\mathbb{Q}$. The Peter--Weyl theorem for representations of finite groups in $\mathsf{M}$ \cite[Theorem 3.2]{GetzlerPreprint} implies that for each $\lambda \vdash n$, there exists a unique mixed Hodge structure $M_\lambda$ such that
\[ V \cong \bigoplus_{\lambda \vdash n} M_{\lambda} \otimes V_\lambda, \]
where $V_\lambda$ is a Specht module, considered as a mixed Hodge structure of type $(0, 0)$. We then set
\begin{equation}\label{eq:frobmhsdefn}
\ch_n^{\mathsf{M}}(V) := \sum_{\lambda \vdash n} [M_\lambda] \otimes s_\lambda \in K_0(\mathsf{M}) \otimes \Lambda_n
\end{equation}
for the enriched \textit{Frobenius characteristic} of $V$. If $X$ is an $S_n$-variety over $\C$, then its compactly-supported singular cohomology  $H^\star_c(X; \mathbb{Q})$ forms a graded $S_n$-representation in the category $\mathsf{M}$. We define the $S_n$-equivariant Serre characteristic of $X$ as
\begin{equation}\label{eq:serrechardefn}
\mathsf{e}^{S_n}(X):= \sum_{i \geq 0} (-1)^i\mathrm{ch}_n^{\mathsf{M}}( H^i_c(X;\QQ))  \in K_0(\mathsf{M}) \otimes \Lambda_{n}.
\end{equation}
If $X$ is proper and has only finite quotient singularities, then $\mathsf{e}^{S_n}(X)$ determines the $S_n$-equivariant Hodge decomposition of $X$. Using (\ref{eq:frobmhsdefn}), we can define unique classes $\mathsf{e}_{\lambda}(X) \in K_0(\mathsf{M})$ which satisfy
\begin{equation}\label{eq:serremultiplicity}
\mathsf{e}^{S_n}(X) = \sum_{\lambda \vdash n} \mathsf{e}_{\lambda}(X) \otimes s_{\lambda},
\end{equation}
and the element $\mathsf{e}_{\lambda}(X)$ should be thought of as the virtual multiplicity of the Specht module $V_{\lambda}$ in $\mathsf{e}^{S_n}(X)$.
There is a homomorphism \begin{equation}\label{eq:dimserre}
K_0(\mathsf{M}) \otimes \Lambda \to 
K_0(\mathsf{Vect}_\QQ) \otimes \Lambda = \Lambda
\end{equation}
determined by taking the rank of a mixed Hodge structure. For $X$ an $S_n$-variety, we define the $S_n$-equivariant Euler characteristic of $X$ (also called the Frobenius characteristic of $X$) as the image $\chi^{S_n}(X)$ of $\mathsf{e}^{S_n}(X)$ under (\ref{eq:dimserre}); this is equivalent to the definition (\ref{eq:frobchardefn}) given in the introduction. We can also define \[\chi_{\lambda}(X) \in \Z = K_0(\mathsf{Vect}_{\QQ})\] to be the image of the virtual mixed Hodge structure $\mathsf{e}_{\lambda}(X)$. Then we have
\[ \chi_{\lambda}(X) = \langle s_{\lambda}, \chi^{S_n}(X) \rangle, \]
with respect to the Hall inner product (\ref{eq:Hall}).

The \emph{plethysm} operation $f \circ g$ is defined for $f, g \in \Lambda$ whenever $f$ has bounded degree or $g \in \oplus_{n > 0} \Lambda_n$. It is defined and characterised by the following properties:\begin{itemize}
    \item for any $g\in \Lambda$, the map $\circ g: \Lambda\to \Lambda$ given by $f\mapsto f\circ g$ is an algebra homomorphism;
    \item for all $n$, the map $p_n\circ: \Lambda\to \Lambda$ given by $f\mapsto p_n\circ f$ is an algebra homomorphism;
    \item $p_n\circ p_m = p_{mn}$.
\end{itemize}
Plethysm extends to $K_0(\mathsf{M}) \otimes \Lambda$ by setting $h_n \circ [V] = [\mathrm{Sym}^{n}(V)]$ for $[V]$ the class of a mixed Hodge structure $V$.

\subsection{$\mathbb{S}$-modules} An $\mathbb{S}$\emph{-module} $\mathcal{V}$ in $\mathsf{M}$ is a sequence of $S_n$-representations $\mathcal{V}(n)$ in $\mathsf{M}$ for each $n \geq 0$. There is a box product operation  \[(\mathcal{V}\boxtimes \mathcal{W})(n):= \bigoplus_{j = 0}^{n} \Ind_{S_{j}\times S_{n-j}}^{S_n} \mathcal{V}(j)\otimes \mathcal{W}(n-j),\] which endows the Grothendieck group of $\mathbb{S}$-modules in $\mathsf{M}$ with a ring structure. 

In fact, the Grothendieck ring is isomorphic to $K_0(\mathsf{M}) \otimes \Lambda$, with the ring isomorphism given by 
\begin{equation}\label{eq:K_0(Smod)}
[\mathcal{V}] \mapsto \ch(\mathcal{V}):= \sum_{n \geq 0} \mathrm{ch}_n^{\mathsf{M}}(\mathcal{V}(n)) \in K_0(\mathsf{M}) \otimes \Lambda,
\end{equation}
where $\ch_n^{\mathsf{M}}$ is defined as in (\ref{eq:frobmhsdefn}). Under the isomorphism (\ref{eq:K_0(Smod)}), plethysm in $K_0(\mathsf{M}) \otimes \Lambda$ corresponds to the following composition structure $\mathcal{V} \circ \mathcal{W}$, defined on $\mathbb{S}$-modules $\mathcal{V}$ and $\mathcal{W}$ whenever $\mathcal{V}$ has bounded degree or $\mathcal{W}(0) = 0$: \[(\mathcal{V}\circ \mathcal{W})(n):= \bigoplus_{j = 0}^{\infty} \left(\mathcal{V}(j)\otimes_{S_j} \mathcal{W}^{\boxtimes j}(n) \right).\] The notation $\otimes_{S_j}$ means taking the $S_j$-coinvariants of the tensor product representation $\mathcal{V}(j)\otimes \mathcal{W}^{\boxtimes j}(n)$. Composition of $\mathbb{S}$-modules is compatible with plethysm of symmetric functions in the sense that $$\ch(\mathcal{V}\circ\mathcal{W}) = \ch(\mathcal{V})\circ \ch(\mathcal{W}).$$

\subsection{$S_2\times S_n$-symmetric functions} 
Consider the tensor product $K_0(\mathsf{M}) \otimes\Lambda\otimes\Lambda$, thought of as a ring of bisymmetric functions. The characteristics associated to $S_2\times \mathbb{S}$-spaces generate the vector subspace $K_0(\mathsf{M}) \otimes \Lambda_2 \otimes \Lambda$ spanned by bihomogeneous polynomials of bidegrees $(2,i)$ for all $i\geq 0$. There are two convenient linear bases for $\Lambda_2$ that we will use: the first is \[\left\{\frac{1}{2}p_1^2, \frac{1}{2}p_2\right\}, \] and the second is $$\left\{\ch_2(\mathrm{Triv}_{S_2})=\frac{1}{2}(p_1^2 + p_2),\,\ch_2(\mathrm{Alt}_{S_2})=\frac{1}{2}(p_1^2 - p_2)\right\}.$$

The Frobenius characteristic of an $S_2\times S_n$-representation in $\mathsf{Vect}_\QQ$ is defined analogously to the case of $S_n$:
\begin{equation}
    \ch_{2,n}(V):=\frac{1}{2n!}\sum_{(\sigma, \tau) \in S_{2} \times S_n} \mathrm{Tr}_V(\sigma, \tau) p_{\lambda(\sigma)} \otimes  p_{\lambda(\tau)} \in  \Lambda_2 \otimes \Lambda.
\end{equation}
Applying the Peter--Weyl theorem to $S_2 \times S_n$-representations in $\mathsf{M}$, there exists a unique mixed Hodge structure $M_{\lambda, \nu}$ for each pair of partitions $\lambda \vdash 2$, $\nu \vdash n$, such that
\[V \cong \bigoplus_{\substack{{\lambda \vdash 2}\\{\nu \vdash n} }} M_{\lambda, \nu} \otimes V_\lambda \otimes V_{\nu},\]
where $V_\lambda$ and $V_\nu$ are Specht modules. We use this to define the Frobenius characteristic
\begin{equation}\label{eq:S2timesSfrobchar}
\ch_{2, n}^{\mathsf{M}}(V) := \sum_{\substack{{\lambda \vdash 2}\\{\nu \vdash n} }} [M_{\lambda, \nu}] \otimes s_\lambda \otimes s_{\nu} \in K_0(\mathsf{M}) \otimes \Lambda_2 \otimes \Lambda.
\end{equation}
For $X$ an $S_2 \times S_n$-space, we define the $S_2 \times S_n$-equivariant Serre characteristic by
$$\mathsf{e}^{S_2 \times S_n}(X):=\sum_{i\geq 0}(-1)^i \ch_{2,n}^{\mathsf{M}}(H^i_c(X; \mathbb{Q})) \in K_0(\mathsf{M}) \otimes \Lambda_2 \otimes \Lambda.$$ Then the Frobenius characteristic $\chi^{S_2 \times S_n}(X) \in \Lambda_2 \otimes \Lambda$ is obtained by taking ranks of mixed Hodge structures. Later we will need the following formula, which computes $\mathsf{e}^{S_2 \times S_n}(Y)$ if \[Y = \Res^{S_{n + 2}}_{S_2 \times S_n} X,\] for an $S_{n + 2}$-space $X$.
\begin{lem}\label{lem-sn+2}
Let $X$ be an $S_{n+2}$-space, then $$\mathsf{e}^{S_2 \times S_n}(\Res_{S_2\times S_n}^{S_{n+2}}X) = \frac{1}{2}p_1^2\otimes\left(\frac{\partial^2 \mathsf{e}^{S_{n+2}}(X)}{\partial p_1^2}\right)+ p_2\otimes \left(\frac{\partial \mathsf{e}^{S_{n+2}}(X)}{\partial p_2}\right).$$
\end{lem}
\begin{proof}
    It suffices to prove the formula for the Frobenius characteristic of an $S_{n+2}$-representation $V$. Since there are only two irreducible representations of $S_2$, we have that $$\Res_{S_2\times S_n}^{S_{n+2}}V\cong \mathrm{Triv}_{S_2}\otimes \mathrm{Hom}_{S_2}(\mathrm{Triv}_{S_2},V)\oplus \mathrm{Alt}_{S_2}\otimes \mathrm{Hom}_{S_2}(\mathrm{Alt}_{S_2},V).$$ Here both $\mathrm{Hom}_{S_2}(\mathrm{Triv}_{S_2},V)$ and $\mathrm{Hom}_{S_2}(\mathrm{Alt}_{S_2},V)$ are $S_n$-representations.

    From \cite[Proposition 8.10]{GetzlerKapranov}, for an $S_2$-representation $W$, we have $$\ch_n(\mathrm{Hom}_{S_2}(W, V))=D(\ch_2 (W))(\ch_{n+2}(V)),$$ where $D(-)$ is a linear operator defined on a basis by $$D(\ch_2(\mathrm{Triv}_{S_2}))=\frac{1}{2}\frac{\partial^2}{\partial p_1^2}+\frac{\partial}{\partial p_2}\quad\mbox{and}\quad D(\ch_2(\mathrm{Alt}_{S_2}))=\frac{1}{2}\frac{\partial^2}{\partial p_1^2}-\frac{\partial}{\partial p_2}.$$ Applying these operators, we see that
    \begin{align*}
    \ch_{2, n}(\mathrm{Res}^{S_{n + 2}}_{S_2 \times S_n} X) &= \frac{1}{2}\left( p_1^2 + p_2 \right) \otimes \left(\frac{1}{2} \frac{\partial^2}{\partial p_1^2} \ch_{n + 2}(X) + \frac{\partial}{\partial p_2}\ch_{n + 2}(X) \right)  \\&+ \frac{1}{2}\left( p_1^2- p_2 \right) \otimes \left(\frac{1}{2} \frac{\partial^2}{\partial p_1^2} \ch_{n + 2}(X) - \frac{\partial}{\partial p_2}\ch_{n + 2}(X) \right).
    \end{align*}
    The formula follows after rearranging.
\end{proof}

\subsection{Type-$B$ symmetric functions}
Symmetric functions associated to representations of the hyperoctahedral group $B_n:= S_2 \wr S_n$ are needed for our study of decorated necklaces which arise in the localisation formula for $\Mbar_{1, n}(\P^r, d)$. Our presentation is synthesised from that of Petersen \cite[\S3]{semiclassicalremark} and Macdonald \cite[Appendix B]{MACDONALD1980, Macdonald}.

Write $S_2$ multiplicatively as $S_2 = \{\pm 1\}$; then the elements of $B_n$ are tuples $(g_1, \ldots, g_n; x)$ where $g_i = \pm1$ and $x \in S_n$. Let $\Lambda(S_2)$ be the free power series ring over $\QQ$ generated by $\pp_n=p_n(1)$ and $\qq_n=p_n(-1)$ for $n\geq 1$. We endow $\Lambda(S_2)$ with a grading by setting $\pp_i$ and $\qq_i$ to have degree $i$. There are two natural inclusion maps \begin{equation}\label{eq:wreathinclusions}
\iota_{\mathfrak{p}}, \iota_{\mathfrak{q}}: \Lambda \to \Lambda(S_2).
\end{equation}
The generalisation of the cycle map to this setting is $\Psi: B_n\to \Lambda(S_2)$ defined by $$(g_1,\dots,g_n; x)\mapsto \prod_{\sigma \text{ a cycle in $x$}} p_{|\sigma|}\left(\prod_{j\in \sigma} g_j\right).$$ Given a $B_n$-representation $V$ in $\mathsf{Vect}_{\QQ}$, the type-$B$ Frobenius characteristic is given by 
\begin{equation}\label{eq:wreathfrobchar}
\ch_{2\wr n}(V):=\frac{1}{2^n\cdot n!}\sum_{(\vec{g}; x) \in B_n}\mathrm{Tr}_V((\vec{g}; x))\Psi((\vec{g}; x)) \in \Lambda(S_2).
\end{equation}
For each ordered pair $(\lambda, \nu)$ of partitions such that $|\lambda| + |\nu| = n$, there is a unique irreducible representation $V_{\lambda, \nu}$ of $B_n$, and these form a complete list of the irreducible representations of $B_n$ \cite[§1.B.9]{Macdonald}. For these representations, we have
\[ \ch_{2\wr n}(V_{\lambda, \nu}) = \iota_\mathfrak{p}(s_{\lambda}) \iota_{\mathfrak{q}}(s_{\nu}), \]
where $s_{\lambda}, s_{\nu} \in \Lambda$ are the ordinary Schur functions and the maps $\iota_\pp, \iota_\qq$ are as in (\ref{eq:wreathinclusions}). Now for $V$ a $B_n$-representation in $\mathsf{M}$, the Peter--Weyl theorem for $\mathsf{M}$ allows us to define the enriched Frobenius characteristic of the mixed Hodge structure $V$ as
\begin{equation}\label{eq:wreathproductschur}
\ch_{2\wr n}^{\mathsf{M}}(V) := \sum_{k = 0}^{n} \sum_{\substack{{\lambda \vdash k}\\{\nu \vdash (n-k)}}} [M_{\lambda, \nu}] \otimes \iota_\pp (s_{\lambda}) \iota_\qq(s_{\nu}) \in K_0(\mathsf{M}) \otimes \Lambda(S_2),
\end{equation}
where the $M_{\lambda, \nu}$ satisfy
 \[ V \cong \bigoplus_{k = 0}^{n} \bigoplus_{\substack{{\lambda \vdash k}\\{\nu \vdash (n-k)}}} M_{\lambda, \nu} \otimes V_{\lambda, \nu}. \] 
For a $B_n$-space $X$, the $B_n$-equivariant Serre characteristic is defined by \[\mathsf{e}^{B_n}(X):= \sum_{i} (-1)^i \ch_{2 \wr n} (H^i_c(X;\QQ)) \in K_0(\mathsf{M}) \otimes \Lambda(S_2). \]
The following formula calculates the $B_n$-equivariant Serre characteristic of induction from a trivial representation.
\begin{prop}\cite[Proposition 3.8]{semiclassicalremark} Let $H < B_n$ be a subgroup. Then
    $$\mathsf{e}^{B_n}\left( \Ind_{H}^{B_n}\Spec(\C)\right) = \frac{1}{|H|} \sum_{h\in H} \Psi(h) \in K_0(\mathsf{M}) \otimes \Lambda(S_2).$$
\end{prop} 

\subsection{$S_2 \times \mathbb{S}$ and $\mathbb{B}$-modules} Similar to $\mathbb{S}$-modules, we define an $S_2\times \mathbb{S}$-module $\mathcal{V}$ (resp. $\mathbb{B}$-module) in $\mathsf{M}$ as a sequence $\mathcal{V}(n)$ of $S_2\times S_n$-representations (resp. $B_n$-representations) in $\mathsf{M}$ for each $n \geq 0$. The Frobenius characteristics of $S_2\times \mathbb{S}$ and $\mathbb{B}$-modules are defined in the same way as for $\mathbb{S}$-modules, using (\ref{eq:S2timesSfrobchar}) and (\ref{eq:wreathproductschur}). They are denoted as $\ch_{2,-}$ and $\ch_{2\wr-}$ respectively.

Let $\mathcal{V}, \mathcal{W}$ be $S_2\times \mathbb{S}$-modules. We set
\[ (\mathcal{V} \boxtimes \mathcal{W})(n) = \bigoplus_{k = 0}^{n} \Ind_{S_k \times S_{n - k}}^{S_n} \mathcal{V}(k) \otimes \mathcal{W}(n - k). \] In particular, the $j$th box power
$\mathcal{W}^{\boxtimes j}$ naturally carries an $B_j$-action. One can thus compose a $\mathbb{B}$-module $\mathcal{V}$ and $S_2\times \mathbb{S}$-module $\mathcal{W}$, as long as either $\mathcal{V}$ has bounded degree or $\mathcal{W}(0) = 0$. It is defined so that $$(\mathcal{V}\circ_{S_2}\mathcal{W})(n):=\bigoplus_{j= 0}^\infty \mathcal{V}(j)\otimes_{B_j} \mathcal{W}^{\boxtimes j}(n).$$ There is a corresponding operation $f \circ_{S_2} g$ for $f \in K_0(\mathsf{M}) \otimes \Lambda(S_2)$ and $g \in K_0(\mathsf{M}) \otimes \Lambda_2 \otimes \Lambda$, whenever $f$ has bounded degree or $g \in K_0(\mathsf{M}) \otimes (\Lambda_2 \oplus (\oplus_{n > 0}  \Lambda_n))$. The operation $\circ_2$ extends the assignment $$(\ch_{2,-}(\mathcal{V}), \ch_{2\wr -}(\mathcal{W}))\mapsto \ch(\mathcal{V}\circ_{S_2}\mathcal{W})$$ for an $S_2\times \mathbb{S}$-module $\mathcal{V}$ and a $\mathbb{B}$-module $\mathcal{W}$, and satisfies analogous properties to the usual plethysm discussed in \S\ref{subsec-sym}: for any $n$, $\pp_n \circ_{S_2} -$ and $\qq_n \circ_{S_2} -$ are algebra homomorphisms.
\begin{lem}\label{lem:frakpleth1} For $f\in K_0(\mathsf{M}) \otimes\Lambda_{2}\otimes \Lambda$, we have the following identities:
    \begin{enumerate}
        \item $\pp_n\circ_{S_2}f = p_n\circ (\pp_1\circ_{S_2} f)$, and $\qq_n\circ_{S_2}f = p_n\circ (\qq_1\circ_{S_2} f)$,
        \item $\pp_n^j\circ_{S_2}f = (\pp_n\circ_{S_2}f)^j,$ $\qq_n^j\circ_{S_2}f = (\qq_n\circ_{S_2}f)^j$.
    \end{enumerate}
\end{lem}
\begin{lem}\label{lem:frakpleth2}
    Writing $f\in K_0(\mathsf{M}) \otimes\Lambda_{2}\otimes \Lambda$ as $f = \frac{1}{2}p_1^2\otimes f_1 + \frac{1}{2}p_2\otimes f_2$ for unique $f_1,f_2\in K_0(\mathsf{M})\otimes \Lambda$, we have $$\mathfrak{p}_1\circ_{S_2} f=f_1,$$ $$\mathfrak{q}_1\circ_{S_2} f=f_2.$$
\end{lem}
\begin{proof}
    We decompose $f$ into homogeneous graded pieces and assume $f\in \Lambda_2\otimes \Lambda_n$ for some $n\geq 0$. Because $\Lambda_n$ is spanned by $\ch_n(V)$ of $S_n$-representations $V$, it suffices to prove the formula for $f=\ch_{2,n}(V)$ where $V$ is some $S_2\times S_n$-representation.

    Following the proof of Lemma \ref{lem-sn+2}, we may decompose $V$ uniquely into $$V= \left(\mathrm{Triv}_{S_2}\otimes V_1\right)\oplus \left(\mathrm{Alt}_{S_2}\otimes V_2\right), $$ where $V_1 = \mathrm{Hom}_{S_2}(\mathrm{Triv}_{S_2},V)$ and $V_2 = \mathrm{Hom}_{S_2}(\mathrm{Alt}_{S_2},V)$ are $S_n$-representations. 
    
    Applying $\ch_{2,n}$ on both sides, we decompose $\ch_{2,n}(V)\in \Lambda_2\otimes \Lambda_n$ into: \begin{align}\
        \ch_{2,n}(V) & = \frac{p_1^2+p_2}{2}\otimes \ch_n(V_1)+ \frac{p_1^2-p_2}{2}\otimes \ch_n(V_2)\\ & = \frac{p_1^2}{2}\otimes (\ch_n(V_1)+\ch_n(V_2)) + \frac{p_2}{2}\otimes (\ch_n(V_1)-\ch_n(V_2)).\label{eq-decomp}
    \end{align} 

    In the following, we identify $B_1\cong S_2$ . Applying the definition of plethysm given above, we find that $$ \mathrm{Triv}_{S_2} \circ_{S_2} V = \mathrm{Triv}_{S_2}\otimes_{S_2} V = \mathrm{Hom}_{S_2}(\mathrm{Triv}_{S_2},V)=V_1,$$ 
    and
    $$\mathrm{Alt}_{S_2} \circ_{S_{2}} V = \mathrm{Alt}_{S_2}\otimes_{S_2}V  = \mathrm{Hom}_{S_2}(\mathrm{Alt}_{S_2},V)=V_2.$$ Applying Frobenius characteristic on both sides and computing $\ch_{2\wr 1}(\mathrm{Triv}_{S_2})$, $\ch_{2\wr 1}(\mathrm{Alt}_{S_2})$ from the definitions, we get $$\frac{1}{2}(\pp_1+\qq_1)\circ_{S_2} \ch_{2,n}(V) = \ch_n (V_1),\quad \frac{1}{2}(\pp_1-\qq_1) \circ_{S_2} \ch_{2, n}(V) = \ch_n (V_2).$$ Adding and subtracting the two equations, we have $$\pp_1 \circ_{S_2} \ch_{2,n}(V)= \ch_n(V_1)+\ch_n(V_2),\quad \qq_1 \circ_{S_2} \ch_{2, n}(V) = \ch_n(V_1)-\ch_n(V_2).$$ This agrees with the decomposition of $\ch_{2,n}(V)$ along $\frac{p_1^2}{2}\otimes \Lambda_{n}\oplus \frac{p_2}{2}\otimes \Lambda_n$ from (\ref{eq-decomp}) and hence proves the formula. 
\end{proof}

\subsection{Plethysm in the graded setting}\label{sec-GradedPlethysm}
All of the above constructions generalise to graded $\mathbb{S}$-, $(S_2 \times \mathbb{S})$-, and $\mathbb{B}$-modules, whose Grothendieck rings are isomorphic to $K_0(\mathsf{M)} \otimes R[[q]]$, where $R = \Lambda$, $\Lambda_2 \otimes \Lambda$, and $\Lambda(S_2)$, respectively. In these settings, we extend plethysm by defining $p_n \circ q = q^n$.
\section{Composition and localisation on graded $\mathbb{S}$-spaces}\label{section-rt}
In this section, we explain how the plethysm operations described in the previous section can be upgraded to \textit{compositions} of (graded) $\mathbb{S}$-spaces. In the context of moduli spaces of marked nodal curves, the composition operation has the appealing geometric interpretation as attaching along marked points. This observation allows us to translate the geometric descriptions of strata in $\Mbar_{1,n}(\mathbb{P}^r,d)$ into a plethystic formula for $\mathsf{b}_{1,r}$. This is the content of Proposition \ref{prop-cut}.

The composition operation on graded $\mathbb{S}$-spaces has been explored in \cite{GetzlerPandharipande}, and we combine their formalism with torus localisation. In Lemma \ref{lem-loc}, we prove the (well-known) fact that the Frobenius characteristic of a graded $\mathbb{S}$-space with $\C^\star$-action can be recovered from the fixed locus. The only requirement is for the $\mathbb{C}^\star$-action to commute with the $S_n$-action; no smoothness or compactness property is needed.


\subsection{$\mathbb{S}$-, $\mathbb{B}$-, and $S_2 \times \fS$-spaces} An $\mathbb{S}$-space $\mathcal{X}$ is a sequence $\mathcal{X}(n)$ of $S_n$-spaces, for each $n \geq 0$. $\mathbb{B}$-spaces and $S_2 \times \fS$ spaces are defined analogously, replacing $S_n$ with $B_n$ and $S_2 \times S_n$, respectively. For each $\mathbb{A} \in \{\mathbb{S}, \mathbb{B}, S_2 \times \mathbb{S} \}$, the Grothendieck groups $K_0(\mathsf{V}, \mathbb{A})$ of $\mathbb{A}$-spaces can be made into a ring using the $\boxtimes$ product:
\[(\mathcal{X}\boxtimes \mathcal{Y})(n):= \coprod_{j = 0}^{n} \Ind_{G_{j}\times G_{n-j}}^{G_n} \mathcal{X}(j)\times \mathcal{Y}(n-j),\] where $G_k = S_k$ for all $k$ if $\mathbb{A} \in\{ \mathbb{S}, S_2 \times \mathbb{S}\}$ or $G_k = B_k$ for all $k$ if $\mathbb{A} = \mathbb{B}$. In the case of $\mathbb{S}$-spaces $\mathcal{X}$ and $\mathcal{Y}$ which are moduli spaces of marked objects (e.g. pointed curves) the operation $\Ind_{S_{j}\times S_{n-j}}^{S_n}$ should be seen as taking all ways of distributing the $n$ marked points to the $j$ markings coming from $\mathcal{X}(j)$ and the remaining $n-j$ ones from $\mathcal{Y}(n -j)$.

We make use of the \emph{composition} of $\mathbb{S}$-spaces. Notice that for an $\mathbb{S}$-space $\mathcal{Y}$ and $j\geq 0$, the product $\mathcal{Y}^{\boxtimes j}$ admits an $S_j$-action (not to be confused with its $\mathbb{S}$-space structure) of coordinate permutations. The composition operation is defined as $$(\mathcal{X}\circ \mathcal{Y})(n):= \coprod_{j = 0}^{\infty} \left(\mathcal{X}(j)\times \mathcal{Y}^{\boxtimes j}(n) \right)/S_j,$$
for $\mathcal{X}$ any $\fS$-space and $\mathcal{Y}$ an $\fS$-space with $\mathcal{Y}(0) = \varnothing$; the product $\mathcal{X}(j)\times \mathcal{Y}^{\boxtimes j}(n)$ receives the diagonal $S_j$-action. If $\mathcal{X}, \mathcal{Y}$ are moduli spaces of marked objects, this operation should be visualised as attaching `tails' coming from $\mathcal{Y}$ to $\mathcal{X}$, which is seen as the `core' component. Similarly, if $\mathcal{X}$ is a $\mathbb{B}$-space and $\mathcal{Y}$ is a $S_2 \times \mathbb{S}$-space, we define a $\mathbb{S}$-space $\mathcal{X} \circ_{S_2} \mathcal{Y}$ by
\[(\mathcal{X} \circ_{S_2} \mathcal{Y})(n) = \coprod_{k \geq 0} (\mathcal{X}(k) \times \mathcal{Y}^{\boxtimes k}(n))/B_k, \]
where $B_k$ acts diagonally. For $[\mathcal{X}] \in K_0(\mathsf{V}, \mathbb{A})$ and $\mathbb{A} \in \{\mathbb{S}, \mathbb{B}, S_2 \times \mathbb{S}\}$, we write \[\mathsf{e}(\mathcal{X}) = \sum_{n \geq 0} \mathsf{e}^{G_n}(\mathcal{X}(n))\] for the \textit{Serre characteristic} of the $\mathbb{A}$-space $\mathcal{X}$; again $G_k$ is the appropriate group among $S_n$, $B_n$, or $S_2 \times S_n$, and $\mathsf{e}$ takes values in $K_0(\mathsf{M}) \otimes R$ where $R = \Lambda$, $R = \Lambda(S_2)$, or $R = \Lambda_2 \otimes \Lambda$, respectively. In each case, $\mathsf{e}$ is a ring homomorphism, which respects the composition structures.

\subsection{Graded $\mathbb{S}$-spaces and $\mathbb{B}$-spaces}
A \textit{graded} $\mathbb{A}$-space $\mathcal{X}(n, d)$ for $\mathbb{A} = \mathbb{S}$ or $\mathbb{A} = \mathbb{B}$ is a
$G_n$-space for each integer $d \geq 0$, where $G_n = S_n$ or $G_n = B_n$, respectively. The Grothendieck ring of graded $\mathbb{A}$-spaces is naturally isomorphic as a ring to $K_0(\mathsf{V},\mathbb{A})[[q]]$ via
\[[\mathcal{X}] \mapsto \sum_{d \geq 0}q^d\sum_{n \geq 0} [\mathcal{X}(n, d)]\in K_0(\mathsf{V}, \mathbb{A})[[q]].\]
 Put differently, the product operation is defined by: \[(\mathcal{X}\boxtimes \mathcal{Y})(n,d):= \coprod_{j = 0}^{n} \coprod_{i = 0}^{d} \Ind_{G_{j}\times G_{n-j}}^{G_n} \mathcal{X}(j,i)\times \mathcal{Y}(n-j,d-i).\]
The composition operation on $K_0(\mathsf{V}, \fS)[[q]]$ is defined by:
 \[(\mathcal{X}\circ \mathcal{Y})(n,d):= \coprod_{i=0}^d\coprod_{j = 0}^{\infty} \left(\mathcal{X}(j,i)\times \mathcal{Y}^{\boxtimes j}(n,d-i) \right)/S_j.\]

The Serre characteristic $\mathsf{e}(\mathcal{X})$ of a graded $\mathbb{A}$-space $\mathcal{X}$ is defined as in (\ref{eq:SerreChar})
and the Frobenius characteristic $\mathsf{ch}(\mathcal{X})$
is the image of $\mathsf{e}(X)$ after taking ranks.

For graded $\mathbb{S}$-spaces $\mathcal{X}$ and $\mathcal{Y}$, we have $$\mathsf{e}(\mathcal{X}\circ \mathcal{Y}) = \mathsf{e}(\mathcal{X})\circ \mathsf{e}(\mathcal{Y}),$$ where the plethysm on the right is as defined in \S\ref{sec-GradedPlethysm}. For a graded ${\mathbb{B}}$-space $\mathcal{X}$ and an $S_2 \times \fS$-space $\mathcal{Y}$, we naturally get a graded $\fS$-space $(\mathcal{X} \circ_2 \mathcal{Y})$ by performing the $S_2$-composition degree-by-degree:
\[(\mathcal{X} \circ_2 \mathcal{Y})(n, d) = (\mathcal{X}(-,d) \circ_2 \mathcal{Y})(n). \]
With this definition, we get that
\[\mathsf{e}(\mathcal{X}) \circ_{S_2} \mathsf{e}(\mathcal{Y}) = \mathsf{e}(\mathcal{X}\circ_{S_2}\mathcal{Y}), \]
for a graded $\mathbb{B}$-space $\mathcal{X}$ and an $S_2 \times \fS$-space $\mathcal{Y}$, where the plethysm on the right is as defined in \S\ref{sec-GradedPlethysm}.

The following elementary result will be applied to calculate the graded plethysm in later sections.
\begin{lem}\cite[§4.14]{tomDieck}\label{lem-frob}
    Let $G$ be a finite group and $H\subset G$ a subgroup. For an $H$-space $X$, there is an isomorphism $X/H\cong (\Ind_{H}^G X)/G$.
\end{lem}

\subsection{Localisation for graded $\mathbb{S}$-spaces}
Let $\mathcal{X}$ be a graded $\mathbb{S}$-space with $\mathbb{C}^\star$-action on each $\mathcal{X}(n,d)$ that commutes with the $S_n$-action: in other words, the $\mathbb{C}^\star$-action is given by morphisms of graded $\mathbb{S}$-spaces. Then we define a graded $\mathbb{S}$-space $\mathcal{X}^{\C^\star}$ by
\[ \mathcal{X}^{\C^\star}(n, d) := \mathcal{X}(n, d)^{\C^\star}. \]

\begin{lem}\label{lem-loc}
With the same assumption as above,
\[ \mathsf{ch}(\mathcal{X} ) = \mathsf{ch}(\mathcal{X}^{\C^\star})\in \Lambda[[q]].  \]
\end{lem}
\begin{proof}
Working with each bidegree $(n,d)$, it suffices to show that if $X$ is a variety with commuting actions of $\C^\star$ and $S_n$, then $\chi^{S_n}(X) = \chi^{S_n}(X^{\C^\star})$. By the Lefschetz fixed-point theorem, we have
\[ \chi^{S_n}(X) = \frac{1}{n!}\sum_{\sigma \in S_n} \chi(X^{\sigma}) \prod_{i = 1}^{n} p_i^{\lambda_i(\sigma)} \]
where $\lambda_i(\sigma)$ is the number of $i$-cycles in $\sigma$. Since the $S_n$-action commutes with the $\C^\star$-action, we have
\[ (X^{\sigma})^{\C^{\star}} = (X^{\C^{\star}})^{\sigma} \]
for all $\sigma \in S_n$. Therefore, to prove the lemma, it suffices to prove that for all $\sigma\in S_n,$ 
\[ \chi(X^{\sigma} ) = \chi\left((X^{\sigma})^{\C^{\star}} \right). \]
To show this, we will prove that for an arbitrary $\C^\star$-variety $Y$, we have $\chi(Y) = \chi(Y^{\C^\star})$, which is likely well-known to experts. Let $Y^{sing} \subset Y$ denote the singular locus of $Y$, and let $Y^{sm} = Y \smallsetminus Y^{sing}$ be the smooth locus of $Y^{sm},$ which is open in $Y.$ Since $\C^\star$ must preserve $Y^{sing}$, both $Y^{sing}$ and $Y^{sm}$ are $\C^\star$-invariant, and we have a partition
\[ Y^{\C^\star} = (Y^{sm})^{\C^\star} \sqcup (Y^{sing})^{\C^\star}. \]
By induction on dimension, we can assume that $\chi((Y^{sing})^{\C^\star}) = \chi(Y^{sing})$, so we are reduced to showing that $\chi(Y^{sm}) = \chi((Y^{sm})^{\C^\star})$. Since $\C^\star$ is connected, each $\C^\star$-orbit is irreducible and is contained in some irreducible component of $Y^{sm}.$ Therefore, the $\C^\star$-action preserves each irreducible component of $Y^{sm},$ so we reduce to the case where $Y^{sm}$ is irreducible. 

In this case, we can find a $\C^\star$-invariant open subset $U \subset Y$ \cite{Sumihiro}. By passing further to an open subset if necessary, we may assume that $\C^\star$-action on $U$ is almost free, so that the quotient $U \to U/\C^\star$ exists, and in particular is \'etale locally a $\C^\star$-bundle. It follows that $\chi(U) = 0$. If $(Y^{sm} \smallsetminus U)^{\C^\star} = Y^{sm} \smallsetminus U$, then $(Y^{sm})^{\C^\star} = Y^{sm} \smallsetminus U$ and $\chi(Y^{sm}) = \chi((Y^{sm})^{\C^\star})$. Otherwise, we can repeat this procedure by replacing $Y^{sm}$ with $Y^{sm} \smallsetminus U$. The proof is complete by Noetherian induction.
\end{proof}
\subsection{Separating rational tails}
We now show how to separate rational tails in the calculation of the class of the graded $\mathbb{S}$-space $\Mbar_{g}(r)$. The first step in our calculation is to separate stable maps into:\begin{itemize}
    \item maps from curves which do not contain any proper subcurve of arithmetic genus $g$, which we call the \emph{core} of the stable map;
    \item maps from nodal rational curves glued onto the core along a single point; we call these the \emph{rational tails} of the stable map.
\end{itemize}
This operation allows us to single out the contribution from the first class of maps, which we calculate in the rest of the work. The precise definitions are given as follows.\begin{defn}
    Let $[f: C\to \mathbb{P}^r]\in \Mbar_{g,n}(\mathbb{P}^r,d)$ be a genus-$g$, $n$-pointed stable map. Let $\mathbf{G}$ be the dual graph of $f$; this is the dual graph of the curve $C$, decorated additionally by the degree of $f$ along each component. We refer to such graphs as  \textit{$(g,n,d)$-graphs} in \cite[Definition 1.1]{vzdualcomplex}.

    The \textit{core} of $\mathbf{G}$, denoted as $\mathbf{G}^{\mathrm{core}}$ is the minimal genus-$g$ subgraph of $\mathbf{G}$, together with a degree decoration restricted from $\mathbf{G}$, as well as one extra leg for each edge that connects $\mathbf{G}^{\mathrm{core}}$ to $\mathbf{G}\smallsetminus \mathbf{G}^{\mathrm{core}}$. A \textit{rational tail} of $\mathbf{G}$ is a connected component of $\mathbf{G} \smallsetminus \mathbf{G}^{\mathrm{core}}$, where $\mathbf{G}$ and $\mathbf{G}^{\mathrm{core}}$ are topologized as $1$-dimensional CW-complexes, ignoring legs and vertex decorations. Each rational tail carries degree decorations restricted from $\mathbf{G}$, and we endow each one with one extra leg for the edge in $\mathbf{G}$ that connects it to $\mathbf{G}^{\mathrm{core}}$.
\end{defn}

\begin{exmp}
    The following figure is a stable map dual graph $\mathbf{G}$ in $\Mbar_{1,5}(\mathbb{P}^r,7)$ and its core $\mathbf{G}^{\mathrm{core}}$ as well as rational tails $\mathbf{T}_1,\dots, \mathbf{T}_4$. Note that in labeling the rational tails, we have made a non-canonical choice of ordering the tails.

    \begin{figure}[H]
        \centering
        \includegraphics[width=0.7\linewidth]{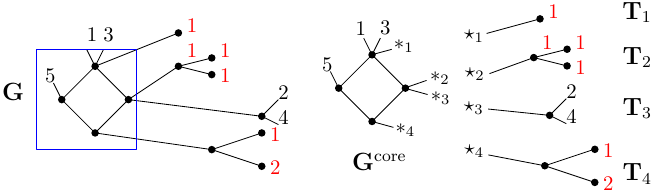}
        \caption{The marking assignment is in black, and the degree assignment in red. Vertices without red labels have degree zero.}
        \label{fig:cut4tails}
    \end{figure}
\end{exmp}

The combinatorial definitions regarding the core and rational tails specify the following locally closed subscheme in the moduli space.

\begin{defn}
    We denote by $\Mbar_{g,n}^{\mathrm{nrt}}(\mathbb{P}^r,d)\subset \Mbar_{g,n}(\mathbb{P}^r,d)$ the locally closed subscheme consisting of stable maps where the dual graph does not contain any rational tails. Equivalently, the dual graph of the stable map is equal to its own core.
\end{defn}

Before stating the main formula of the section, we see how the above definitions allow us to express each stratum in $\Mbar_{g,n}(\mathbb{P}^r,d)$ as a fibre product. Let $\mathbf{G}$ be an $(g,n,d)$-graph, so that it specifies a locally closed stratum $\mathcal{M}_{\mathbf{G}}\subset \Mbar_{g,n}(\mathbb{P}^r,d)$. Let $\mathbf{G}^{\mathrm{core}}$ be the core, and let $\mathbf{T}_1,\dots, \mathbf{T}_k$ be the rational tails of $\mathbf{G}$. The graphs $\mathbf{G}^{\mathrm{core}}$ and $\mathbf{T}_1,\dots, \mathbf{T}_k$ specify strata $$\mathcal{M}_{\mathbf{G}^{\mathrm{core}}}\subset \Mbar_{g,n'}(\mathbb{P}^r,d'),\mathcal{M}_{\mathbf{T}_i}\subset \Mbar_{0,n_i}(\mathbb{P}^r,d_i)$$ for $i=1,\dots, k$. The degrees $d',d_i$ and valences $n', n_i$ in general differ from $d$ and $n$ and can be determined from $\mathbf{G}$. Let $\star_i$ denote the extra leg in $\mathbf{T}_i$, and let $*_i$ denote the extra legs in $\mathbf{G}^{\mathrm{core}}$ connecting to $\mathbf{T}_i$. These legs define evaluation maps $$\mathrm{ev}_{\star_i}: \mathcal{M}_{\mathbf{T}_i}\to \mathbb{P}^r,\mathrm{ev}_{*_i}: \mathcal{M}_{\mathbf{G}^{\mathrm{core}}}\to \mathbb{P}^r$$ that record the image of the marked points associated to the legs.

With notation as earlier, consider the quotient space $\hat{\mathbf{G}}:=\mathbf{G}/\mathbf{G}^{\mathrm{core}}$ obtained by contracting the core to a point. $\hat{\mathbf{G}}$ is a rooted tree with trees $\mathbf{T}_1,\dots, \mathbf{T}_k$ glued along the root. We assign the central vertex degree $d'$ and restrict the degree and marked point assignment from $\mathbf{G}$ to $\hat{\mathbf{G}}$. 

Define $\mathrm{Aut}(\hat{\mathbf{G}})$ to be the group of rooted tree automorphisms of $\hat{\mathbf{G}}$ that respect the decorations. It follows from the definitions that $\mathrm{Aut}(\hat{\mathbf{G}})\subset S_k$ is the Young subgroup that permutes blocks of subsets $\{\mathbf{T}_{i_1}, \dots ,\mathbf{T}_{i_j}\}$ where $\mathbf{T}_{i_1}, \dots ,\mathbf{T}_{i_j}$ are isomorphic, unmarked decorated trees. $\mathrm{Aut}(\hat{\mathbf{G}})$ thus naturally acts on $\mathcal{M}_{\mathbf{G}^{\mathrm{core}}}$ by permuting relevant markings and acts on $\prod_{i=1}^k \mathcal{M}_{\mathbf{T}_i}$ by permuting the entries. After describing the group action, we state the description of the stratum $\mathcal{M}_{\mathbf{G}}$ as follows.

\begin{lem}\label{lem-tilfib}
The stratum $\mathcal{M}_{\mathbf{G}}$ is isomorphic to the quotient $\widetilde{\mathcal{M}}_{\mathbf{G}}/\mathrm{Aut}(\hat{\mathbf{G}})$, where $\widetilde{\mathcal{M}}_{\mathbf{G}}$ is the fibre product \begin{equation*} 
    \begin{tikzcd} 
     \widetilde{\mathcal{M}}_{\mathbf{G}} \ar[d] \ar[dr, phantom, "\square"] \ar[r] & \prod_{i=1}^k \mathcal{M}_{\mathbf{T}_i} \ar[d, "\prod_{i} \mathrm{ev}_{\star_i}"]  \\
     \mathcal{M}_{\mathbf{G}^{\mathrm{core}}} \ar[r, "\prod_{i}\mathrm{ev}_{*_i}"] & \prod_{i=1}^k \mathbb{P}^r
    \end{tikzcd}
\end{equation*}
\end{lem}

There is an alternative description of $\widetilde{\mathcal{M}}_{\mathbf{G}}$ that will be useful. We note that $\mathrm{PGL}_{r+1}$ acts transitively on the target $\mathbb{P}^r$ and hence has an induced action on spaces of maps $\mathcal{M}_{\mathbf{T}_i}$, so that the evaluation maps $\mathrm{ev}_{\star_i}: \mathcal{M}_{\mathbf{T}_i}\to \mathbb{P}^r$ are $\mathrm{PGL}_{r+1}$-equivariant. The $\mathrm{PGL}_{r+1}$-equivariance implies that $\mathrm{ev}_{\star_i}$ are Zariski locally trivial fibrations. We define $\mathcal{M}^\star_{\mathbf{T}_i} := \mathrm{ev}_{\star_i}^{-1}(p)$ to be the fibre over any $p\in \mathbb{P}^r$.

\begin{lem}\label{lem-*fib}
The map $\widetilde{\mathcal{M}}_{\mathbf{G}}\to \mathcal{M}_{\mathbf{G}^{\mathrm{core}}}$ is Zariski locally trivial with fibres isomorphic to $\prod_{i=1}^k \mathcal{M}^\star_{\mathbf{T}_i}$.
\end{lem}
\begin{proof}
    In the previous diagram, the vertical map $\prod_i \mathrm{ev}_{\star i}$ on the right is a Zariski locally trivial fibration with fibres $\prod_{i=1}^k \mathcal{M}^\star_{\mathbf{T}_i}$. Pulling back some Zariski local trivialisation of $\prod_i \mathrm{ev}_{\star i}$ along the fibre diagram trivialises the vertical map $\widetilde{\mathcal{M}}_{\mathbf{G}}\to \mathcal{M}_{\mathbf{G}^{\mathrm{core}}}$ on the left with the same fibres.
\end{proof}


We recall the following notation, used in the introduction, to state the formulas on $S_n$-equivariant Euler characteristics.

\begin{defn}
    Let $\Mbar_g(r)$ be the graded $\mathbb{S}$-space given by $(n,d)\mapsto \Mbar_{g,n}(\mathbb{P}^r,d)$, and let $\Mbar^{\mathrm{nrt}}_g(r)$ be the graded $\mathbb{S}$-space given by $(n,d)\mapsto \Mbar_{g,n}^{{\mathrm{nrt}}}(\mathbb{P}^r,d)$. As earlier, the Frobenius characteristics of these spaces are denoted as $\mathsf{b}_{g,r}$ and  $\mathsf{b}_{g,r}^{\mathrm{nrt}}$, respectively.
\end{defn}
Consider the evaluation map $\mathrm{ev}_{n+1}: \Mbar_{0,n+1}(\mathbb{P}^r,d)\to \mathbb{P}^r$. By the same reasoning as earlier, $\mathrm{ev}_{n+1}$ is $\mathrm{PGL}_{r+1}$-equivariant and hence a Zariski locally trivial fibration. Let \[\mathrm{ev}_{n+1}^{-1}(p)=:\Mbar^\star_{0,n}(\mathbb{P}^r,d) \subset \Mbar_{0,n+1}(\mathbb{P}^r,d)\] be a fibre of the map, which is an $S_n$-space. Define $\Mbar^\star_{0}(r)$ as the graded $\mathbb{S}$-space given by $$(n,d)\mapsto \Mbar_{0,n}^\star(\mathbb{P}^r,d).$$ Note that by Zariski local triviality, we have
    \[ [\Mbar^\star_0(r)] = \frac{D\Mbar_{0}(r)}{[\P^r]}, \] in $K_0(\mathsf{V}, \mathbb{S})[[q]]$. From this it follows that
    \[ \mathsf{ch}(\Mbar^\star_0(r)) = \frac{1}{r + 1} \mathsf{b}_{0, r}'. \]

The following is the main result of the section. It states how the Serre characteristic of $\Mbar_{1}(r)$ is determined by those of the no-rational-tails locus $\Mbar_{1}^{\mathrm{nrt}}(r)$ and the space $\Mbar^\star_0(r)$ of maps from rational tails.
\begin{prop}\label{prop-cut}
    There is an equality in $K_0(\mathsf{V}, \mathbb{S})[[q]]$: $$\Mbar_{1}(r)= \Mbar_{1}^{\mathrm{nrt}}(r)\circ \left([\Spec \C] + \frac{D\Mbar_{0}(r)}{[\P^r]}\right).$$ In particular, $$\mathsf{b}_{1,r} = \mathsf{b}_{1,r}^{\mathrm{nrt}}\circ \left(p_1+\frac{\mathsf{b}'_{0,r}}{r+1}\right)\in \Lambda[[q]].$$
\end{prop}
\begin{rem}
    The term $[\Spec \C]$ in the formula has bidegree $n=1, d=0$ and is interpreted as providing legs on the core. In the following proof, it will be convenient to interpret the term by formally setting $\Mbar_{0}^\star(r)(1,0):=\Spec \C$, whereas $\Mbar_{0}^\star(r)(0,0)=\varnothing$ because of the stability condition.
\end{rem}
\begin{proof} 
    Unwrapping the definition of plethysms of graded $\mathbb{S}$-spaces, the assignment of $(n,d)$ on the right hand side is equal to $$\coprod_{\delta=0}^d \coprod_{j=0}^\infty \left[ \Mbar_{1,j}^{\mathrm{nrt}}(\mathbb{P}^r,\delta)\times \left(\Mbar_{0}^{*}(r)\right)^{\boxtimes j}(n, d-\delta)\right]/S_j.$$ 
    
    Denote the item in the above formula indexed by $\delta,j$ as $X_{\delta,j}$. We claim that $X_{\delta,j}$ has the same class as the locally closed stratum $\Mbar^{(\delta,j)}_{1,n}(\mathbb{P}^r,d)$ in $\Mbar_{1,n}(\mathbb{P}^r,d)$ consisting of stable maps where the core has degree $\delta$, and the core has valency $j$: in other words, the number of rational tails and the number of legs on the core sum to $j$. 
    
    To prove the claim, we further stratify $\Mbar^{(\delta,j)}_{1,n}(\mathbb{P}^r,d)$ by the unordered partitions $\lambda \vdash n$ and $\mu\vdash (d-\delta)$ into $j$ non-negative parts, which record the distribution of the number of marked points and degrees onto each leg or rational tail. Let $\Mbar^{(\lambda, \mu)}_{1,n}(\mathbb{P}^r,d)$ be such a refined stratum. From Lemmas \ref{lem-tilfib}, \ref{lem-*fib}, we have that $$\Mbar^{(\lambda, \mu)}_{1,n}(\mathbb{P}^r,d) = \Ind^{S_n}_{S_{\lambda}}\left[\Mbar^{\mathrm{nrt}}_{1,j}(\mathbb{P}^r,\delta)\times \prod_{i=1}^j \Mbar^\star_{0}(r)(\lambda_i, \mu_i)\right]/\mathrm{Aut}(\lambda, \mu),$$where $\mathrm{Aut}(\lambda, \mu)\subset S_j$ is the subgroup which permutes the maximal blocks $I\subset \{1,\dots,j\}$ such that for all $i\in I$, $\mu_i$ are equal and $\lambda_i=0$. Geometrically, $\mathrm{Aut}(\lambda, \mu)$ permutes the unmarked rational tails of the same degree.

    Returning to the terms $X_{\delta,j}$ in the plethystic formula, we notice that it can be partitioned along pairs of unordered partitions $(\lambda, \mu)$ as  $$X_{\delta,j}=\coprod_{(\lambda\vdash n, \mu\vdash n-j)}\left[\Mbar_{1,j}^{\mathrm{nrt}}(\mathbb{P}^r,\delta)\times Y_{(\lambda, \mu)}\right]/S_j,$$ where each $Y_{(\lambda, \mu)}$ is a union of $j$-fold products in $\left(\mathsf{pt}_{S_1}+\Mbar_{0}^{*}(r)\right)^{\boxtimes j}(n, d-\delta)$ such that the unordered partitions on marked points and degrees recover $\lambda, \mu$. More explicitly, consider the set $L_{(\lambda, \mu)}$ of pairs of ordered partitions that forget to $(\lambda, \mu)$, then $Y_{(\lambda, \mu)}$ is a disjoint union $$Y_{(\lambda, \mu)} = \coprod_{(\tilde{\lambda}, \tilde{\mu})\in L_{(\lambda, \mu)}} \prod_{i=1}^j \Mbar_{0}^\star(r)(\lambda_i, \mu_i).$$  The $S_j$-action on $Y_{(\lambda, \mu)}$ permutes terms indexed by $(\tilde{\lambda}, \tilde{\mu})\in L_{(\lambda, \mu)}$ via composition. The stabilisers of $L_{(\lambda, \mu)}$ under the $S_j$-action are uniformly $\mathrm{Aut}(\lambda, \mu)$. Therefore, the same holds for $Y_{(\lambda, \mu)}$.
    
    Applying Lemma \ref{lem-frob} to the $S_j$-action on $\Mbar^{\mathrm{nrt}}_{1,j}(\mathbb{P}^r,\delta)\times Y_{(\lambda, \mu)}$, we hence identify $$\left[\Mbar_{1,j}^{\mathrm{nrt}}(\mathbb{P}^r,\delta)\times Y_{(\lambda, \mu)}\right]/S_j\cong \Mbar^{(\lambda, \mu)}_{1,n}(\mathbb{P}^r,d).$$ Notice that all the equalities and isomorphisms above are $S_n$-equivariant. Summing over all partitions $(\lambda, \mu)$ and all $(\delta,j)$, we get the desired formula in the equivariant Grothendieck ring.
\end{proof}

\section{Localisation on the Kontsevich space}\label{sec-localisation}
With the formula in Proposition \ref{prop-cut} in hand, we turn to the calculation of the characteristics of maps with no rational tails, namely that of 
\[ \mathsf{b}_{1, r}^{\mathrm{nrt}}:= \mathsf{ch}(\Mbar_{1}^{\mathrm{nrt}}(r) ). \]

Pick a generic $\mathbb{C}^\star$-action on $\mathbb{P}^r$ so that its fixed points are $(r+1)$ isolated points. This induces an action of $\C^\star$ on the $\mathbb{S}$-space $\Mbar_{g}(r)$, which preserves the subspace $\Mbar_{g}^{\mathrm{nrt}}(r)$. Lemma \ref{lem-loc} gives 
\[ \sfch\left(\Mbar_{g}^{\mathrm{nrt}}(r)\right) =  \sfch\left(\Mbar_{g}^{\mathrm{nrt}}(r)^{\C^\star}\right) .\]
Thus we focus our attention on computing
\[\mathscr{B}_r := \mathsf{e}(\Mbar_{1}^{\mathrm{nrt}}(r)). \]
In this section, we introduce a combinatorial graded $\mathbb{B}$-space $\mathsf{Dih}_r$ and an $(S_2 \times \mathbb{S})$-space $\mathsf{Cat}$ of chains of rational curves, and prove the following theorem for the motive $[\Mbar_{1}^{\mathrm{nrt}}(r)] \in K_0(\mathsf{V}, \fS)[[q]]$. 
\begin{thm}\label{thm-locnrt} For an integer $r \geq 1$, we have
\[[\Mbar_{1}^{\mathrm{nrt}}(r)^{\C^\star}] = [({\P^r})^{\C^\star}][\Mbar_{1}^{\mathrm{nrt}}(0)] +  [\mathsf{Dih}_r] \circ_{S_2} [\mathsf{Cat}]. \]
\end{thm}
Graber and Pandharipande \cite[\S 4]{graberpand} have described the components
\[\Mbar_{g}(r)^{\C^\star}(n, d) = \Mbar_{g, n}(\P^r, d)^{\C^\star} \] of the graded $\mathbb{S}$-space $\Mbar_{g}(r)^{\C^\star}$ in terms of certain decorated graphs, as we recall shortly. Theorem \ref{thm-locnrt} is a motivic interpretation of their description when $g = 1$. 
\subsection{Localisation graphs}
Consider the groupoid $\mathsf{LG}_{g,n,d}$ of tuples \[\mathbf{G} = (G, w, m, s),\]
where $G$ is a connected graph with at least one edge, and: \begin{itemize}\item $w: V(G)\to \mathbb{Z}_{\geq 0}$ is a genus assignment, with total genus $h^1(G)+ \sum_{v\in V(G)}w(v)=g$; \item $m: \{1,\dots, n\}\mapsto V(G)$ is a marking/leg assignment;
\item $s : E(G) \to \Z_{> 0}$ is a degree assignment, with $\sum_{e \in E(G)} s(e) = d$.
\end{itemize}
We exclude graphs with loop edges from $\mathsf{LG}_{g, n, d}$, but pairs of parallel edges are allowed. An isomorphism in $\mathsf{LG}_{g,n,d}$ is an isomorphism of graphs that respects the decorations. In particular, isomorphisms of graphs are viewed as bijections on the underlying sets of vertices and half-edges. Note that the positivity requirement on $s$ implies that $\mathsf{LG}_{g, n, d}$ has finitely many isomorphism classes of objects, which are connected graphs with at most $d$ edges. Now, for each integer $r > 0$, we define a groupoid $\mathsf{LG}_{g, n, d}(r)$ of pairs $(\mathbf{G}, c)$ where $\mathbf{G} \in \mathrm{Ob}(\mathsf{LG}_{g, n,d})$ and \[c : V(\mathbf{G}) \to (\P^r)^{\C^\star} \cong \{0, 1, \ldots, r\} \]
is a graph colouring of $\mathbf{G}$ with $r+1$ colours. The isomorphisms in $\mathsf{LG}_{g, n, d}(r)$ are $\mathsf{LG}_{g, n, d}$-isomorphisms that also respect the colouring. For convenience, we set $\mathsf{LG}_{g, n, d}(0):= \mathsf{LG}_{g,n,d}$. 
An object $(\mathbf{G},c)$ of the groupoid $\mathsf{LG}_{g,n,d}(r)$ corresponds to a $\C^\star$-fixed closed subscheme $\Mbar_{\mathbf{G}, c}\subset \Mbar_{g,n}(\mathbb{P}^r,d)$ with
\begin{equation}\label{eq:nonequivariantproduct}
\Mbar_{\mathbf{G}, c} \cong \prod_{v\in V(\mathbf{G})}\Mbar_{w(v), n(v) +  \mathrm{val}(v)}/\mathrm{Aut}_{\mathsf{LG}_{g, n, d}(r)}(\mathbf{G},c),
\end{equation}
where for a vertex $v \in V(\mathbf{G})$ we put $n(v):=|m^{-1}(v)|$, and $\val(v)$ denotes the graph-theoretic valence of $v$. The geometric picture is as follows: each edge $e$ in the graph $\mathbf{G}$ represents a $\C^\star$-fixed ramified covering of degree $s(e)$ from $\mathbb{P}^1$ onto a $\mathbb{C}^\star$-invariant line on the target $\mathbb{P}^r$, between the two points of $(\P^r)^{\C^\star}$ indicated by the graph colouring $c$. These coverings do not carry any moduli. At each vertex of $\mathbf{G}$, arbitrary Deligne--Mumford stable curves can be glued in, which are contracted in the stable map to the points of $(\mathbb{P}^r)^{\C^\star}$ determined by the colouring. Graber and Pandharipande prove that 
\begin{equation}\label{eq:grabpandupshot}
\Mbar_{g, n}(\P^r, d)^{\C^\star} = \coprod_{(\mathbf{G}, c) \in \pi_0(\mathsf{LG}_{g,n,d}(r))} \Mbar_{\mathbf{G}, c},
\end{equation}
where $\pi_0$ denotes a choice of skeleton of a groupoid. We now make both (\ref{eq:nonequivariantproduct}) and (\ref{eq:grabpandupshot}) equivariant with respect to the symmetric group, as follows. First observe that $S_n$ acts on the groupoids $\mathsf{LG}_{g, n, d}(r)$ for all $r \geq 0$. We write 
\[ \widetilde{\mathsf{LG}}_{g, n, d}(r):= \frac{\mathsf{LG}_{g, n, d}(r)}{S_n} \]
for the quotient: concretely, the marking $m$ is coarsened to an assignment of an integral number of marked points to each vertex, and there are additional graph isomorphisms that reflect this. We adopt the calligraphic symbol $\mathcal{G}$ for objects of the category $\widetilde{\mathsf{LG}}_{g, n, d}$. In the following, if $c$ is a $(\P^r)^{\C^\star}$-colouring of $\mathcal{G}$ in $\widetilde{\mathsf{LG}}_{g,n,d}$ we write
\[ \Aut(\mathcal{G}, c) := \Aut_{\widetilde{\mathsf{LG}}_{g, n, d}(r)}(\mathcal{G}, c), \]
\[ \Aut(\mathcal{G}):= \Aut_{\widetilde{\mathsf{LG}}_{g,n,d}}(\mathcal{G}).\]
Then $\Aut(\mathcal{G}, c)$ is a subgroup of $\Aut(\mathcal{G})$. Each $(\mathcal{G}, c)$ in the groupoid $\widetilde{\mathsf{LG}}_{g, n, d}(r)$ defines an $S_n$-equivariant $\C^\star$-fixed subscheme $\Mbar_{\mathcal{G}, c}$ of $\Mbar_{g, n}(\P^r, d)$, leading to $S_n$-equivariant versions of equations (\ref{eq:nonequivariantproduct}) and (\ref{eq:grabpandupshot}):
\begin{equation}\label{eq:equivproduct}
    \Mbar_{\mathcal{G}, c} \cong \left(\Ind_{\prod_{v \in V(\mathcal{G})} S_{n(v)}}^{S_n} \prod_{v \in V(\mathcal{G})} \Mbar_{g(v), n(v) + \val(v)}\right)/\Aut(\mathcal{G}, c),
\end{equation}
\begin{equation}\label{eq:equivgrabpand}
    \Mbar_{g, n}(\P^r, d)^{\C^\star} = \coprod_{(\mathcal{G}, c) \in \pi_0\left(\widetilde{\mathsf{LG}}_{g,n,d}(r)\right)} \Mbar_{\mathcal{G}, c}.
\end{equation}
It is important to note that the action of $\Aut(\mathcal{G}, c)$ in (\ref{eq:equivproduct}) is induced by its action on the half-edges $H(\mathcal{G})$. This action commutes with the action of $S_n$, so there is a residual $S_n$-action on the quotient.


Now we specialise to genus one, in order to prove Theorem \ref{thm-locnrt}. 
\subsection{Localisation graphs for $\Mbar_{1, n}^{\mathrm{nrt}}(\P^r, d)$}\label{subsec: locgraphdefn} In genus one, $\C^\star$-fixed maps without rational tails correspond to localisation graphs which are cycles: for each $r \geq 0$, let
\[\widetilde{\mathsf{LG}}_{1, n, d}^\circ(r) \subset \widetilde{\mathsf{LG}}_{1, n, d}(r) \]
denote the full subgroupoid consisting of those localisation graphs whose underlying graph is a $k$-cycle, for some integer $k \geq 2$. There is a natural decomposition 
\[\widetilde{\mathsf{LG}}_{1, n, d}^\circ(r) = \coprod_{k \geq 2} \widetilde{\mathsf{LG}}_{1, n, d}^{\circ}(r)^{k} \]
into subgroupoids where $\widetilde{\mathsf{LG}}_{1, n, d}^{\circ}(r)^{k}$ consists of the cycles with exactly $k$ edges. A special role is played by the case $n = 0$. We now choose a preferred set of representatives for $\pi_0(\mathsf{LG}_{1, 0, d}(r)^k)$. We define a cycle graph $C_k$ of length $k$ by setting \[V(C_k) := \{1, \ldots, k\}\] for the set of vertices, and  \[H(C_k) := \{ h_1^{+}, h_1^{-}, h_2^{_+}, h_2^{-}, \ldots, h_k^{+}, h_k^{-}\}\] for the set of half-edges. In our convention, each half-edge $h_i^{\pm}$ is incident to $v_i$, and for $i \leq k - 1$, the half-edge $h_i^{+}$ is glued to $h_{i + 1}^{-}$, while $h_k^{+}$ is glued to $h_1^{-}$; this determines the set of edges $E(C_k)$.

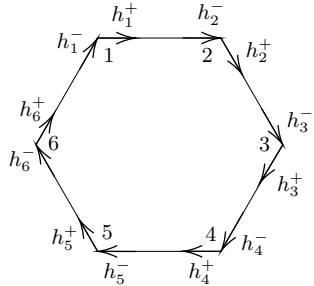
\begin{figure}[H]
    \centering
    \begin{tikzpicture}[x=0.75pt,y=0.75pt,yscale=-1,xscale=1]

\draw   (184.8,98.1) -- (150.35,157.77) -- (81.45,157.77) -- (47,98.1) -- (81.45,38.43) -- (150.35,38.43) -- cycle ;
\draw [line width=0.55]    (173,77.7) -- (183.99,96.48) ;
\draw [shift={(185,98.2)}, rotate = 239.66] [color={rgb, 255:red, 0; green, 0; blue, 0 }  ][line width=0.55]    (10.93,-3.29) .. controls (6.95,-1.4) and (3.31,-0.3) .. (0,0) .. controls (3.31,0.3) and (6.95,1.4) .. (10.93,3.29)   ;
\draw [line width=0.55]    (185,98.2) -- (173.03,117.99) ;
\draw [shift={(172,119.7)}, rotate = 301.16] [color={rgb, 255:red, 0; green, 0; blue, 0 }  ][line width=0.55]    (10.93,-3.29) .. controls (6.95,-1.4) and (3.31,-0.3) .. (0,0) .. controls (3.31,0.3) and (6.95,1.4) .. (10.93,3.29)   ;
\draw [line width=0.55]    (163,135.7) -- (151.34,156.03) ;
\draw [shift={(150.35,157.77)}, rotate = 299.83] [color={rgb, 255:red, 0; green, 0; blue, 0 }  ][line width=0.55]    (10.93,-3.29) .. controls (6.95,-1.4) and (3.31,-0.3) .. (0,0) .. controls (3.31,0.3) and (6.95,1.4) .. (10.93,3.29)   ;
\draw [line width=0.55]    (150.35,157.77) -- (131,157.71) ;
\draw [shift={(129,157.7)}, rotate = 0.18] [color={rgb, 255:red, 0; green, 0; blue, 0 }  ][line width=0.55]    (10.93,-3.29) .. controls (6.95,-1.4) and (3.31,-0.3) .. (0,0) .. controls (3.31,0.3) and (6.95,1.4) .. (10.93,3.29)   ;
\draw [line width=0.55]    (102.8,157.84) -- (83.45,157.78) ;
\draw [shift={(81.45,157.77)}, rotate = 0.18] [color={rgb, 255:red, 0; green, 0; blue, 0 }  ][line width=0.55]    (10.93,-3.29) .. controls (6.95,-1.4) and (3.31,-0.3) .. (0,0) .. controls (3.31,0.3) and (6.95,1.4) .. (10.93,3.29)   ;
\draw [line width=0.55]    (81.45,157.77) -- (72,141.43) ;
\draw [shift={(71,139.7)}, rotate = 59.96] [color={rgb, 255:red, 0; green, 0; blue, 0 }  ][line width=0.55]    (10.93,-3.29) .. controls (6.95,-1.4) and (3.31,-0.3) .. (0,0) .. controls (3.31,0.3) and (6.95,1.4) .. (10.93,3.29)   ;
\draw [line width=0.55]    (57.45,116.17) -- (48,99.83) ;
\draw [shift={(47,98.1)}, rotate = 59.96] [color={rgb, 255:red, 0; green, 0; blue, 0 }  ][line width=0.55]    (10.93,-3.29) .. controls (6.95,-1.4) and (3.31,-0.3) .. (0,0) .. controls (3.31,0.3) and (6.95,1.4) .. (10.93,3.29)   ;
\draw [line width=0.55]    (47,98.1) -- (56,82.44) ;
\draw [shift={(57,80.7)}, rotate = 119.89] [color={rgb, 255:red, 0; green, 0; blue, 0 }  ][line width=0.55]    (10.93,-3.29) .. controls (6.95,-1.4) and (3.31,-0.3) .. (0,0) .. controls (3.31,0.3) and (6.95,1.4) .. (10.93,3.29)   ;
\draw [line width=0.55]    (71.45,55.83) -- (80.45,40.17) ;
\draw [shift={(81.45,38.43)}, rotate = 119.89] [color={rgb, 255:red, 0; green, 0; blue, 0 }  ][line width=0.55]    (10.93,-3.29) .. controls (6.95,-1.4) and (3.31,-0.3) .. (0,0) .. controls (3.31,0.3) and (6.95,1.4) .. (10.93,3.29)   ;
\draw [line width=0.55]    (148.35,38.43) -- (129,38.37) ;
\draw [shift={(150.35,38.43)}, rotate = 180.18] [color={rgb, 255:red, 0; green, 0; blue, 0 }  ][line width=0.55]    (10.93,-3.29) .. controls (6.95,-1.4) and (3.31,-0.3) .. (0,0) .. controls (3.31,0.3) and (6.95,1.4) .. (10.93,3.29)   ;
\draw [line width=0.55]    (100.8,38.49) -- (81.45,38.43) ;
\draw [shift={(102.8,38.5)}, rotate = 180.18] [color={rgb, 255:red, 0; green, 0; blue, 0 }  ][line width=0.55]    (10.93,-3.29) .. controls (6.95,-1.4) and (3.31,-0.3) .. (0,0) .. controls (3.31,0.3) and (6.95,1.4) .. (10.93,3.29)   ;
\draw [line width=0.55]    (150.35,38.43) -- (161.34,57.21) ;
\draw [shift={(162.35,58.93)}, rotate = 239.66] [color={rgb, 255:red, 0; green, 0; blue, 0 }  ][line width=0.55]    (10.93,-3.29) .. controls (6.95,-1.4) and (3.31,-0.3) .. (0,0) .. controls (3.31,0.3) and (6.95,1.4) .. (10.93,3.29)   ;

\draw (83.45,42.83) node [anchor=north west][inner sep=0.75pt]  [font=\footnotesize]  {$1$};
\draw (138.45,42.83) node [anchor=north west][inner sep=0.75pt]  [font=\footnotesize]  {$2$};
\draw (88,17.4) node [anchor=north west][inner sep=0.75pt]  [font=\footnotesize]  {$h_{1}^{+}$};
\draw (170.45,92.83) node [anchor=north west][inner sep=0.75pt]  [font=\footnotesize]  {$3$};
\draw (140.45,142.83) node [anchor=north west][inner sep=0.75pt]  [font=\footnotesize]  {$4$};
\draw (82.45,141.83) node [anchor=north west][inner sep=0.75pt]  [font=\footnotesize]  {$5$};
\draw (52.45,91.83) node [anchor=north west][inner sep=0.75pt]  [font=\footnotesize]  {$6$};
\draw (136,16.4) node [anchor=north west][inner sep=0.75pt]  [font=\footnotesize]  {$h_{2}^{-}$};
\draw (163,37.4) node [anchor=north west][inner sep=0.75pt]  [font=\footnotesize]  {$h_{2}^{+}$};
\draw (186,76.4) node [anchor=north west][inner sep=0.75pt]  [font=\footnotesize]  {$h_{3}^{-}$};
\draw (180.5,112.35) node [anchor=north west][inner sep=0.75pt]  [font=\footnotesize]  {$h_{3}^{+}$};
\draw (160.5,145.35) node [anchor=north west][inner sep=0.75pt]  [font=\footnotesize]  {$h_{4}^{-}$};
\draw (131,161.1) node [anchor=north west][inner sep=0.75pt]  [font=\footnotesize]  {$h_{4}^{+}$};
\draw (83.45,161.17) node [anchor=north west][inner sep=0.75pt]  [font=\footnotesize]  {$h_{5}^{-}$};
\draw (54.45,142.17) node [anchor=north west][inner sep=0.75pt]  [font=\footnotesize]  {$h_{5}^{+}$};
\draw (30.45,98.17) node [anchor=north west][inner sep=0.75pt]  [font=\footnotesize]  {$h_{6}^{-}$};
\draw (36.45,70.17) node [anchor=north west][inner sep=0.75pt]  [font=\footnotesize]  {$h_{6}^{+}$};
\draw (57.45,32.17) node [anchor=north west][inner sep=0.75pt]  [font=\footnotesize]  {$h_{1}^{-}$};

\end{tikzpicture}
    \caption{A six-cycle with vertex and half-edge decorations. We use $i$ in place of $v_i$ in the figure.}
    \label{fig:6cycle}
\end{figure}


Continue to let $B_k = S_2 \wr S_k$ denote the wreath product of $S_2$ with $S_k$. Then $B_k$ naturally acts on $V(C_k)$ and $H(C_k)$. The action on $V(C_k)$ factors through $B_k \to S_k$, but the action on $H(C_k)$ does not. We write $S_2 = \{\pm1 \}$ multiplicatively, and use the usual cycle notation for elements of $S_k$. We identify the dihedral group $D_k = \Aut(C_k)$ with the subgroup of $B_k$ generated by the elements $\tau$ and $\rho$, where \[\tau = (-1, -1, \ldots, -1; (1, k)(2,k-1)(3, k-2)\cdots )\] and \[\rho = (1, 1, \ldots, 1; (1,2,3,\ldots,k)).\] The element $\tau$ is a flip of $C_k$ about the perpendicular bisector of the edge between $1$ and $k$ and $\rho$ is a rotation. All future references to the dihedral group $D_k$ are to this particular subgroup of $B_k$.

\begin{rem}
   The hyperoctahedral group $B_k$ is essential to the geometric picture. Our computation of $\mathsf{b}_{1,r}^{\mathrm{nrt}}$ sees $\Mbar_{1}^{\mathrm{nrt}}(r)$ as gluing caterpillar curves along marked points to a cycle and then decorating the cycle. The wreath product is necessary for the first step as it records the half-edge orientations on each vertex, which is not seen by the standard embedding\footnote{Further, the dihedral group fails to be a subgroup of the symmetric group when $k=2$, providing another motivation for the appearance of the wreath product.} $D_k< S_k$. This point of view has been taken by Petersen in \cite{semiclassicalremark} for the calculation of $\chi^{S_n}(\Mbar_{1,n}^{\mathrm{nrt}})$.
\end{rem}

\begin{defn}
Let $d \geq 2$. An $(r, d)$-\textit{\textbf{decoration}} $\delta$ of $C_k$ is a pair of functions
\[ s:  E(C_k) \to \Z_{>0} \text{ and } c: V(C_k) \to (\P^r)^{\C^\star} \]
such that $c$ is a graph colouring, and
\[ \sum_{e \in E(C_k)} s(e) = d. \]
We put $\tilde{\Gamma}_{r, d}(k)$ for the set of all $(r, d)$-decorations of $C_k$.
\end{defn}
\begin{rem}
    The definition of $(r,d)$-decorations makes sense on an arbitrary graph. In the case of $C_k$, an $(r, d)$-decoration defines an object of the groupoid $\mathsf{LG}_{1, 0, d}(r)^k$.
\end{rem}

The dihedral group $D_k$, thought of as a subgroup of $B_k$ as above, acts on $\tilde{\Gamma}_{r, d}(k)$. We put \[\Gamma_{r, d}(k): = \tilde{\Gamma}_{r, d}(k) / D_k. \]
Then $\Gamma_{r, d}(k)$ is a set of isomorphism class representatives for the groupoid $\mathsf{LG}_{1, 0, d}^\circ(r)^k$, and given $[\delta] \in \Gamma_{r, d}(k)$, a choice of lift to $\delta \in \tilde{\Gamma}_{r, d}(k)$ defines a subgroup $D_k(\delta)$ of $D_k$, by taking the stabiliser of $\delta$. We encode the combinatorics of these decorations as follows: for each integer $r \geq 1$, define a graded $\mathbb{B}$-space $\mathsf{Dih}_r$ by
\[ \mathsf{Dih}_r(k, d) := \coprod_{[\delta] \in \Gamma_{r, d}(k)} \Ind_{D_k(\delta)}^{B_k} \Spec(\C). \]
Different choices of lifts of decorations $\delta$ result in conjugate subgroups, so $\mathsf{Dih}_{r}(k, d)$ is well-defined. Observe that $\mathsf{Dih}_r$ is supported in degrees $k, d \geq 2$. We will now explain how marked chains of $\P^1$'s can be attached at the vertices of cycles in $\mathbb{B}$-space $\mathsf{Dih}_r$.

\subsection{Caterpillars} We refer to the marked genus-zero curves contracted by genus-one maps without rational tails as \textit{caterpillars}.
\begin{defn}
Let $n \geq 1$. The moduli space $\text{Cat}_{2,n}$ of $n$-pointed caterpillars is the locally closed subscheme of $\Mbar_{0, n+2}$ consisting of pointed rational stable curves $(C, \{q_i\}_{i = 1}^{n + 2})$ that satisfy the following two conditions: 
\begin{itemize}
    \item the dual graph of $C$ is a path,
    \item the first two marked points $q_1,q_2$ are on the two ends\footnote{The dual graph of $(C, \{q_i\}_{i = 1}^{n + 2})$ admits a orientation from $q_1$ to $q_2$.} of the path.
\end{itemize}
\end{defn}

Let $S_{2,n}\subset S_{n+2}$ be the Young subgroup of the partition $[n+2]=\{1,2\}\sqcup\{3, \dots, n+2\}$, so that $S_{2,n}\cong S_2\times S_n$. We observe that $\mathrm{Cat}_{2,n}$ is an $S_2\times S_n$-invariant subspace of $\Mbar_{0,n+2}$. 
\begin{defn}
Define the $S_2\times \mathbb{S}$-space $\mathsf{Cat}$ as $\mathsf{Cat}(n) := \text{Cat}_{2,n}$ for $n>0$ and set $\Cat(0)$ to be a point with the trivial action of $S_2 \times S_0$. \end{defn}
We recall Bergstr\"om--Minabe's calculation of the $S_2\times \mathbb{S}$-equivariant Euler characteristic of $\mathsf{Cat}$ in \S\ref{sec-cat}. Theorem \ref{thm-locnrt} now follows immediately from the following lemma.
\begin{lem}\label{lem-locnrt}
When $d = 0$, we have
\[\Mbar_{1, n}^{\mathrm{nrt}}(\P^r, 0)^{\C^\star} \cong \Mbar_{1, n}^{\mathrm{nrt}} \times \left(\P^r\right)^{\C^\star}. \]
For $d > 0$, we have an $S_n$-equivariant partition
\[ \Mbar_{1, n}^{\mathrm{nrt}}(\P^r, d)^{\C^\star} \cong \coprod_{k \geq 2} \coprod_{[\delta] \in \Gamma_{r, d}(k)} \left(\coprod_{\substack{{n_1 + \cdots + n_k = n} \\{n_i  \geq  0} }} \Ind_{S_{n_1} \times \cdots \times S_{n_k}}^{S_n} \Cat(n_1) \times \cdots \times \Cat(n_k)\right)/D_k(\delta), \]
where $D_k(\delta)$ acts as a subgroup of $D_k \leq B_k$.
\end{lem}
\begin{proof}
The statement for $d = 0$ is clear. For $d > 0$, we first adapt equation 
(\ref{eq:equivgrabpand}) to the case of maps without rational tails, writing
\[ \Mbar_{1, n}^{\mathrm{nrt}}(\P^r, d)^{\C^\star} \cong \coprod_{k \geq 2}\coprod_{(\mathcal{G}, c) \in \pi_0\left(\widetilde{\mathsf{LG}}^{\circ}_{1, n, d}(r)^{k}\right)} \Mbar^\circ_{\mathcal{G}, c}, \]
where we set 
\[\Mbar^\circ_{\mathcal{G}, c} \subset \Mbar_{\mathcal{G}, c} \] to be the subset of $\C^\star$-fixed maps without rational tails. Now adapting equation (\ref{eq:equivproduct}) to this setting, we write
\[\Mbar^\circ_{\mathcal{G}, c} \cong \left(\Ind_{\prod_{v \in V(\mathcal{G})}S_{n(v)}}^{S_n} \prod_{v \in V(\mathcal{G})} \mathsf{Cat}(n(v)) \right)/\Aut(\mathcal{G}, c), \]
where $\Aut(\mathcal{G}, c)$ acts on the half-edges. Now define a set map
\[ \pi_0\left(\widetilde{\mathsf{LG}}_{1, n, d}^\circ(r)^k\right) \to \Gamma_{r, d}(k) \]
by forgetting the marking function. If we let $[\delta] \in \Gamma_{r, d}(k)$ denote the image of $(\mathcal{G}, c)$ under this functor and $\delta \in \tilde{\Gamma}_{r, d}(k)$ its chosen lift, we get an identification of the underlying graph of $\mathcal{G}$ with $C_k$, which also gives the data of a function $f: V(C_k) \to \Z_{\geq 0}$ with $\sum_{i = 1}^{k} f(i) = n$, by recording how many markings each vertex of $\mathcal{G}$ supports. Let $\mathrm{Part}(C_k, n)$ denote the set of such functions. The group $D_k(\delta)$ acts on $\mathrm{Part}(C_k, n)$, and our choice of lift $\delta$ induces an identification of $\Aut(\mathcal{G}, c)$ with the stabiliser $\mathrm{Stab}_{D_k(\delta)}(f)$. If we let $n_i := f(i)$ under the identification of $V(C_k)$ with $\{1, \ldots, k\}$,  we can write
\[\Mbar^\circ_{\mathcal{G}, c} \cong \left(\Ind_{\mathrm{Stab}_{D_k(\delta)}(f)}^{D_k(\delta)} \Ind_{S_{f(1)} \times \cdots \times S_{f(k)}}^{S_n} \prod_{i = 1}^k \mathsf{Cat}(f(i)) \right)/D_k(\delta). \]
The fibre over $[\delta]$ is identified with 
\[  \frac{\mathrm{Part}(C_k, n)}{D_k(\delta)}.\] Hence
\begin{align*}
    \Mbar_{1, n}^{\mathrm{nrt}}(\P^r, d) &\cong \coprod_{[\delta] \in \Gamma_{r, d}(k)} \coprod_{[f] \in  \frac{\mathrm{Part}(C_k, n)}{D_k(\delta)} }\left(\Ind_{\mathrm{Stab}_{D_k(\delta)}(f)}^{D_k(\delta)} \Ind_{S_{f(1)} \times \cdots \times S_{f(k)}}^{S_n} \prod_{i = 1}^k \mathsf{Cat}(f(i)) \right)/D_k(\delta) \\& \cong \coprod_{[\delta] \in \Gamma_{r, d}(k)} \left(\coprod_{f \in \mathrm{Part}(C_k, n)} \Ind_{S_{f(1)}\times \cdots \times S_{f(k)}}^{S_n} \prod_{i = 1}^{k} \Cat(f(i))\right) / D_k(\delta),
\end{align*}
as we wanted to show.
\end{proof}

\section{Calculation of $\mathsf{Dih}_r$}\label{sec-calculation}
In this section we compute the class of the graded $\mathbb{B}$-space $\mathsf{Dih}_r$ for $r \geq 1$. Recall that it is defined as in \S\ref{sec-localisation} by:
\[ \mathsf{Dih}_r(k, d) =  \coprod_{[\delta] \in \Gamma_{r, d}(k)}\Ind_{D_k(\delta)}^{B_k} \Spec(\C). \]
To state our theorem, we first define some polynomials in $r$. Fix an integer $d \geq 2$. Given an integer $k$ with $2 \leq k \leq d$ and a divisor $j \mid k$, we set
\begin{equation}\label{eq:Fequation}
F_{j, k, d}(r) : = \sum_{\substack{{\ell}\\{\ell \mid j\text{ and }k\mid d\ell}}} \mu\left( \frac{j}{\ell} \right) \binom{\frac{d\ell}{k} - 1}{\ell - 1}(r^\ell + (-1)^\ell r).
\end{equation}
Similarly, we set
\begin{equation}\label{eq:Gequation}G_{j,k, d}(r) : = \sum_{\substack{{\ell} \\{2\ell \mid j \text{ and }k\mid d\ell}}} \mu\left(\frac{j}{2\ell} \right) \binom{\frac{d\ell}{k} - 1}{\ell - 1}(r^{\ell + 1} + r^{\ell}).
\end{equation}
Note that both polynomials vanish when the sum is empty; for instance, $G_{j, k, d}(r) = 0$ if $j$ is odd. Now, for $j$ and $k$ as above, we set
\begin{equation}\label{eq:Aeqn}
\theta_{j, k, d}(r) := \frac{\varphi(j)}{2k} \sum_{i \mid \frac{k}{j}} F_{i, k, d}(r)
\end{equation}
and
\begin{equation}\label{eq:Beqn}
\eta_{k, d}(r) := \sum_{j \mid k} \frac{G_{j,k, d}(r)}{4}.
\end{equation}
These are the same factors $\theta_{j, k, d}(r)$ and $\eta_{k, d}(r)$ appearing in Theorems \ref{thm:main_formula} and \ref{thm:main-tech}. In this section, we prove the following formula for the graded wreath product symmetric function \[{\mathsf{e}(\mathsf{Dih}_r) \in K_0(\mathsf{M}) \otimes\Lambda(S_2)[[q]]},\] in terms of the combinatorial factors $\theta_{j,k,d}$ and $\eta_{k,d}$.
\begin{thm}\label{thm:dihedral_char}
Fix an integer $r \geq 1$. We have
 \begin{align*}
\mathsf{e}(\mathsf{Dih}_r) &= 
\sum_{d \geq 2} q^d \left(\sum_{k = 2}^{d} \left(\eta_{k ,d}(r)\mathfrak{p}_2^{k/2 - 1} \mathfrak{q}_1^2 +  \sum_{j \mid k} 
 \theta_{j,k, d}(r)  \mathfrak{p}_j^{k/j} \right)\right).
 \end{align*}

\end{thm}
 The proof of Theorem \ref{thm:dihedral_char} amounts to enumerating localisation graphs with fixed isomorphism type of automorphism group. Given $[\delta] \in \Gamma_{r, d}(k)$, a choice of lift $\delta$ determines a subgroup $D_k(\delta)$ of $D_k$, well-defined up to conjugacy. Given a conjugacy class of subgroups $[H]$ of $D_k$, we let
\[\Gamma^{[H]}_{r, d}(k) = \{[\delta] \in \Gamma_{r, d}(k) \mid [D_k(\delta)] = [H] \} \]
and
\[\gamma^{[H]}_{r, d}(k) := \left| \Gamma^{[H]}_{r, d}(k)  \right|. \]
The relevance of these definitions is as follows. Let $O_k$ denote a set of orbit representatives for the action of $D_k$ on its subgroups by conjugation. Grouping $[\delta]\in \Gamma_{r,d}(k)$ according to the conjugacy class types of $D_k(\delta)$, we have
\begin{equation}\label{eq:dih_decomp}
\mathsf{e}(\mathsf{Dih}_r(k, d)) = \sum_{[\delta] \in \Gamma_{r, d}(k)} \mathsf{e}^{2 \wr k}\left(\Ind_{D_k(\delta)}^{B_k} \Spec (\C)\right) = \sum_{[H] \in O_k} \gamma_{r, d}^{[H]}(k) \mathsf{e}^{2 \wr k}\left(\Ind_{H}^{B_k} \Spec(\C) \right),
\end{equation}
where $H$ is any choice of representative of $[H]$. We will spend the rest of this section explaining how to understand the set $O_k$ and giving formulas for $\gamma^{[H]}_{r, d}(k)$ as $[H]$ varies over all elements of $O_k$. 
\subsection{Conjugacy classes of dihedral subgroups} The subgroup structure of $D_k$ is well-known, see e.g. \cite[§3]{Con19}. 
\begin{prop}\label{prop-conjsub} Every subgroup of $D_k$ is uniquely of the form $\langle \rho^j, \rho^i\tau\rangle$ for a divisor $j|k$ and $0\leq i< j$.

If $k$ is odd, then every subgroup of $D_k$ is conjugate to a subgroup of the form $\langle \rho^j\rangle$ or $\langle \rho^j, \tau\rangle$, for a divisor $j \mid k$. If $k$ is even, then every subgroup of $D_k$ is conjugate to one of:
\begin{itemize}
\item $\langle \rho^j \rangle$ for $j \mid k$;
\item $\langle \rho^j, \tau\rangle$ for $j \mid k$;
\item $\langle \rho^j, \rho \tau \rangle$ for $j \mid k$ and $j$ even.
\end{itemize}
Therefore, the above list is a set of representatives of $O_k$.
\end{prop}
For a decorated graph in $\tilde{\Gamma}_{r, d}(k)$, adjacent vertices are always distinguished by their colours. Therefore, the subgroups $D_k(\delta)$ never contain elements which stabilise an edge $e \in C_k$. This observation is summarised by the following lemma.
 
\begin{lem}
Suppose $j$ divides $k$, and suppose $H = \langle \rho^j, \rho^c\tau\rangle$ for some $1 \leq c \leq k$. Then
\[ \gamma^{[H]}_{r, d}(k) = 0 \]
if any of the following conditions hold:
\begin{enumerate}
\item $k$ is odd;
\item $k$ is even, and $j$ is odd;
\item $k$, $j$, and $c$ are all even.
\end{enumerate}
\end{lem}
If $k$ and $j$ are even while $c$ is odd, then $[\langle \rho^j, \rho^c \tau \rangle] = [\langle \rho^{j'}, \rho \tau \rangle]$ where $j'$ is even. Thus we focus on computing $\gamma^{[\langle \rho^j \rangle]}_{r, d}(k)$ for all $k$ and $j \mid k$, as well as the numbers $\gamma^{[\langle \rho^j, \rho \tau\rangle]}_{r, d}(k)$ when $k$ is even and $j \mid k$ is an even divisor.

\subsection{Frobenius reciprocity and normalisers} 
To count the unlabeled decorations \[\gamma^{[H]}_{r,d}(k) = |\Gamma^{[H]}_{r,d}(k)|,\] it is convenient to introduce and count certain auxiliary sets of labeled decorations as follows. For a fixed subgroup $H \leq D_k$, let
\[ \tilde{\Gamma}^H_{r, d}(k) \subset \tilde{\Gamma}_{r, d}(k) \]
be the set of labeled decorations with stabiliser equal to $H$, and let
\[\tilde{\Gamma}^{[H]}_{r, d}(k) \subset \tilde{\Gamma}_{r, d}(k)\] be the set of labeled decorations with stabiliser conjugate to $H$. The following identifications follow from unwinding the definitions. They relate the sets of labeled decorations defined above to  $\Gamma^{[H]}_{r,d}(k)$.
\begin{lem}
    $D_k$ acts on $\tilde{\Gamma}^{[H]}_{r, d}(k)$, and \[ \Gamma^{[H]}_{r, d}(k) = \tilde{\Gamma}^{[H]}_{r, d}(k) / D_k .\] Let $N_{D_k}(H)$ be the normaliser of $H$ in $D_k$. Then we have an isomorphism of $D_k$-sets \[ \tilde{\Gamma}^{[H]}_{r, d}(k) \cong \Ind_{N_{D_k}(H)}^{D_k} \tilde{\Gamma}^H_{r, d}(k). \] By Frobenius reciprocity, this is equivalent to \[ \Gamma^{[H]}_{r, d}(k) \cong \tilde{\Gamma}^H_{r, d}(k) / N_{D_k}(H). \]
\end{lem}

Using the orbit counting theorem, we have
\[ \gamma^{[H]}_{r, d}(k) = \frac{1}{|N_{D_k}(H)|} \sum_{g \in N_{D_{k}}(H)} \tilde{\Gamma}^H_{r, d}(k)^g  
 =  \frac{|H|}{|N_{D_k}(H)|}\tilde{\gamma}^H_{r, d}(k).\]
The indices $|N_{D_k}(H):H|$ for the list of subgroups appearing in Proposition \ref{prop-conjsub} can be computed by elementary methods, leading to the following result.
\begin{lem}\label{lem:normaliser_indices}
When $H = \langle \rho^j \rangle$ for a divisor $j$ of $k$, we have $N_{D_k}(H) = D_k$ and in particular,
\[ \gamma_{r, d}^{[\langle \rho^j \rangle]}(k) = \frac{\tilde{\gamma}^{\langle \rho^j \rangle}_{r, d}(k)}{2j}. \]
When $k$ is even and $j$ is an even divisor of $k$, then for $H = \langle \rho^j, \rho \tau \rangle$, we have $N_{D_k}(H) = \langle \rho^{j/2}, \rho \tau\rangle$. In particular we have
\[  \gamma_{r, d}^{[\langle \rho^j, \rho \tau \rangle]}(k) = \frac{\tilde{\gamma}^{\langle \rho^j , \rho \tau\rangle}_{r, d}(k)}{2}. \]
\end{lem}
The lemma allows us to focus our attention on computing the numbers $\tilde{\gamma}^{H}_{r, d}(k)$ of labeled decorations with stabiliser equal to a fixed subgroup $H \leq D_k$.
\subsection{M\"obius inversion}
The sets 
$\tilde{\Gamma}_{r,d}^H(k)$ are defined as decorations of $C_k$ with stabiliser in $D_k$ strictly equal to $H$. We instead enumerate decorations with stabiliser containing $H$, and then apply M\"obius inversion to the subgroup lattice $\mathcal{L}(D_k)$ of $D_k$. We set the notation
\[ \tilde{\Gamma}_{r, d}(k)^H := \coprod_{H \subset H'} \tilde{\Gamma}_{r, d}^{H'}(k) \]
for decorations fixed by $H$. Let $\mu_{D_k}$ denote the M\"obius function of $\mathcal{L}(D_k)$. Then M\"obius inversion gives
\[ \tilde{\gamma}^H_{r, d}(k) = \sum_{H \leq H'} \mu_{\mathcal{L}}(H, H') \left|\tilde{\Gamma}_{r, d}(k)^{H'}\right|. \]
The following lemma computes $\mu_{D_k}$ in terms of the standard M\"obius function $\mu$ on the natural numbers.
\begin{lem}\label{lem:mobius}
    Suppose that $a$ and $b$ are divisors of $k$, such that $b \mid a$. Then we have
\[ \mu_{D_k}(\langle \rho^a\rangle, \langle \rho^b, \rho^c \tau \rangle) = -\frac{a}{b} \mu\left( \frac{a}{b}\right), \]
\[ \mu_{D_k}(\langle \rho^a\rangle, \langle \rho^b \rangle) = \mu\left( \frac{a}{b}\right), \]
\[ \mu_{D_k}(\langle \rho^a, \rho^c \tau \rangle, \langle \rho^b , \rho^{c} \tau \rangle) = \mu\left( \frac{a}{b}\right). \]
\end{lem}
\begin{proof}
    The second two cases follow from the fact that the intervals \[[\langle \rho^a \rangle, \langle \rho^b \rangle] \quad \mbox{and}\quad [\langle \rho^a, \rho^c \tau \rangle, \langle \rho^b, \rho^c\tau \rangle]\] are each isomorphic as lattices to the lattice of divisors of $a/b$. For the first case, we argue as follows. If $G$ is an arbitrary finite group, let $\mu_{G}$ denote the M\"obius function of the poset $\mathcal{L}(G)$ of subgroups of $G$. If $N$ is a normal subgroup of $G$ which is contained in another subgroup $H$ of $G$. Then
    \[\mu_{G}(H, G) = \mu_{G/N}(H/N, G/N), \]
    since the intervals $[H, G]$ in $\mathcal{L}(G)$ and $[H/N, G/N]$ in $\mathcal{L}(G/N, H/N)$ are isomorphic as lattices. Since cyclic subgroups of $D_k$ are normal, we have
    \[ \mu_{D_k}(\langle \rho^a \rangle , \langle \rho^{b}, \rho^c \tau \rangle) = \mu_{D_{a/b}}(1, D_{a/b}), \]
    and we have
    \begin{equation}\label{eq:dihmobius}
        \mu_{D_{s}}(1, D_{s}) = -s\mu(s),
    \end{equation}
    as is likely well-known. For completeness, we provide a proof of (\ref{eq:dihmobius}). When $s = 1$, we can calculate by hand that  $\mu_{D_1}(1, 1) + \mu_{D_1}(1, D_1)=0$ implies that $\mu_{D_1}(1, D_1) = -1$. Perform induction on $s$, and suppose that the formula holds for all $s'|s$. We apply the identity $$\sum_{1\leq H\leq D_s} \mu_{D_s}(1, H) = 0.$$ Using the list of subgroups of $D_s$, this can be expanded as $$\sum_{i|s}\mu_{D_s}(1, \langle \rho^i \rangle) + \sum_{\substack{i|s, i\neq 1\\ 0\leq j<i}}\mu_{D_s}(1, \langle \rho^i, \rho^j \tau\rangle) + \mu_{D_s}(1, D_s)=0.$$ We have computed $\mu_{D_s}(1, \langle \rho^i\rangle) = \mu(s/i)$. For each divisor $i|s, i\neq 1$, we have $\langle \rho^i, \rho^j \tau\rangle\cong D_{s/i}$,
    so from the inductive hypothesis there is $\mu_{D_s}(1, \langle \rho^i, \rho^j \tau\rangle) = -(s/i)\mu(s/i)$. Summing up, we have $$\sum_{i|s}\mu(s/i) - s\sum_{{i|s, i\neq 1}}\mu(s/i)+ \mu_{D_s}(1, D_s)=0.$$ There is an identity $\sum_{i|s}\mu(s/i)=0$, hence $$\sum_{i|s, i\neq 1}\mu(s/i)=-\mu(s)$$ so rearranging we have $\mu_{D_s}(1, D_s) = -s\mu(s)$ as desired. 
\end{proof}

\subsection{Quotient graphs} 
We enumerate the sets $\tilde{\Gamma}_{r, d}(k)^{H'}$ via the following construction of \emph{quotient graphs}. This is an established technique for enumerating graphs with automorphisms \cite[§1]{Hanlon}.
\begin{defn}
    Fix a $k$-cycle $C_k$ with vertex and half-edge labeling as in \S\ref{sec-localisation}. For a subgroup $H\subset D_k$ that is in the list from Proposition \ref{prop-conjsub}, the \textit{quotient graph} of $C_k$ under $H$, denoted as $C_{k, H}$, is the graph with vertices (resp. half-edges) as the quotient set $V(C_k)/H$ (resp. $H(C_k)/H$), which determines the set of edges as $E(C_k)/H$.

    A quotient graph $C_{k, H}$ is intuitively a `fundamental domain' of the $k$-cycle under the action of $H$. It is convenient to pick the following vertex labeling.
    \begin{itemize}
        \item For $H = \langle \rho^j\rangle$ where $j \mid k$, we label the vertices as $\{1,\dots, j\}$, with each element representing their class modulo $j$; the quotient graph is a cycle on the $j$ vertices.
        \item For $H = \langle \rho^j, \rho^c\tau\rangle$ where $1\leq c\leq k$ is odd, $k$ is even and $j$ is an even divisor of $k$, we label the vertices as $\{1,\dots, j/2+1\}$; the quotient graph is a path with $j/2$ edges on the $j/2+1$ vertices.
    \end{itemize}
\end{defn}

\begin{figure}[h]
    \centering
    \tikzset{every picture/.style={line width=0.55pt}} 

\begin{tikzpicture}[x=0.75pt,y=0.75pt,yscale=-1,xscale=1]

\draw    (230.39,77.13) -- (292.74,77.03) ;
\draw   (169.99,77.98) -- (138.91,132.02) -- (76.56,132.12) -- (45.3,78.17) -- (76.39,24.13) -- (138.74,24.03) -- cycle ;
\draw  [draw opacity=0][fill={rgb, 255:red, 248; green, 8; blue, 18 }  ,fill opacity=1 ] (135.86,24.03) .. controls (135.86,22.44) and (137.15,21.16) .. (138.74,21.16) .. controls (140.33,21.16) and (141.61,22.44) .. (141.61,24.03) .. controls (141.61,25.62) and (140.33,26.91) .. (138.74,26.91) .. controls (137.15,26.91) and (135.86,25.62) .. (135.86,24.03) -- cycle ;
\draw  [draw opacity=0][fill={rgb, 255:red, 248; green, 8; blue, 18 }  ,fill opacity=1 ] (136.03,132.02) .. controls (136.03,130.43) and (137.32,129.14) .. (138.91,129.14) .. controls (140.49,129.14) and (141.78,130.43) .. (141.78,132.02) .. controls (141.78,133.61) and (140.49,134.9) .. (138.91,134.9) .. controls (137.32,134.9) and (136.03,133.61) .. (136.03,132.02) -- cycle ;
\draw  [draw opacity=0][fill={rgb, 255:red, 248; green, 8; blue, 18 }  ,fill opacity=1 ] (42.42,78.17) .. controls (42.42,76.59) and (43.71,75.3) .. (45.3,75.3) .. controls (46.89,75.3) and (48.18,76.59) .. (48.18,78.17) .. controls (48.18,79.76) and (46.89,81.05) .. (45.3,81.05) .. controls (43.71,81.05) and (42.42,79.76) .. (42.42,78.17) -- cycle ;
\draw  [draw opacity=0][fill={rgb, 255:red, 0; green, 0; blue, 0 }  ,fill opacity=1 ] (73.51,24.13) .. controls (73.51,22.54) and (74.8,21.25) .. (76.39,21.25) .. controls (77.98,21.25) and (79.27,22.54) .. (79.27,24.13) .. controls (79.27,25.72) and (77.98,27.01) .. (76.39,27.01) .. controls (74.8,27.01) and (73.51,25.72) .. (73.51,24.13) -- cycle ;
\draw  [draw opacity=0][fill={rgb, 255:red, 0; green, 0; blue, 0 }  ,fill opacity=1 ] (167.12,77.98) .. controls (167.12,76.39) and (168.4,75.1) .. (169.99,75.1) .. controls (171.58,75.1) and (172.87,76.39) .. (172.87,77.98) .. controls (172.87,79.57) and (171.58,80.86) .. (169.99,80.86) .. controls (168.4,80.86) and (167.12,79.57) .. (167.12,77.98) -- cycle ;
\draw  [draw opacity=0][fill={rgb, 255:red, 0; green, 0; blue, 0 }  ,fill opacity=1 ] (73.68,132.12) .. controls (73.68,130.53) and (74.97,129.24) .. (76.56,129.24) .. controls (78.15,129.24) and (79.44,130.53) .. (79.44,132.12) .. controls (79.44,133.71) and (78.15,135) .. (76.56,135) .. controls (74.97,135) and (73.68,133.71) .. (73.68,132.12) -- cycle ;
\draw  [draw opacity=0][fill={rgb, 255:red, 248; green, 8; blue, 18 }  ,fill opacity=1 ] (289.86,77.03) .. controls (289.86,75.44) and (291.15,74.16) .. (292.74,74.16) .. controls (294.33,74.16) and (295.61,75.44) .. (295.61,77.03) .. controls (295.61,78.62) and (294.33,79.91) .. (292.74,79.91) .. controls (291.15,79.91) and (289.86,78.62) .. (289.86,77.03) -- cycle ;
\draw  [draw opacity=0][fill={rgb, 255:red, 0; green, 0; blue, 0 }  ,fill opacity=1 ] (227.51,77.13) .. controls (227.51,75.54) and (228.8,74.25) .. (230.39,74.25) .. controls (231.98,74.25) and (233.27,75.54) .. (233.27,77.13) .. controls (233.27,78.72) and (231.98,80.01) .. (230.39,80.01) .. controls (228.8,80.01) and (227.51,78.72) .. (227.51,77.13) -- cycle ;

\draw (65,9.4) node [anchor=north west][inner sep=0.75pt]  [font=\footnotesize]  {$1$};
\draw (142,8.4) node [anchor=north west][inner sep=0.75pt]  [font=\footnotesize]  {$2$};
\draw (175,70.4) node [anchor=north west][inner sep=0.75pt]  [font=\footnotesize]  {$3$};
\draw (65.91,134.42) node [anchor=north west][inner sep=0.75pt]  [font=\footnotesize]  {$5$};
\draw (31.91,70.42) node [anchor=north west][inner sep=0.75pt]  [font=\footnotesize]  {$6$};
\draw (101,25.4) node [anchor=north west][inner sep=0.75pt]  [color={rgb, 255:red, 248; green, 8; blue, 18 }  ,opacity=1 ]  {$1$};
\draw (143,49.4) node [anchor=north west][inner sep=0.75pt]  [color={rgb, 255:red, 248; green, 8; blue, 18 }  ,opacity=1 ]  {$1$};
\draw (143,89.4) node [anchor=north west][inner sep=0.75pt]  [color={rgb, 255:red, 248; green, 8; blue, 18 }  ,opacity=1 ]  {$1$};
\draw (101,114.4) node [anchor=north west][inner sep=0.75pt]  [color={rgb, 255:red, 248; green, 8; blue, 18 }  ,opacity=1 ]  {$1$};
\draw (61,86.4) node [anchor=north west][inner sep=0.75pt]  [color={rgb, 255:red, 248; green, 8; blue, 18 }  ,opacity=1 ]  {$1$};
\draw (63,46.4) node [anchor=north west][inner sep=0.75pt]  [color={rgb, 255:red, 248; green, 8; blue, 18 }  ,opacity=1 ]  {$1$};
\draw (140.91,132.54) node [anchor=north west][inner sep=0.75pt]  [font=\footnotesize]  {$4$};
\draw (219,61.4) node [anchor=north west][inner sep=0.75pt]  [font=\footnotesize]  {$1$};
\draw (294,61.4) node [anchor=north west][inner sep=0.75pt]  [font=\footnotesize]  {$2$};
\draw (255,78.4) node [anchor=north west][inner sep=0.75pt]  [color={rgb, 255:red, 248; green, 8; blue, 18 }  ,opacity=1 ]  {$1$};

\end{tikzpicture}
    \caption{A labeled decorated graph in $\tilde{\Gamma}_{r, 6}(6)^{\langle \rho^2, \rho\tau\rangle}$ with degree decoration in red and vertex colouring is given by red and black dots. Its quotient graph under $\langle \rho^2, \rho\tau\rangle$-symmetry is on the right.}
    \label{fig:quotient6cycle}
\end{figure}
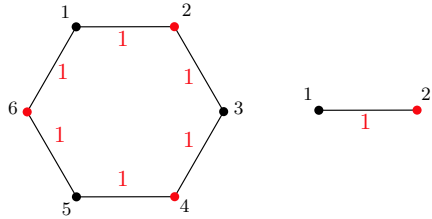

\begin{rem}
    The quotient half-edge ordering on $C_{k, H}$ gives a natural ordering of edges that is compatible with the one on $C_k$. When $C_{k, H}$ is a path, such an ordering is determined by the natural ordering on the vertices. When $C_{k, H}$ is a cycle, such an ordering is equivalent to labeling a half-edge, say $h_1^+$, on the cycle $C_i$. The ordering is not intrinsic to the quotient graph when $H = \langle \rho^2\rangle$, as illustrated below in Figure \ref{fig:2cycle}.
\end{rem}

\begin{figure}[h]
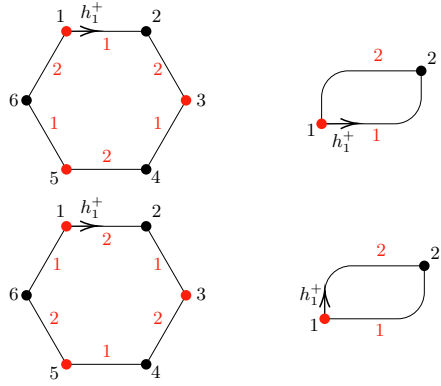

    \centering
    \includestandalone[width=0.35\textwidth]{2cycle}
    \caption{Two distinct labeled colourings with symmetry group $\langle \rho^2\rangle$. Their quotient graphs become the same after forgetting the edge orderings.}
    \label{fig:2cycle}
\end{figure}


We write $\tilde{\varepsilon}_{r, d}(k)$ for the set of $(r, d)$-decorations of a labelled path $P_k$ with $k + 1$ ordered vertices; formally, such a decoration is an $(r + 1)$-colouring $\delta_V: V(P_k) \to \{0, \ldots, r + 1\}$ of $P_k$, together with a function $\delta_E: E(P_k) \to \Z_{> 0}$ with $\sum_{e \in E(P_k)} \delta_E(e) = d$.



\begin{lem}\label{lem:quotient-bijections}
Suppose that $d \geq 2$ and $2 \leq k \leq d$. Assume that $j$ is a positive divisor of $k$.

\begin{enumerate}
    \item If $k \nmid dj$, then
\[ \tilde{\Gamma}_{r, d}(k)^{\langle \rho^j \rangle} = \varnothing. \]
If $k \mid dj$, then there is a bijection
\[ \tilde{\Gamma}_{r, d}(k)^{\langle \rho^j \rangle} \to \tilde{\Gamma}_{r, dj/k}(j). \]
\item Suppose $1 \leq c \leq k$. If $c$ is odd, $k$ and $j$ are both even, and $2k \mid dj$, then there is a bijection
\[ \tilde{\Gamma}_{r, d}^{\langle \rho^j, \rho^c \tau \rangle}(k) \to \tilde{\varepsilon}_{r, dj/2k}(j/2). \]
Otherwise,
\[ \tilde{\Gamma}_{r, d}(k)^{\langle \rho^j, \rho^c \tau\rangle} = \varnothing. \]

\end{enumerate}
\end{lem}
\begin{proof}
    In both cases, for any element in $\tilde{\Gamma}_{r,d}(k)^H$, the colouring and degree decoration are invariant under the induced $H$-action. Therefore, they descend along quotient maps \begin{align*}
        V(C_k)\to V(C_k)/H = V(C_{k, H}),\\ H(C_k)\to H(C_k)/H = H(C_{k, H}),\\E(C_k)\to E(C_k)/H = E(C_{k, H})
    \end{align*} and give colourings and degree decorations on the quotient graph $C_{k,H}$. This defines maps $$\tilde{\Gamma}_{r, d}(k)^{\langle \rho^j \rangle} \to \tilde{\Gamma}_{r, dj/k}(j) \quad \mbox{and}\quad \tilde{\Gamma}_{r, d}(k)^{\langle \rho^j, \rho^c \tau \rangle} \to \tilde{\varepsilon}_{r, dj/2k}(j/2).$$ Inverses of the maps are given by pulling back vertex colourings and degree decorations on $C_{k, H}$ along the quotient map of vertices and edges. One can check that the pullback does define a $(r,d)$-decoration on $C_k$. Therefore, the pairs of sets are in bijection.
\end{proof}

\subsection{Strata enumeration} Set
\[\omega_k(x) := x(x-1)^k \]
for the chromatic polynomial of a path with $k$ edges, and let
\[ \varpi_k(x) := (x - 1)^k + (-1)^k(x - 1) \]
denote the chromatic polynomial of $C_k$ \cite{birkhoff, whitney}. 
They yield straightforward formulas for the numbers $\tilde{\gamma}_{r, d}(k) = |\tilde{\Gamma}_{r, d}(k)|$ and $|\tilde{\varepsilon}_{r, d}(k)|$.
\begin{prop}
We have formulas
\[ |\tilde{\Gamma}_{r, d}(k)| =  \binom{d - 1}{k - 1} \varpi_{k}(r + 1)\quad\mbox{and}\quad |\tilde{\varepsilon}_{r, d}(k)| = \binom{d - 1}{k - 1} \omega_k(r + 1). \]
\end{prop}
\begin{proof}
In both cases, a decoration is chosen in two steps: choose a colouring, and then choose an ordered partition of $d$ into $k$ positive parts.
\end{proof}
The following corollary is now an immediate consequence of Lemma \ref{lem:quotient-bijections}.
\begin{cor}\label{cor:mobiusenumeration}
Suppose $d \geq 2$ and $2 \leq k \leq d$. Suppose $j \mid k$ and $k \mid dj$. Then
\[ \left|\tilde{\Gamma}_{r, d}(k)^{\langle \rho^j \rangle} \right| = \binom{\frac{dj}{k} - 1}{j - 1} \varpi_j(r + 1) \]
If $2\ell \mid k$ and $k \mid d \ell$, and in addition $1 \leq c \leq k$ is odd, then
\[\left| \tilde{\Gamma}_{r, d}(k)^{\langle \rho^{2\ell}, \rho^c \tau \rangle} \right| = \binom{\frac{d\ell}{k} - 1}{\ell - 1} \omega_{\ell}(r + 1). \]
\end{cor}
Now we apply M\"obius inversion to compute the number of decorations of $C_k$ with fixed stabiliser group. For the statement of the formula, recall that the polynomials $F_{j,k,d}(r)$ and $G_{j,k,d}(r)$ are defined in (\ref{eq:Fequation}) and (\ref{eq:Gequation}), respectively.
\begin{prop}\label{prop:labelled-decorations}
Suppose that $k \geq 2$ and $j$ divides $k$. Then
\[\tilde{\gamma}_{r, d}^{\langle \rho^j \rangle}(k) = F_{j, k, d}(r) - \frac{j}{2}G_{j, k, d}(r). \]
If $k$ is even and $j$ is an even divisor of $k$, then
\[\tilde{\gamma}_{r, d}^{\langle \rho^j ,\rho\tau\rangle}(k) = G_{j, k, d}(r).\]
\end{prop}

\begin{proof}
By M\"obius inversion, Lemma \ref{lem:mobius}, and Corollary \ref{cor:mobiusenumeration}, we have
\begin{align*}
\tilde{\gamma}^{\langle \rho^j \rangle}_{r, d}(k) &= \sum_{\ell \mid j} \mu_{\mathcal{L}}(\langle \rho^j\rangle, \langle\rho^{\ell}\rangle)|\tilde{\Gamma}_{r, d}(k)^{\langle \rho^\ell \rangle}| + \sum_{2\ell \mid j} \sum_{s = 1}^{\ell} \mu_{\mathcal{L}}(\langle\rho^j \rangle, \langle \rho^{2\ell}, \rho^{2s - 1} \tau \rangle)|\tilde{\Gamma}_{r, d}(k)^{\langle \rho^{2\ell}, \rho^{2s - 1} \tau \rangle}|
\\&= \sum_{\ell \mid j} \mu\left( \frac{j}{\ell}\right) 
\binom{\frac{d\ell}{k} - 1}{\ell - 1} \varpi_{\ell}(r + 1) - \sum_{2 \ell \mid j} \frac{j}{2\ell} \sum_{s = 1}^{\ell} \mu \left( \frac{j}{2\ell} \right) \binom{\frac{d\ell}{k} - 1}{\ell - 1} \omega_{\ell}(r + 1)
\\&= \sum_{\ell \mid j} \mu\left( \frac{j}{\ell}\right) 
\binom{\frac{d\ell}{k} - 1}{\ell - 1} \varpi_{\ell}(r + 1) - \frac{j}{2}\sum_{2 \ell \mid j}  \mu \left( \frac{j}{2\ell} \right) \binom{\frac{d\ell}{k} - 1}{\ell - 1} \omega_{\ell}(r + 1)
\\&= F_{j, k, d}(r) - \frac{j}{2} G_{j, k, d}(r),
\end{align*}
as desired. The second count is similar, except only the terms with $2\ell \mid j$ survive:
\begin{align*}
\tilde{\gamma}^{\langle \rho^j, \rho\tau \rangle}_{r, d}(k) &=  \sum_{2\ell \mid j} \mu_{P_k}(\langle\rho^j, \rho\tau \rangle, \langle \rho^{2\ell}, \rho \tau \rangle)|\tilde{\Gamma}_{r, d}(k)^{\langle \rho^{2\ell}, \rho \tau \rangle}|
\\&= \sum_{2 \ell \mid j} \mu \left( \frac{j}{2\ell} \right) \binom{\frac{d\ell}{k} - 1}{\ell - 1} \omega_{\ell}(r + 1)
\\&= G_{j, k, d}(r).
\end{align*}\end{proof}
The following corollary is immediate from Proposition \ref{prop:labelled-decorations} and Lemma \ref{lem:normaliser_indices}.
\begin{cor}\label{cor:enumeration}
Suppose that $k \geq 2$ and $j$ divides $k$. Then
\[\gamma_{r, d}^{[\langle \rho^j \rangle]}(k)=  \frac{F_{j, k, d}(r)}{2j} - \frac{G_{j, k, d}(r)}{4} \]
If $k$ is even and $j$ is an even divisor of $k$, then
\[ \gamma_{r, d}^{[\langle \rho^j, \rho \tau\rangle]}(k) = \frac{G_{j, k,d}(r)}{2}. \]
\end{cor}

\subsection{Computing the characters} The final step in our calculation of $\mathsf{e}(\Dih_r)$ is the computation of the characters $\mathsf{e}(\Ind_{H}^{B_k} \Spec \C)$.
\begin{prop}\label{prop:subgroup_chars}
Suppose $j \mid k$. 
\begin{enumerate}
\item If $H = \langle \rho^j\rangle$, then
\[ \mathsf{e}^{2 \wr k}(\Ind_H^{B_k} \Spec\C) = \frac{j}{k} \sum_{s \mid \frac{k}{j}} \varphi(s) \mathfrak{p}_{s}^{k/s} .\]
\item If $H = \langle \rho^j, \rho \tau\rangle$, then
\[ \mathsf{e}^{2 \wr k}(\Ind_H^{B_k} \Spec \C) = \frac{1}{2} \mathfrak{p}_2^{k/2 - 1} \mathfrak{q}_1^2 +  \frac{j}{2k} \sum_{s \mid \frac{k}{j}} \varphi(s) \mathfrak{p}_{s}^{k/s}. \]
\end{enumerate}
\end{prop}

\begin{proof}
In both cases we have $\mathsf{e}^{2 \wr k}(\Ind_H^{B_k} \Spec\C) = \ch_{2 \wr k}(\Ind_{H}^{B_k} \mathrm{Triv}_H). $ As explained in \cite{semiclassicalremark}, 
\begin{align*}
    \ch_{2 \wr k}(\Ind_{H}^{B_k} \mathrm{Triv}_H) = \frac{1}{|H|} \sum_{h \in H} \Psi(h),
\end{align*}
where the cycle map $\Psi$ is defined in \S\ref{section-sym}. When $H = \langle \rho^j \rangle$, we have $|H| = \frac{k}{j}$, and for each divisor $s$ of $k/j$, the group $H$ contains exactly $\varphi(s)$ rotations of order $s$. The $S_k$ component of such a rotation admits a cycle decomposition as $k/s$ many $s$-cycles, while the $S_2^k$ component acts trivially. Therefore, $\Psi(h) = \mathfrak{p}_{s}^{k/s}$, and the first formula follows. In the second case, the order of the subgroup is $\frac{2k}{j}$. Each rotation acts in the same way as earlier, but each reflection $f$ fixes two vertices, and switches the orientations of the two half-edges at each of these vertices. Hence we have $\Psi(f) = \mathfrak{p}_2^{k/2 - 1}\mathfrak{q}_1^2$. Since there are $k/j$ distinct reflections, the second formula follows.
\end{proof}

\begin{proof}[Proof of Theorem \ref{thm:dihedral_char}]
Plugging Proposition \ref{prop:subgroup_chars} and Corollary \ref{cor:enumeration} into (\ref{eq:dih_decomp}), and then switching the order of summation, we obtain
\begin{align*}
    \mathsf{e}(\mathsf{Dih}_r) &= \sum_{d \geq 2} q^d \sum_{k = 2}^{d} \sum_{j \mid k} \left(\frac{G_{j, k, d}(r)}{4}\mathfrak{p}_2^{k/2 - 1} \mathfrak{q}_1^2 +  \left( \frac{F_{j, k, d}(r)}{2k} - \frac{jG_{j, k, d}(r)}{4k} + \frac{jG_{j, k, d}(r)}{4k}\right) \sum_{s \mid \frac{k}{j}}\varphi(s)\mathfrak{p}_{s}^{k/s}\right) \\&=
    \sum_{d \geq 2} q^d \sum_{k = 2}^{d} \sum_{j \mid k} \left(\frac{G_{j, k, d}(r)}{4}\mathfrak{p}_2^{k/2 - 1} \mathfrak{q}_1^2 +   \frac{F_{j, k, d}(r)}{2k} \sum_{s \mid \frac{k}{j}}\varphi(s)\mathfrak{p}_{s}^{k/s}\right) \\&= \sum_{d \geq 2} q^d \left(\sum_{k = 2}^{d} \eta_{k, d}(r) \pp_2^{k/2 - 1} \qq_1^2 + \sum_{s \mid k} \left(\mathfrak{p}_{s}^{k/s} \sum_{j \mid \frac{k}{s}} \frac{\varphi(s)}{2k}  F_{j, k, d}(r) \right)\right)
    \\&= \sum_{d \geq 0} q^d \sum_{k = 2}^{d} \left(\eta_{k, d}(r) \mathfrak{p}_2^{k/2 - 1} \mathfrak{q}_1^2 + \left(\sum_{s \mid k} \theta_{s, k, d}(r) \mathfrak{p}_{s}^{k/s}\right)\right)
\end{align*}
which is equivalent to the statement of the Theorem.
\end{proof}

\section{Calculation of $\mathscr{B}_r$} \label{sec-cat}
In this section, we will assemble the results of previous sections to derive our formula for $\mathscr{B}_r$. From Theorem \ref{thm-locnrt} we have that
\begin{equation} \mathscr{B}_r = \mathscr{B}_0 + \mathsf{e}(\mathsf{Dih}_r) \circ_{S_2} \mathsf{e}(\Cat).
\end{equation}
The Serre characteristic $\mathsf{e}(\mathsf{Dih}_r)$ was computed in \S\ref{sec-calculation}. In this section, we will explain a formula for $\mathsf{e}(\Cat)$ that appears in \cite[Remark 4.2]{BergstromMinabe1}, and then compute the plethysm $\mathsf{e}(\mathsf{Dih}_r) \circ_{S_2} \mathsf{e}(\Cat)$. Combining this with Getzler's calculation of $\mathscr{B}_{0}$ will give our main technical Theorem \ref{thm:main-tech}, recalled below, which implies Theorem \ref{thm:main_formula}.

\begin{customthm}{B}\label{thm:main-tech}
We have
\begin{align*}
\mathscr{B}_r &= (r + 1)\left(\mathscr{A}_1 + \frac{\dot{\mathscr{A}}_0(\dot{\mathscr{A}}_0 + 1) + \frac{1}{4}\psi_2(\mathscr{A}_0'')}{1 - \psi_2(\mathscr{A}_0'')} - \frac{1}{2} \sum_{n \geq 1} \frac{\varphi(n)}{n}\log(1 - \psi_n(\mathscr{A}_0''))\right) \\&\phantom{space}+ \sum_{d \geq 2} q^d \sum_{k = 2}^{d} \left(\eta_{k, d}(r) \frac{(1 + 2 \dot{\mathscr{A}}_0)^2}{(1 - \psi_2(\mathscr{A}_0''))^{k/2 + 1}} + \sum_{j \mid k} \theta_{j,k,d}(r)\frac{1}{(1 - \psi_j(\mathscr{A}_0''))^{k/j}}\right).
\end{align*}
\end{customthm}

\subsection{The Serre characteristic of $\Cat$.} Recall that $\Cat$ is the $S_2 \times \mathbb{S}$ space such that $\Cat(n) = \mathrm{Cat}_{2,n}$ for $n \geq 1$, and $\Cat(0)$ is a point. In this section we explain a formula for
$$\mathsf{e}(\mathsf{Cat}):=\sum_{n\geq 0}\mathsf{e}^{S_2 \times S_n}(\mathsf{Cat}(n)) \in 
K_0(\mathsf{M}) \otimes \Lambda_{2} \otimes \Lambda$$
given by Bergstr\"om--Minabe \cite{BergstromMinabe1}.
The formula, recalled in Theorem \ref{thm-chcat} below, expresses $\sfe(\mathsf{Cat})$ in terms of the known series
\[\mathscr{A}_0 := \sum_{n \geq 3} \mathsf{e}^{S_n}(\M_{0, n}), \]
and its derivatives $$\mathscr{A}_0'':=\frac{\partial^2 \mathscr{A}_0}{\partial p_1^2} \quad \mbox{and} \quad \dot{\mathscr{A}}_0:=\frac{\partial \mathscr{A}_0}{\partial p_2}.$$

\begin{thm}\label{thm-chcat}
We have
\[ \mathsf{e}(\Cat) = \frac{1}{2}p_1^2\otimes  \frac{1}{1-\mathscr{A}_0''}+\frac{1}{2}p_2\otimes \frac{1 + 2\dot{\mathscr{A}}_0}{1 - \psi_2(\mathscr{A}_0'')}. \]
\end{thm}

The formula is stated in \cite[Remark 4.2]{BergstromMinabe2}. Here we explain how it follows from their main result on Serre characteristics of the Losev--Manin spaces $\overline{\M}_{0,2|n}$. In Proposition 4.1 of loc. lit., they prove that \begin{equation}\label{eq:BM}
h_2\otimes 1 + \sum_{n=1}^\infty \mathsf{e}^{S_2 \times S_n}(\overline{\M}_{0,2|n}) = \frac{1}{2}p_1^2\otimes \frac{1}{1 - \sum_{i = 1}^{\infty} f_i} +\frac{1}{2}p_2\otimes \frac{1 + \sum_{i = 1}^{\infty} g_i}{1 - \psi_2\left( \sum_{i = 1}^{\infty}f_i \right)} \end{equation} where $f_n, g_n \in K_0(\mathsf{M}) \otimes \Lambda$ encode the Serre characteristic of the dense torus $\M_{0,2|n}\subset\overline{\M}_{0,2|n}$: $$\sum_{n=1}^\infty \mathsf{e}^{S_2\times S_n}({\M}_{0,2|n}) = \frac{1}{2}p_1^2\otimes \sum_{n=1}^\infty f_n+ \frac{1}{2}p_2\otimes \sum_{n=1}^\infty g_n.$$
The modular interpretation of $\mathsf{Cat}$ implies that its characteristic is obtained from $\sum_{n} \mathsf{e}^{S_2 \times S_n}(\Mbar_{0,2|n})$ via `distinguishing' the light marked points, i.e., by applying the plethystic logarithm as in \cite[Proposition 8.6]{GetzlerKapranov}:
\[ \mathsf{e}(\mathsf{Cat}) = \left(h_2 \otimes 1 + \sum_{n \geq 1} \mathsf{e}^{S_2 \times S_n}(\Mbar_{0, 2|n})\right) \circ(1 \otimes \Log(p_1)). \]
To compute $\mathsf{e}(\mathsf{Cat})$, it then suffices to compute $ \sum_{i=1}^\infty f_i\circ \Log(p_1)$ and $\sum_{i=1}^\infty g_i \circ \Log(p_1)$. We calculate \begin{align*}
    \sum_{n \geq 1} \mathsf{e}^{S_2 \times S_n}(\M_{0,2|n})\circ (1\otimes \Log(p_1)) &= \sum_{n \geq 1} \mathsf{e}^{S_2 \times S_n}(\Res_{S_2\times S_n}^{S_{n+2}}\M_{0,n+2}) \\&= \frac{1}{2}p_1^2 \otimes \mathscr{A}_0'' + p_2 \otimes \dot{\mathscr{A}}_0,
\end{align*}    
where the first equality follows from \cite[Proposition 4.3]{hodgehl}, and the second follows from Lemma \ref{lem-sn+2}. From the way $f_i$ and $g_i$ are defined, it follows that
\[ \sum_{i = 1}^{\infty} f_i \circ \Log(p_1) = \mathscr{A}_0'' \quad \mbox{and}\quad \sum_{i = 1}^{\infty} g_i \circ \Log(p_1) = 2 \dot{\mathscr{A}}_0. \]

Plugging this back into (\ref{eq:BM}) gives the formula for $\mathsf{e}(\mathsf{Cat})$ as claimed in Theorem \ref{thm-chcat}. With this formula in hand, we can now compute
\[\mathsf{e}(\mathsf{Dih}_r) \circ_{S_2} \mathsf{e}(\Cat), \]
using the formula for $\mathsf{e}(\mathsf{Dih}_r)$ in Theorem \ref{thm:dihedral_char} and basic properties of wreath product plethysm.

\begin{thm}\label{thm:e-nrt}
For an integer $r \geq 1$, we have 
\begin{align*}
\mathsf{e}(\mathsf{Dih}_r) \circ_{S_2} \mathsf{e}(\Cat) = \sum_{d \geq 2} q^d \sum_{k = 2}^{d} \left(\eta_{k, d}(r) \frac{(1 + 2 \dot{\mathscr{A}}_0)^2}{(1 - \psi_2(\mathscr{A}_0''))^{k/2 + 1}} + \sum_{j \mid k} \theta_{j,k,d}(r)\frac{1}{(1 - \psi_j(\mathscr{A}_0''))^{k/j}}\right).
\end{align*}

\end{thm}



\begin{proof}
This reduces to computing
\[ \mathfrak{p}_{j}^{k/j} \circ_{S_2} \mathsf{e}(\Cat) \quad \mbox{and} \quad  \mathfrak{p}_2^{k/2 - 1}\mathfrak{q}_1^2\circ_{S_2} \mathsf{e}(\Cat) \]
for integers $j, k \geq 1$ with $j \mid k$. From Lemmas \ref{lem:frakpleth1} and \ref{lem:frakpleth2}, we have
\[ \mathfrak{p}_{j}^{k/j} \circ_{S_2} \mathsf{e}(\Cat) = \frac{1}{(1 - \psi_{j}(\mathscr{A}_0''))^{k/j}},  \]
\begin{align*} \mathfrak{p}_2^{k/2 - 1}\mathfrak{q}_1^2\circ_{S_2} \mathsf{e}(\Cat) &= \left( \frac{1}{1 - \psi_2(\mathscr{A}_0'')} \right)^{k/2 - 1}\left( \frac{1 + 2\dot{\mathscr{A}}_0}{1 - \psi_2(\mathscr{A}_0'')} \right)^2 \\ &= \frac{(1 + 2 \dot{\mathscr{A}}_0)^2}{(1 - \psi_2(\mathscr{A}_0''))^{k/2 + 1}}
\end{align*}
The result follows by plugging these expressions into Theorem \ref{thm:dihedral_char}.
\end{proof}
We now conclude by indicating how Theorems \ref{thm:main_formula} and \ref{thm:main-tech} follows from Theorem \ref{thm:e-nrt}: in \cite{GetzlerSemiClassical} (with a typo corrected in \cite[Page 1]{semiclassicalremark}), Getzler proves 
\begin{equation}\label{eq:nrtdeg0}
\sum_{n \geq 1} \mathsf{e}^{S_n}(\Mbar_{1, n}^{\mathrm{nrt}}) = \mathscr{A}_1 + \frac{\dot{\mathscr{A}}_0(1 + \dot{\mathscr{A}}_0 ) + \frac{1}{4} \psi_2(\mathscr{A}_0'') }{1 - \psi_2(\mathscr{A}_0'')} -\frac{1}{2} \sum_{n \geq 1} \frac{\varphi(n)}{n}\log(1 - \psi_n(\mathscr{A}_0'')).
\end{equation}
Theorem \ref{thm:main-tech}, and hence Theorem \ref{thm:main_formula}, follows from (\ref{eq:nrtdeg0}) with Lemma \ref{lem-locnrt}, which implies that 
\[ \mathscr{B}_0 = (r +  1)\sum_{n \geq 1} \mathsf{e}^{S_n}(\Mbar_{1, n}^{\mathrm{nrt}}). \]

\section{Calculations}\label{sec:data}
Calculations of $\chi^{S_{n}}(\Mbar_{1, n}(\P^r, d))$ for small $d$ and $n$, performed in SageMath\footnote{The Sage notebook used to produce these calculations can be downloaded \href{https://drive.google.com/file/d/1Q-TanqsALZII2Vrbl_wLsF8lhZ56Sx7r/view}{here}.} using Theorem \ref{thm:main_formula}, are found in Tables \ref{table:degree1-frob-chars}, \ref{table:degree2-frob-chars}, and $\ref{table:degree3-frob-chars}$. The first step in producing these tables is to compute $\mathsf{b}_{1, r}^{\mathrm{nrt}}$ as follows. First, Theorem \ref{thm-locnrt} gives the formula \[\mathscr{B}_r = \mathscr{B}_0 + \mathsf{e}(\mathsf{Dih}_r) \circ_{S_2} \mathsf{e}(\Cat).\] Theorem \ref{thm:e-nrt} gives a closed formula for the plethysm $\mathsf{e}(\mathsf{Dih}_r)\circ_{S_2}\mathsf{e}(\Cat)$ in terms of the derivatives $$\mathscr{A}_0''=\frac{\partial^2 \mathscr{A}_0}{\partial p_1^2}\quad \mbox{and} \quad \dot{\mathscr{A}}_0=\frac{\partial \mathscr{A}_0}{\partial p_2}.$$ Their specialisations to symmetric functions are denoted as $\mathsf{a}_0''$ and $\dot{\mathsf{a}}_0$, respectively. Therefore, to determine $\mathsf{e}(\mathsf{Dih}_r)\circ_{S_2}\mathsf{e}(\Cat)$, it suffices to plug in formulas for $\mathsf{a}_0''$ and $\dot{\mathsf{a}}_0$ into Theorem \ref{thm:e-nrt}. To do this, we differentiate the following formula of
Getzler \cite{GetzlerGenusZero}: $$\mathsf{a}_0 = \frac{1}{2}(1+p_1)^2\sum_{n=1}^\infty \frac{\varphi(n)}{n}\log(1+p_n)-\frac{1}{4}(2p_1 + 3p_1^2 + p_2).$$ This leads to $$\mathsf{a}_0'' = \sum_{n=1}^\infty \frac{\varphi(n)}{n}\log(1+p_n)\quad \mbox{and}\quad\dot{\mathsf{a}}_0 = \frac{(1+p_1)^2}{4(1+p_2)}-\frac{1}{4}.$$

The term $\mathscr{B}_0$ specialises to $(r+1)\sum_{n\geq 1}\chi^{S_n}(\Mbar_{1,n}^{\mathrm{nrt}})$, which is determined in \cite{GetzlerSemiClassical} and \cite{semiclassicalremark} in terms of $\mathsf{a}_0''$, $\dot{\mathsf{a}}_0$, and $\mathsf{a}_1= \sum_{n=1}^\infty \chi^{S_n}(\mathcal{M}_{1,n})$. A formula for $\mathsf{a}_1$ is derived by Getzler in \cite[(5.5)]{GetzlerMHM}. Plugging these formulas into Theorem \ref{thm:e-nrt} and taking the plethysm with $p_1 + \mathsf{b}_{0, r}'/(r + 1)$, we obtain Tables \ref{table:degree1-frob-chars}, \ref{table:degree2-frob-chars}, and \ref{table:degree3-frob-chars}. Several additional tables, as well as a discussion of the contribution from excess components to $\chi^{S_n}(\Mbar_{1, n}(\P^r, d))$ can be found in the first arXiv version of the present article.

\begin{table}[H]
\def\arraystretch{1.5}
\begin{tabular}{|l|l|}
\hline
$n$ & $\chi^{S_n} \left(\Mbar_{1, n}(\P^r, 1)\right)$                                                           \\ \hline
$0$ & $4\binom{r + 1}{2}$                                                                                                                       \\ \hline
$1$ & $12 \binom{r + 1}{2}s_1$                                      \\ \hline
$2$ & $32\binom{r + 1}{2}s_2 + 12 \binom{r + 1}{2}s_{1,1}$ \\ \hline
$3$ & $86 \binom{r + 1}{2}s_3 + 56 \binom{r + 1}{2}s_{2,1} + 4 \binom{r + 1}{2}s_{1,1,1}$   \\ \hline
\end{tabular}
\caption{The $S_n$-equivariant topological Euler characteristic of $\Mbar_{1, n}(\P^r, 1)$, for $n \leq 3$.}
\label{table:degree1-frob-chars}
\end{table}

\begin{table}[H]
\def\arraystretch{1.5}
\begin{tabular}{|l|l|}
\hline
$n$ & $\chi^{S_n} \left(\Mbar_{1, n}(\P^r, 2)\right)$                                                           \\ \hline
$0$ & 

$24 \binom{r + 1}{3} + 17 \binom{r + 1}{2}$                             \\ \hline
$1$ & 
$\left[108\binom{r + 1}{3} + 58 \binom{r + 1}{2}\right]s_1$                                    \\ \hline
$2$ & $\left[339 \binom{r + 1}{3} + 171 \binom{r + 1}{2} + 6 \binom{r + 1}{1}\right]s_2 + \left[168 \binom{r + 1}{3} + 65 \binom{r + 1}{2} + 6\binom{r + 1}{1}\right]s_{1,1}$ \\ \hline

$3$ & \begin{tabular}{@{}c@{}}$\left[1176\binom{r + 1}{3} + 498 \binom{r + 1}{2} + 6\binom{r + 1}{1} \right]s_3 + \left[996\binom{r + 1}{3} + 396 \binom{r + 1}{2} + 12 \binom{r + 1}{1}\right]s_{2,1}$ \\ 
 $+ \left[144 \binom{r + 1}{3} + 40\binom{r + 1}{2} + 6 \binom{r + 1}{1}\right]s_{1,1,1}$\end{tabular}    \\ \hline

\end{tabular}

\caption{The $S_n$-equivariant topological Euler characteristic of $\Mbar_{1, n}(\P^r, 2)$, for $n \leq 3$.}
\label{table:degree2-frob-chars}
\end{table}

\begin{table}[H]
\def\arraystretch{1.5}
\begin{tabular}{|l|l|}
\hline
$n$ & $\chi^{S_n} \left(\Mbar_{1, n}(\P^r, 3)\right)$                                                           \\ \hline
$0$ & 

$216\binom{r + 1}{4} +247\binom{r + 1}{3} +55\binom{r + 1}{2}$                             \\ \hline
$1$ & 
$\left[1300\binom{r + 1}{4} +1365\binom{r + 1}{3} +260\binom{r + 1}{2}\right]s_1$                                    \\ \hline
$2$ & $\left[5380\binom{r + 1}{4} +5319\binom{r + 1}{3} +945\binom{r + 1}{2} \right]s_2 + \left[3156\binom{r + 1}{4} +2991\binom{r + 1}{3} +503\binom{r + 1}{2}\right]s_{1,1}$ \\ \hline
\end{tabular}
\caption{The $S_n$-equivariant topological Euler characteristic of $\Mbar_{1, n}(\P^r, 3)$, for $n \leq 2$.}
\label{table:degree3-frob-chars}
\end{table}

\bibliographystyle{amsalpha}
\bibliography{reference}

@ARTICLE{KSGraphEnumeration,
       author = {{Kannan}, Siddarth and {Song}, Terry Dekun},
        title = "{Graph enumeration for moduli spaces of curves and maps}",
      journal = {arXiv e-prints},
         year = 2025,
        month = sep,
          eid = {arXiv:2509.18298},
        pages = {arXiv:2509.18298},
          doi = {10.48550/arXiv.2509.18298},
archivePrefix = {arXiv},
       eprint = {2509.18298},
 primaryClass = {math.AG},
       adsurl = {https://ui.adsabs.harvard.edu/abs/2025arXiv250918298K}
}

@ARTICLE{tosteson,
       author = {{Tosteson}, Philip},
        title = "{Representation stability for the Kontsevich space of stable maps}",
      journal = {arXiv e-prints},
         year = 2022,
        month = jan,
          eid = {arXiv:2201.05100},
        pages = {arXiv:2201.05100},
          doi = {10.48550/arXiv.2201.05100},
archivePrefix = {arXiv},
       eprint = {2201.05100},
 primaryClass = {math.AG},
       adsurl = {https://ui.adsabs.harvard.edu/abs/2022arXiv220105100T}
}

@Article{GetzlerMHM,
 Author = {Getzler, E.},
 Title = {Resolving mixed {Hodge} modules on configuration spaces},
 FJournal = {Duke Mathematical Journal},
 Journal = {Duke Math. J.},
 ISSN = {0012-7094},
 Volume = {96},
 Number = {1},
 Pages = {175--203},
 Year = {1999},
 Language = {English},
 DOI = {10.1215/S0012-7094-99-09605-9},
 zbMATH = {1425170},
 Zbl = {0986.14005}
}

@Article{Gorsky,
 Author = {Gorsky, Eugene},
 Title = {The equivariant {Euler} characteristic of moduli spaces of curves},
 FJournal = {Advances in Mathematics},
 Journal = {Adv. Math.},
 ISSN = {0001-8708},
 Volume = {250},
 Pages = {588--595},
 Year = {2014},
 Language = {English},
 DOI = {10.1016/j.aim.2013.10.003},
 zbMATH = {6284418},
 Zbl = {1291.14043}
}

@ARTICLE{GetzlerPreprint,
       author = {{Getzler}, Ezra},
        title = "{Mixed Hodge structures of configuration spaces}",
      journal = {arXiv e-prints},
         year = 1995,
        month = oct,
          eid = {alg-geom/9510018},
        pages = {alg-geom/9510018},
          doi = {10.48550/arXiv.alg-geom/9510018},
archivePrefix = {arXiv},
       eprint = {alg-geom/9510018},
 primaryClass = {math.AG},
       adsurl = {https://ui.adsabs.harvard.edu/abs/1995alg.geom.10018G}
}

@book {Macdonald,
    AUTHOR = {Macdonald, I. G.},
     TITLE = {Symmetric functions and {H}all polynomials},
    SERIES = {Oxford Mathematical Monographs},
   EDITION = {Second},
      NOTE = {With contributions by A. Zelevinsky,
              Oxford Science Publications},
 PUBLISHER = {The Clarendon Press, Oxford University Press, New York},
      YEAR = {1995},
     PAGES = {x+475},
      ISBN = {0-19-853489-2},
   MRCLASS = {05E05 (05-02 20C30 20C33 20K01 33C80 33D80)},
  MRNUMBER = {1354144},
MRREVIEWER = {John R. Stembridge},
}

@article {BergstromMinabe2,
    AUTHOR = {Bergstr\"{o}m, Jonas and Minabe, Satoshi},
     TITLE = {On the cohomology of the {L}osev-{M}anin moduli space},
   JOURNAL = {Manuscripta Math.},
  FJOURNAL = {Manuscripta Mathematica},
    VOLUME = {144},
      YEAR = {2014},
    NUMBER = {1-2},
     PAGES = {241--252},
      ISSN = {0025-2611},
   MRCLASS = {14H10 (14M25)},
  MRNUMBER = {3193775},
MRREVIEWER = {Dawei Chen},
       DOI = {10.1007/s00229-013-0647-5},
       URL = {https://doi.org/10.1007/s00229-013-0647-5},
}

@article {BergstromMinabe1,
    AUTHOR = {Bergstr\"{o}m, Jonas and Minabe, Satoshi},
     TITLE = {On the cohomology of moduli spaces of (weighted) stable
              rational curves},
   JOURNAL = {Math. Z.},
  FJOURNAL = {Mathematische Zeitschrift},
    VOLUME = {275},
      YEAR = {2013},
    NUMBER = {3-4},
     PAGES = {1095--1108},
      ISSN = {0025-5874},
   MRCLASS = {14H10 (20C30)},
  MRNUMBER = {3127048},
MRREVIEWER = {Montserrat Teixidor i Bigas},
       DOI = {10.1007/s00209-013-1171-8},
       URL = {https://doi.org/10.1007/s00209-013-1171-8},
}

@incollection {GetzlerGenusZero,
    AUTHOR = {Getzler, E.},
     TITLE = {Operads and moduli spaces of genus {$0$} {R}iemann surfaces},
 BOOKTITLE = {The moduli space of curves ({T}exel {I}sland, 1994)},
    SERIES = {Progr. Math.},
    VOLUME = {129},
     PAGES = {199--230},
 PUBLISHER = {Birkh\"{a}user Boston, Boston, MA},
      YEAR = {1995},
   MRCLASS = {18C10 (14H10 18D99 18G10 55P99)},
  MRNUMBER = {1363058},
MRREVIEWER = {J. Stasheff},
       DOI = {10.1007/978-1-4612-4264-2\_8},
       URL = {https://doi.org/10.1007/978-1-4612-4264-2_8},
}

@article {GetzlerSemiClassical,
    AUTHOR = {Getzler, E.},
     TITLE = {The semi-classical approximation for modular operads},
   JOURNAL = {Comm. Math. Phys.},
  FJOURNAL = {Communications in Mathematical Physics},
    VOLUME = {194},
      YEAR = {1998},
    NUMBER = {2},
     PAGES = {481--492},
      ISSN = {0010-3616},
   MRCLASS = {14H10 (14H15 14H52)},
  MRNUMBER = {1627677},
MRREVIEWER = {Andreas Gathmann},
       DOI = {10.1007/s002200050365},
       URL = {https://doi.org/10.1007/s002200050365},
}

@article {GetzlerKapranov,
    AUTHOR = {Getzler, E. and Kapranov, M. M.},
     TITLE = {Modular operads},
   JOURNAL = {Compositio Math.},
  FJOURNAL = {Compositio Mathematica},
    VOLUME = {110},
      YEAR = {1998},
    NUMBER = {1},
     PAGES = {65--126},
      ISSN = {0010-437X},
   MRCLASS = {18C15 (08A02 14H10 57M50)},
  MRNUMBER = {1601666},
MRREVIEWER = {Alexandre I. Kabanov},
       DOI = {10.1023/A:1000245600345},
       URL = {https://doi.org/10.1023/A:1000245600345},
}

@article {GetzlerPandharipande,
    AUTHOR = {Getzler, Ezra and Pandharipande, Rahul},
     TITLE = {The {B}etti numbers of {$\overline{\mathscr{M}}_{0,n}(r,d)$}},
   JOURNAL = {J. Algebraic Geom.},
  FJOURNAL = {Journal of Algebraic Geometry},
    VOLUME = {15},
      YEAR = {2006},
    NUMBER = {4},
     PAGES = {709--732},
      ISSN = {1056-3911},
   MRCLASS = {14D22 (14C35 14F25 14N35)},
  MRNUMBER = {2237267},
MRREVIEWER = {Hsian-Hua Tseng},
       DOI = {10.1090/S1056-3911-06-00425-5},
       URL = {https://doi.org/10.1090/S1056-3911-06-00425-5},
}

@Article{CFGP,
 Author = {Chan, Melody and Faber, Carel and Galatius, S{\o}ren and Payne, Sam},
 Title = {The {{\(S_n\)}}-equivariant top weight Euler characteristic of {{\(\mathcal{M}_{g,n}\)}}},
 FJournal = {American Journal of Mathematics},
 Journal = {Am. J. Math.},
 ISSN = {0002-9327},
 Volume = {145},
 Number = {5},
 Pages = {1549--1585},
 Year = {2023},
 Language = {English},
 DOI = {10.1353/ajm.2023.a907705},
 zbMATH = {7771203}
}

@Article{hodgehl,
 Author = {Kannan, Siddarth and Serpente, Stefano and Yun, Claudia He},
 Title = {Equivariant {Hodge} polynomials of heavy/light moduli spaces},
 FJournal = {Forum of Mathematics, Sigma},
 Journal = {Forum Math. Sigma},
 ISSN = {2050-5094},
 Volume = {12},
 Pages = {19},
 Note = {Id/No e34},
 Year = {2024},
 Language = {English},
 DOI = {10.1017/fms.2024.20},
 zbMATH = {7819868},
 Zbl = {1544.14029}
}

@article{PW1,
 author = {Payne, Sam and Willwacher, Thomas},
 title = {Weight 2 compactly supported cohomology of moduli spaces of curves},
 fjournal = {Duke Mathematical Journal},
 journal = {Duke Math. J.},
 issn = {0012-7094},
 volume = {173},
 number = {16},
 pages = {3107--3178},
 year = {2024},
 language = {English},
 doi = {10.1215/00127094-2024-0003},
 url = {projecteuclid.org/journals/duke-mathematical-journal/volume-173/issue-16/Weight-2-compactly-supported-cohomology-of-moduli-spaces-of-curves/10.1215/00127094-2024-0003.full},
 zbMATH = {7974309}
}

@incollection {semiclassicalremark,
    AUTHOR = {Petersen, Dan},
     TITLE = {A remark on {G}etzler's semi-classical approximation},
 BOOKTITLE = {Geometry and arithmetic},
    SERIES = {EMS Ser. Congr. Rep.},
     PAGES = {309--316},
 PUBLISHER = {Eur. Math. Soc., Z\"{u}rich},
      YEAR = {2012},
   MRCLASS = {18D50 (14H10 20C30)},
  MRNUMBER = {2987667},
MRREVIEWER = {Domenico Senato},
       DOI = {10.4171/119-1/18},
       URL = {https://doi.org/10.4171/119-1/18},
}

@article{MarianOpreaPandharipande,
	title = {The moduli space of stable quotients},
	volume = {15},
	issn = {1364-0380, 1465-3060},
	url = {http://www.msp.org/gt/2011/15-3/p10.xhtml},
	doi = {10.2140/gt.2011.15.1651},
	language = {en},
	number = {3},
	urldate = {2023-08-17},
	journal = {Geometry \& Topology},
	author = {Marian, Alina and Oprea, Dragos and Pandharipande, Rahul},
	month = oct,
	year = {2011},
	pages = {1651--1706},
}

@article{Sumihiro,
author = {Hideyasu Sumihiro},
title = {{Equivariant completion}},
volume = {14},
journal = {Journal of Mathematics of Kyoto University},
number = {1},
publisher = {Duke University Press},
pages = {1 -- 28},
year = {1974},
doi = {10.1215/kjm/1250523277},
URL = {https://doi.org/10.1215/kjm/1250523277}
}

@article{graberpand,
	title = {Localization of virtual classes},
	volume = {135},
	copyright = {http://www.springer.com/tdm},
	issn = {0020-9910, 1432-1297},
	url = {http://link.springer.com/10.1007/s002220050293},
	doi = {10.1007/s002220050293},
	number = {2},
	urldate = {2024-11-21},
	journal = {Inventiones Mathematicae},
	author = {Graber, T. and Pandharipande, R.},
	month = jan,
	year = {1999},
	pages = {487--518}
}

@article{AluMarNas,
 author = {Aluffi, Paolo and Marcolli, Matilde and Nascimento, Eduardo},
 title = {Explicit formulas for the {Grothendieck} class of {{\(\overline{\mathcal{M}}_{0,n}\)}}},
 fjournal = {Moduli},
 journal = {Moduli},
 issn = {2949-7647},
 volume = {2},
 pages = {32},
 note = {Id/No e19},
 year = {2025},
 language = {English},
 doi = {10.1112/mod.2025.10012},
 zbMATH = {8137203}
}

@article{vzdualcomplex,
      title={The dual complex of $\mathcal{M}_{1,n}(\mathbb{P}^r,d)$ via the geometry of the {V}akil--{Z}inger moduli space}, 
      author={Siddarth Kannan and Terry Dekun Song},
      journal = {to appear in Adv. Math.},
      year={2026}
}

@misc{Con19,
    author = "Keith Conrad",
    title = "Dihedral {G}roups {II}",
    howpublished = "online expository notes", 
    year = "2019"}

@article{birkhoff,
 ISSN = {0003486X, 19398980},
 URL = {http://www.jstor.org/stable/1967597},
 author = {George D. Birkhoff},
 journal = {Annals of Mathematics},
 number = {1/4},
 pages = {42--46},
 publisher = {[Annals of Mathematics, Trustees of Princeton University on Behalf of the Annals of Mathematics, Mathematics Department, Princeton University]},
 title = {A Determinant Formula for the Number of Ways of Coloring a Map},
 urldate = {2024-12-03},
 volume = {14},
 year = {1912}
}

@article{whitney,
author = {Hassler Whitney},
title = {{A logical expansion in mathematics}},
volume = {38},
journal = {Bulletin of the American Mathematical Society},
number = {8},
publisher = {American Mathematical Society},
pages = {572 -- 579},
year = {1932},
}

@article{PetersenOperadGCovers,
author = {Dan Petersen},
title = {{The operad structure of admissible $G$-covers}},
volume = {7},
journal = {Algebra \& Number Theory},
number = {8},
publisher = {MSP},
pages = {1953 -- 1975},
year = {2013},
doi = {10.2140/ant.2013.7.1953},
URL = {https://doi.org/10.2140/ant.2013.7.1953}
}

@article{Robinson,
title = {Enumeration of colored graphs},
journal = {Journal of Combinatorial Theory},
volume = {4},
number = {2},
pages = {181-190},
year = {1968},
issn = {0021-9800},
doi = {https://doi.org/10.1016/S0021-9800(68)80040-7},
author = {Robert W. Robinson}
}

@book{PolyaRead,
    AUTHOR = {P\'olya, G. and Read, R. C.},
     TITLE = {Combinatorial enumeration of groups, graphs, and chemical
              compounds},
      NOTE = {P\'olya's contribution translated from the German by Dorothee
              Aeppli},
 PUBLISHER = {Springer-Verlag, New York},
      YEAR = {1987},
     PAGES = {viii+148},
      ISBN = {0-387-96413-4},
   MRCLASS = {05A15 (01A75 05C30 20B05 92A40)},
  MRNUMBER = {884155},
MRREVIEWER = {Daniel\ Turz\'ik},
       DOI = {10.1007/978-1-4612-4664-0},
       URL = {https://doi.org/10.1007/978-1-4612-4664-0},
}

@article{Moreau,
    author = {Moreau, C.},
    title = {Sur les permutations circulaires distinctes},
    journal = {Nouv. Ann. (2)},
    year = {1872},
    pages = {309--314},
    volume = {11}
}

@article{Polya,
author = {G. P{\'o}lya},
title = {{Kombinatorische Anzahlbestimmungen für Gruppen, Graphen und chemische Verbindungen}},
volume = {68},
journal = {Acta Mathematica},
number = {none},
publisher = {Institut Mittag-Leffler},
pages = {145 -- 254},
year = {1937},
doi = {10.1007/BF02546665},
URL = {https://doi.org/10.1007/BF02546665}
}

@article{VakilZinger,
author = {Ravi Vakil and Aleksey Zinger},
title = {{A desingularization of the main component of the moduli space of genus-one stable maps into $\mathbb P^n$}},
volume = {12},
journal = {Geometry \& Topology},
number = {1},
publisher = {MSP},
pages = {1 -- 95},
year = {2008},
doi = {10.2140/gt.2008.12.1},
URL = {https://doi.org/10.2140/gt.2008.12.1}
}

@article{rspw,
	title = {Moduli of stable maps in genus one and logarithmic geometry, {I}},
	volume = {23},
	issn = {1364-0380, 1465-3060},
	url = {https://msp.org/gt/2019/23-7/p03.xhtml},
	doi = {10.2140/gt.2019.23.3315},
	language = {en},
	number = {7},
	urldate = {2023-08-17},
	journal = {Geometry \& Topology},
	author = {Ranganathan, Dhruv and Santos-Parker, Keli and Wise, Jonathan},
	month = dec,
	year = {2019},
	pages = {3315--3366},
}

@article{Lee,
author = {Y.-P. Lee},
title = {{Quantum $K$-theory, I: Foundations}},
volume = {121},
journal = {Duke Mathematical Journal},
number = {3},
publisher = {Duke University Press},
pages = {389 -- 424},
year = {2004},
doi = {10.1215/S0012-7094-04-12131-1},
URL = {https://doi.org/10.1215/S0012-7094-04-12131-1}
}

@article{GiventalWDVV,
author = {Alexander Givental},
title = {{On the WDVV equation in quantum K-theory.}},
volume = {48},
journal = {Michigan Mathematical Journal},
number = {1},
publisher = {University of Michigan, Department of Mathematics},
pages = {295 -- 304},
year = {2000},
doi = {10.1307/mmj/1030132720},
URL = {https://doi.org/10.1307/mmj/1030132720}
}

@misc{GiventalII,
      title={Permutation-equivariant quantum {K}-theory {I}{I}. Fixed point localization}, 
      author={Alexander Givental},
      year={2015},
      eprint={1508.04374},
      archivePrefix={arXiv},
      primaryClass={math.AG},
      url={https://arxiv.org/abs/1508.04374}, 
}

@article{AndersonChenTseng,
    author = {Anderson, David and Chen, Linda and Tseng, Hsian-Hua},
    title = {On the Finiteness of Quantum {K}-Theory of a Homogeneous Space},
    journal = {International Mathematics Research Notices},
    volume = {2022},
    number = {2},
    pages = {1313-1349},
    year = {2020},
    month = {06},
    issn = {1073-7928},
    doi = {10.1093/imrn/rnaa108},
    url = {https://doi.org/10.1093/imrn/rnaa108},
    eprint = {https://academic.oup.com/imrn/article-pdf/2022/2/1313/42236033/rnaa108.pdf},
}

@misc{GiventalI,
      title={Permutation-equivariant quantum {K}-theory {I}. {D}efinitions. {E}lementary {K}-theory of $\overline{M}_{0, n}/{S}_n$}, 
      author={Alexander Givental},
      year={2015},
      eprint={1508.02690},
      archivePrefix={arXiv},
      primaryClass={math.AG},
      url={https://arxiv.org/abs/1508.02690}, 
}

@InCollection{KontsevichTorusActions,
 Author = {Kontsevich, Maxim},
 Title = {Enumeration of rational curves via torus actions},
 BookTitle = {The moduli space of curves. Proceedings of the conference held on Texel Island, Netherlands during the last week of April 1994},
 ISBN = {0-8176-3784-2},
 Pages = {335--368},
 Year = {1995},
 Publisher = {Basel: Birkh{\"a}user},
 Language = {English},
 zbMATH = {826032},
 Zbl = {0885.14028}
}

@Article{MustataMustata1,
 Author = {Musta{\c{t}}{\u{a}}, Andrei and Musta{\c{t}}{\u{a}}, Magdalena Anca},
 Title = {Intermediate moduli spaces of stable maps},
 FJournal = {Inventiones Mathematicae},
 Journal = {Invent. Math.},
 ISSN = {0020-9910},
 Volume = {167},
 Number = {1},
 Pages = {47--90},
 Year = {2007},
 Language = {English},
 DOI = {10.1007/s00222-006-0006-1},
 zbMATH = {5118882},
 Zbl = {1111.14018}
}

@Article{MustataMustata2,
 Author = {Musta{\c{t}}{\u{a}}, Anca M. and Musta{\c{t}}{\u{a}}, Andrei},
 Title = {The {Chow} ring of {{\({\overline M}_{0,m}({\mathbb P}^n,d)\)}}},
 FJournal = {Journal f{\"u}r die Reine und Angewandte Mathematik},
 Journal = {J. Reine Angew. Math.},
 ISSN = {0075-4102},
 Volume = {615},
 Pages = {93--119},
 Year = {2008},
 Language = {English},
 DOI = {10.1515/CRELLE.2008.011},
 zbMATH = {5267834},
 Zbl = {1139.14043}
}

@ARTICLE{Banerjee,
       author = {{Banerjee}, Oishee},
        title = "{Moduli of curves on toric varieties and their stable cohomology}",
      journal = {arXiv e-prints},
         year = 2022,
        month = oct,
          eid = {arXiv:2210.05826},
        pages = {arXiv:2210.05826},
          doi = {10.48550/arXiv.2210.05826},
archivePrefix = {arXiv},
       eprint = {2210.05826},
 primaryClass = {math.AG},
       adsurl = {https://ui.adsabs.harvard.edu/abs/2022arXiv221005826B}
}

@Article{OpreaPrTaur,
 Author = {Oprea, Dragos},
 Title = {Tautological classes on the moduli spaces of stable maps to {{\(\mathbb {P}^r\)}} via torus actions},
 FJournal = {Advances in Mathematics},
 Journal = {Adv. Math.},
 ISSN = {0001-8708},
 Volume = {207},
 Number = {2},
 Pages = {661--690},
 Year = {2006},
 Language = {English},
 DOI = {10.1016/j.aim.2006.01.002},
 zbMATH = {5077972},
 Zbl = {1117.14056}
}

@Article{Diaconu,
 Author = {Diaconu, Adrian},
 Title = {Equivariant {Euler} characteristics of {{\(\overline{\mathscr{M}}_{g,n}\)}}},
 FJournal = {Algebraic Geometry},
 Journal = {Algebr. Geom.},
 ISSN = {2313-1691},
 Volume = {7},
 Number = {5},
 Pages = {523--543},
 Year = {2020},
 Language = {English},
 DOI = {10.14231/AG-2020-018},
 zbMATH = {7262976},
 Zbl = {1452.14022}
}

@Article{NonlinearGrassmannian,
 Author = {Pandharipande, Rahul},
 Title = {The {Chow} ring of the nonlinear {Grassmannian}},
 FJournal = {Journal of Algebraic Geometry},
 Journal = {J. Algebr. Geom.},
 ISSN = {1056-3911},
 Volume = {7},
 Number = {1},
 Pages = {123--140},
 Year = {1998},
 Language = {English},
 zbMATH = {1353490},
 Zbl = {0941.14017}
}

@Article{Lambdag,
 Author = {Faber, C. and Pandharipande, R.},
 Title = {Hodge integrals, partition matrices, and the {{\(\lambda_g\)}} conjecture},
 FJournal = {Annals of Mathematics. Second Series},
 Journal = {Ann. Math. (2)},
 ISSN = {0003-486X},
 Volume = {157},
 Number = {1},
 Pages = {97--124},
 Year = {2003},
 Language = {English},
 DOI = {10.4007/annals.2003.157.97},
 zbMATH = {1933156},
 Zbl = {1058.14046}
}

@Article{SamSnowden,
 Author = {Sam, Steven V. and Snowden, Andrew},
 Title = {Gr{\"o}bner methods for representations of combinatorial categories},
 FJournal = {Journal of the American Mathematical Society},
 Journal = {J. Am. Math. Soc.},
 ISSN = {0894-0347},
 Volume = {30},
 Number = {1},
 Pages = {159--203},
 Year = {2017},
 Language = {English},
 DOI = {10.1090/jams/859},
 zbMATH = {6640528},
 Zbl = {1347.05010}
}

@Article{BattistellaCarocci,
 Author = {Battistella, Luca and Carocci, Francesca},
 Title = {A smooth compactification of the space of genus two curves in projective space: via logarithmic geometry and {Gorenstein} curves},
 FJournal = {Geometry \& Topology},
 Journal = {Geom. Topol.},
 ISSN = {1465-3060},
 Volume = {27},
 Number = {3},
 Pages = {1203--1272},
 Year = {2023},
 Language = {English},
 DOI = {10.2140/gt.2023.27.1203},
 zbMATH = {7699411},
 Zbl = {1521.14052}
}

@Article{DivisorsAndCurves,
 Author = {Kennedy-Hunt, Patrick and Nabijou, Navid and Shafi, Qaasim and Zheng, Wanlong},
 Title = {Divisors and curves on logarithmic mapping spaces},
 FJournal = {Selecta Mathematica. New Series},
 Journal = {Sel. Math., New Ser.},
 ISSN = {1022-1824},
 Volume = {30},
 Number = {4},
 Pages = {30},
 Note = {Id/No 75},
 Year = {2024},
 Language = {English},
 DOI = {10.1007/s00029-024-00956-0},
 zbMATH = {7896248}
}

@article{Devkota,
 author = {Devkota, Prabhat},
 title = {Cohomology of moduli space of multiscale differentials in genus 0},
 fjournal = {IMRN. International Mathematics Research Notices},
 journal = {Int. Math. Res. Not.},
 issn = {1073-7928},
 volume = {2026},
 number = {3},
 pages = {23},
 note = {Id/No rnag006},
 year = {2026},
 language = {English},
 doi = {10.1093/imrn/rnag006},
 zbMATH = {8159198}
}

@Article{KannanP1,
 Author = {Kannan, Siddarth},
 Title = {Moduli of relative stable maps to {{\(\mathbb{P}^1\)}}: cut-and-paste invariants},
 FJournal = {Selecta Mathematica. New Series},
 Journal = {Sel. Math., New Ser.},
 ISSN = {1022-1824},
 Volume = {29},
 Number = {4},
 Pages = {26},
 Note = {Id/No 54},
 Year = {2023},
 Language = {English},
 DOI = {10.1007/s00029-023-00857-8},
 zbMATH = {7722293},
 Zbl = {1530.14052}
}

@Article{Bae,
 Author = {Bae, Younghan},
 Title = {Tautological relations for stable maps to a target variety},
 FJournal = {Arkiv f{\"o}r Matematik},
 Journal = {Ark. Mat.},
 ISSN = {0004-2080},
 Volume = {58},
 Number = {1},
 Pages = {19--38},
 Year = {2020},
 Language = {English},
 DOI = {10.4310/ARKIV.2020.v58.n1.a2},
 zbMATH = {7201276},
 Zbl = {1436.14048}
}

@Article{OpreaFlags,
 Author = {Oprea, Dragos},
 Title = {The tautological rings of the moduli spaces of stable maps to flag varieties},
 FJournal = {Journal of Algebraic Geometry},
 Journal = {J. Algebr. Geom.},
 ISSN = {1056-3911},
 Volume = {15},
 Number = {4},
 Pages = {623--655},
 Year = {2006},
 Language = {English},
 DOI = {10.1090/S1056-3911-06-00452-8},
 zbMATH = {5141401},
 Zbl = {1195.14037}
}

@Article{HodgeIntegrals,
 Author = {Faber, C. and Pandharipande, R.},
 Title = {Hodge integrals and {Gromov}-{Witten} theory},
 FJournal = {Inventiones Mathematicae},
 Journal = {Invent. Math.},
 ISSN = {0020-9910},
 Volume = {139},
 Number = {1},
 Pages = {173--199},
 Year = {2000},
 Language = {English},
 DOI = {10.1007/s002229900028},
 zbMATH = {1388910},
 Zbl = {0960.14031}
}

@Article{FarbWolfson,
 Author = {Farb, Benson and Wolfson, Jesse},
 Title = {Topology and arithmetic of resultants. {I}.},
 FJournal = {The New York Journal of Mathematics},
 Journal = {New York J. Math.},
 ISSN = {1076-9803},
 Volume = {22},
 Pages = {801--821},
 Year = {2016},
 Language = {English},
 zbMATH = {6636844},
 Zbl = {1379.55016}
}

@article{CavalieriFulghesu,
 author = {Cavalieri, Renzo and Fulghesu, Damiano},
 title = {The integral {Chow} ring of {{\(\mathcal{M}_0(\mathbb{P}^r, d)\)}}, for {{\(d\)}} odd},
 fjournal = {Compositio Mathematica},
 journal = {Compos. Math.},
 issn = {0010-437X},
 volume = {159},
 number = {1},
 pages = {184--206},
 year = {2023},
 language = {English},
 doi = {10.1112/S0010437X22007898},
 zbMATH = {7656399},
 Zbl = {1516.14014}
}

@Article{BergstromTommasi,
 Author = {Bergstr{\"o}m, Jonas and Tommasi, Orsola},
 Title = {The rational cohomology of {{\(\overline{\mathcal{M}}_4\)}}},
 FJournal = {Mathematische Annalen},
 Journal = {Math. Ann.},
 ISSN = {0025-5831},
 Volume = {338},
 Number = {1},
 Pages = {207--239},
 Year = {2007},
 Language = {English},
 DOI = {10.1007/s00208-006-0073-z},
 zbMATH = {5162076},
 Zbl = {1126.14030}
}

@Article{HarerZagier,
 Author = {Harer, J. and Zagier, Don},
 Title = {The {Euler} characteristic of the moduli space of curves},
 FJournal = {Inventiones Mathematicae},
 Journal = {Invent. Math.},
 ISSN = {0020-9910},
 Volume = {85},
 Pages = {457--485},
 Year = {1986},
 Language = {English},
 DOI = {10.1007/BF01390325},
 URL = {https://eudml.org/doc/143377},
 zbMATH = {3997972},
 Zbl = {0616.14017}
}

@Article{BehrendOH,
 Author = {Behrend, K. and O'Halloran, A.},
 Title = {On the cohomology of stable map spaces},
 FJournal = {Inventiones Mathematicae},
 Journal = {Invent. Math.},
 ISSN = {0020-9910},
 Volume = {154},
 Number = {2},
 Pages = {385--450},
 Year = {2003},
 Language = {English},
 DOI = {10.1007/s00222-003-0308-5},
 zbMATH = {2021107},
 Zbl = {1092.14019}
}

@Article{Guest,
 Author = {Guest, Martin A.},
 Title = {The topology of the space of rational curves on a toric variety},
 FJournal = {Acta Mathematica},
 Journal = {Acta Math.},
 ISSN = {0001-5962},
 Volume = {174},
 Number = {1},
 Pages = {119--145},
 Year = {1995},
 Language = {English},
 DOI = {10.1007/BF02392803},
 zbMATH = {734731},
 Zbl = {0826.14035}
}

@Article{Segal,
 Author = {Segal, Graeme},
 Title = {The topology of spaces of rational functions},
 FJournal = {Acta Mathematica},
 Journal = {Acta Math.},
 ISSN = {0001-5962},
 Volume = {143},
 Pages = {39--72},
 Year = {1979},
 Language = {English},
 DOI = {10.1007/BF02392088},
 zbMATH = {3665822},
 Zbl = {0427.55006}
}

@Article{Hanlon,
 Author = {Hanlon, P.},
 Title = {The chromatic polynomial of an unlabeled graph},
 FJournal = {Journal of Combinatorial Theory. Series B},
 Journal = {J. Comb. Theory, Ser. B},
 ISSN = {0095-8956},
 Volume = {38},
 Pages = {226--239},
 Year = {1985},
 Language = {English},
 DOI = {10.1016/0095-8956(85)90068-1},
 zbMATH = {3904613},
 Zbl = {0567.05025}
}

@Article{Qmapwallcrossing,
 Author = {Ciocan-Fontanine, Ionu{\c{t}} and Kim, Bumsig},
 Title = {Quasimap wall-crossings and mirror symmetry},
 FJournal = {Publications Math{\'e}matiques},
 Journal = {Publ. Math., Inst. Hautes {\'E}tud. Sci.},
 ISSN = {0073-8301},
 Volume = {131},
 Pages = {201--260},
 Year = {2020},
 Language = {English},
 DOI = {10.1007/s10240-020-00114-0},
 zbMATH = {7209676},
 Zbl = {1475.14109}
}

@ARTICLE{ds,
       author = {{Song}, Terry Dekun},
        title = "{Topology of the Vakil--Zinger moduli space}",
      journal = {arXiv e-prints},
         year = 2026,
        month = feb,
          eid = {arXiv:2602.16103},
        pages = {arXiv:2602.16103},
          doi = {10.48550/arXiv.2602.16103},
archivePrefix = {arXiv},
       eprint = {2602.16103},
 primaryClass = {math.AG},
       adsurl = {https://ui.adsabs.harvard.edu/abs/2026arXiv260216103D}
}

@book{tomDieck,
url = {https://doi.org/10.1515/9783110858372},
title = {{T}ransformation {G}roups},
author = {Tammo tom Dieck},
publisher = {De Gruyter},
address = {Berlin, New York},
doi = {doi:10.1515/9783110858372},
isbn = {9783110858372},
year = {1987},
lastchecked = {2024-12-13}
}

@article{BF, title={Calculating cohomology groups of {$\overline{\mathcal{M}}_{0,n}(\mathbb{P}^r,d)$}}, volume={182}, DOI={10.1007/s10231-003-0071-7}, number={3}, journal={Annali di Matematica Pura ed Applicata}, author={Bini, Gilberto and Fontanari, Claudio}, year={2003}, month={Aug}, pages={337–344}}

@article{PandharipandeQ,
 ISSN = {00029947},
 URL = {http://www.jstor.org/stable/117856},
 author = {Rahul Pandharipande},
 journal = {Transactions of the American Mathematical Society},
 number = {4},
 pages = {1481--1505},
 publisher = {American Mathematical Society},
 title = {Intersections of {$\mathbb{Q}$}-{D}ivisors on {K}ontsevich's Moduli Space {$\overline{M}_{0,n}(\mathbb{P}^r,d)$} and Enumerative Geometry},
 urldate = {2024-12-15},
 volume = {351},
 year = {1999}
}

@article{KM,
author = {Kiem, Young-Hoon and Moon, Han-Bom},
title = {Moduli space of stable maps to projective space via {G}{I}{T}},
journal = {International Journal of Mathematics},
volume = {21},
number = {05},
pages = {639-664},
year = {2010},
doi = {10.1142/S0129167X10006264},

URL = { 
        https://doi.org/10.1142/S0129167X10006264
},
eprint = { 
        https://doi.org/10.1142/S0129167X10006264  
}
}

@article{MACDONALD1980,
title = {Polynomial functors and wreath products},
journal = {Journal of Pure and Applied Algebra},
volume = {18},
number = {2},
pages = {173-204},
year = {1980},
issn = {0022-4049},
doi = {https://doi.org/10.1016/0022-4049(80)90128-0},
url = {https://www.sciencedirect.com/science/article/pii/0022404980901280},
author = {I.G. MacDonald}
}

\end{document}